\documentclass[a4paper, 10pt]{article}
\usepackage{showlabels}
\usepackage[T1]{fontenc}			
\usepackage[english]{babel}
\usepackage{array}						
\usepackage{graphics}					
\usepackage{indentfirst} 			
\usepackage{vmargin}					
\usepackage{amsmath,amsthm,amssymb,stmaryrd,wasysym}
\usepackage{color,colortbl}
\usepackage[table]{xcolor}
\usepackage{stackrel}
\usepackage{mathrsfs}
\usepackage{fancyref}
\usepackage{graphicx}
\usepackage{enumerate}
\usepackage{mdwlist}
\usepackage{caption}
\usepackage{fancybox}
\usepackage{pifont}
\usepackage{makecell}
\usepackage{epigraph}
\usepackage{pifont}
\usepackage{dsfont}
\usepackage{relsize}
\usepackage{rotating}
\usepackage{multirow}
\newcolumntype{?}{!{\vrule width 1pt}}
\usepackage{csquotes}
\usepackage{chngcntr}
\usepackage{tikz-cd}
\usepackage{changepage}
\usepackage{todonotes}
\usepackage{mathtools}
\usepackage{stmaryrd}
\usepackage[all,cmtip]{xy}
\usepackage{amssymb}
\usepackage{amsmath}
\usepackage{mathrsfs} 
\usepackage{booktabs}
\usepackage{tkz-graph, tikz-cd}

\numberwithin{equation}{section}
\counterwithin*{equation}{section}

\input amssym.def

\usepackage{fancyhdr}
\makeindex

\def\bqn{\begin{eqnarray}}
\def\eqn{\end{eqnarray}}
\newcommand*\xbar[1]{
  \hbox{
    \vbox{
      \hrule height 0.5pt 
      \kern0.5ex         
      \hbox{
        \kern-0.25em      
        \ensuremath{#1}
        \kern-0.25em
      }
    }
  }
}

\newtheorem{Theorem}{Theorem}[section]
\newtheorem*{Theorem*}{Theorem}
\newtheorem{Corollary}[Theorem]{Corollary}
\newtheorem{Lemma}[Theorem]{Lemma}
\newtheorem{Proposition}[Theorem]{Proposition}

 { \theoremstyle{definition}
\newtheorem{Definition}[Theorem]{Definition}

\newtheorem{Example}[Theorem]{Example}

}

\usepackage{amsmath}
\usepackage{empheq}
\usepackage[most]{tcolorbox}

\newtcbox{\mymath}[1][]{    nobeforeafter, math upper, tcbox raise base,
    enhanced, colframe=blue!30!black,
       colback=blue!30!red!30!white, boxrule=1pt,   
    #1}

\usepackage[	pdftex,
			colorlinks=true,	
			pdfstartview=FitV,
			linkcolor= color3,
			citecolor= lightblue,
			urlcolor= color3,
			hyperindex=true,
			hyperfigures=false,unicode]
                      			{hyperref}

\definecolor{lightblue}{rgb}{0.22,0.45,0.70}
  \definecolor{color3}{RGB}{  6, 99, 112}

\newcommand{\p}{\partial}

\newcommand{\nn}{\nonumber}

\newcommand{\pl}{\left(}
\newcommand{\pr}{\right)}

\newcommand{\crl}{\left[}
\newcommand{\crr}{\right]}

\newcommand{\dl}{\left \langle}
\newcommand{\dr}{\right \rangle}

\newcommand{\w}{\wedge}
\newcommand{\half}{\frac{1}{2}}

\newcommand{\Lag}{\mathcal{L}}

\newcommand{\alg}{\mathfrak g}
\newcommand{\alh}{\mathfrak h}

\newcommand\fonc[1]{\mathscr C^\infty\pl#1\pr}

\newcommand{\vf}{\field{T\M}}

\newcommand{\M}{\mathscr{M}}

\newcommand{\Om}{\Omega}
\newcommand{\om}{\omega}

\newcommand{\un}{^{-1}}
\newcommand{\ie}{\textit{i.e. }}
\newcommand{\cf}{\textit{cf. }}
\newcommand{\eg}{\textit{e.g. }}

\newcommand{\etc}{\textit{etc.}\, }
\newcommand{\br}[2]{\crl #1,#2 \crr}

\newcommand{\mD}{\mathcal D}

\newcommand\pset[1]{\left \lbrace #1\right \rbrace}

\newcommand{\Ker}{\text{\rm Ker }}

\newcommand{\End}{\mathsf{End}}

\newcommand\field[1]{\Gamma\pl #1\pr}

\newcommand\Prop[1]{Proposition \ref{#1}}
\newcommand\Defi[1]{Definition \ref{#1}}

\newcommand\Hom{\text{\rm Hom}}

\newcommand{\eps}{{\vcenter{\hbox{{$[[\epsilon]]$}}}}}

\newcommand{\brdot}{\br{\cdot}{\cdot}}
\newcommand{\pbdot}{\left\{{\cdot},{\cdot}\right\}}
\newcommand{\pb}[2]{\big\{{#1},{#2}\big\}}

\newcommand{\poly}{\text{\rm poly}}

\newcommand{\Id}{\mathrm{id}}
\newcommand{\mP}{\mathscr{P}}

\newcommand{\iso}{\overset{\sim}{\longrightarrow}}

\newcommand{\V}{\mathcal{V}}

\newcommand{\Lieinf}{$\mathsf{Lie}_\infty$}

\newcommand{\mAss}{\mathsf{Ass}}

\newcommand{\mLie}{\mathsf{Lie}}

\newcommand{\mCom}{\mathsf{Com}}

\newcommand{\mGer}{\mathsf{Ger}}

\newcommand{\Gra}{\mathsf{Gra}}
\newcommand{\dGra}{\mathsf{dGra}}

\newcommand{\gra}{\mathsf{gra}}
\newcommand{\fGC}{\mathsf{fGC}}
\newcommand{\dfGC}{\mathsf{dfGC}}
\newcommand{\dfGCconn}{\mathsf{dfGC}^{\text{\rm con}}}
\newcommand{\GC}{\mathsf{GC}}
\newcommand{\GCor}{\mathsf{GC}^\text{\rm or}}

\newcommand{\corps}{{\mathbb K}}
\newcommand{\fE}{\field{E}}
\newcommand{\foncm}{{\mathscr C}^\infty(\M)}
\newcommand{\ala}{{\it\`a la} }

\newcommand{\dldot}{{\dl\cdot,\cdot\dr_E}}
\newcommand{\dlr}[2]{\dl#1,#2\dr}

\newcommand{\Tpoly}{{\mathcal T_\text{\rm poly}}}
\newcommand{\Dpoly}{{\mathcal D_\text{\rm poly}}}
\newcommand{\Tpolyd}[1]{{\mathcal T^{#1}_\text{\rm poly}}}
\newcommand{\Dpolyd}[1]{{\mathcal D^{#1}_\text{\rm poly}}}

\newcommand{\cH}{{\mathscr H}}
\newcommand{\Q}{{\mathsf Q}}
\newcommand{\Or}{{\mathsf{O}\vec{\mathsf r}}}
\newcommand{\tOr}{s^*{\mathsf{O}\vec{\mathsf r}}}
\newcommand{\sRep}{s^*{\mathsf{Rep}}}

\newcommand{\sgn}{{\text{\rm sgn}}}
\newcommand\foncg[1]{{\mathscr C^{\infty|#1}(\mathcal V)}}
\newcommand{\mS}{{\mathbb S}}
\newcommand{\dime}{{d}}
\newcommand{\NR}{{\mathsf {NR}}}
\newcommand{\m}{{\mathfrak {m}}}
\newcommand{\Ham}{{\rm\mathsf{Ham}}}
\newcommand{\FHam}{{\rm\mathsf{FHam}}}
\newcommand{\MC}{{\rm\mathsf{MC}}}
\newcommand{\CE}{{\rm\mathsf{CE}}}
\newcommand{\NP}{{\rm\mathsf{NP}}}
\newcommand{\deltaS}{\delta_{\rm\mathsf{S}}}
\newcommand{\NPQ}{{\rm\mathsf{NPQ}}}

\newcommand{\SQI}{{\rm\mathsf{SQI}}}
\newcommand{\Tpolyn}{\Tpolyd{(n)}}
\newcommand{\fGCconn}{\fGC^\text{\rm con}}
\newcommand{\grt}{\mathfrak{grt}}
\newcommand{\GRT}{\mathsf{GRT}}
\newcommand{\mGerd}{\mGer_{\hspace{-0.3mm} d}} 
\newcommand{\DAss}{\mathsf{DAss}}

\newcommand{\Rep}{{\rm\mathsf{Rep}}}
\newcommand{\dRep}{{\rm\mathsf{dRep}}}
\newcommand{\npu}{{d}}
\newcommand{\Gamgraph}{{\gamma}}
\newcommand{\Gamgraphgroup}{{\Gamma}}
\newcommand*{\longhookrightarrow}{\ensuremath{\lhook\joinrel\relbar\joinrel\rightarrow}}
\newcommand*{\longtwoheadrightarrow}{\ensuremath{\joinrel\relbar\joinrel\twoheadrightarrow}}
\tikzstyle{n}=[circle, draw, fill, minimum size=6, inner sep=0]
\tikzstyle{ngr}=[circle, draw, fill, gray, minimum size=6, inner sep=0]
\tikzstyle{int}=[draw, circle, fill,  minimum size=4, inner sep=1]
\tikzstyle{root}=[circle, draw, fill, minimum size=0, inner sep=0]
\tikzstyle{vertex}=[circle, draw, fill, minimum size=0, inner sep=0.7]
\tikzstyle{leaf}=[circle, draw, fill, minimum size=0, inner sep=0]
\tikzstyle{intgr}=[draw, circle, fill, gray,  minimum size=3, inner sep=1]
\tikzstyle{ext}=[draw, circle,  minimum size=1, inner sep=1]
\tikzstyle{exttiny}=[draw, circle,  minimum size=1, inner sep=0.65]
\tikzset{lab/.style={draw, circle, minimum size=5, inner sep=1.1}, 
n/.style={draw, circle, fill, minimum size=5, inner sep=1}}
\tikzstyle{fleche}=[->,>= latex?,thick]

\newcommand{\itemnum}{\hfill\refstepcounter{equation}\textup{(\theequation)}}

\newcommand{\Gambr}{{\Gamma_{\vcenter{\hbox{\begin{tikzpicture}[scale=0.2, >=stealth']
\tikzstyle{w}=[circle, draw, minimum size=4, inner sep=1]
\tikzstyle{b}=[circle, draw, fill, minimum size=4, inner sep=1]
\node [b] (b1) at (0,0) {};
\node [b] (b2) at (1.3,0) {};
\draw (b1) to[line width = 0.8]  (b2);
\end{tikzpicture}}}}}}

\newcommand{\deltabr}{\delta}

\newcommand{\Gammult}{\Gamma_{\vcenter{\hbox{\begin{tikzpicture}[scale=0.2, >=stealth']
\tikzstyle{w}=[circle, draw, minimum size=4, inner sep=1]
\tikzstyle{b}=[circle, draw, fill, minimum size=4, inner sep=1]
\node [b] (b1) at (0,0) {};
\node [b] (b2) at (1.3,0) {};
\end{tikzpicture}}}}}

\newcommand{\SectionRef}[1]{\unskip~[\hyperref[#1]{\S\thinspace\ref{#1}}]}

\begin{document}
\tcbset{highlight math style={colback=blue!30!red!30!white}}
\thispagestyle{empty}

\vspace{1cm}

 \begin{centering}

{\large {\bfseries 
M.~Kontsevich's graph complexes and universal structures\\ on graded symplectic manifolds
}
}
\\\vspace{0.5cm}
Kevin Morand
\\
\vspace{0.5cm}
Department of Physics, Sogang University, Seoul  04107, South Korea\\
Center for Quantum Spacetime, Sogang University, Seoul  04107, South Korea\\
\vspace{0.5cm}
{\tt morand@sogang.ac.kr}

\vspace{0.5cm}

\end{centering}
\begin{abstract}
\medskip

\centering\begin{minipage}{\dimexpr\paperwidth-5.85cm}
\noindent
In the formulation of his celebrated {\it Formality conjecture}, M. Kontsevich introduced a universal version of the deformation theory for the Schouten algebra of polyvector fields on affine manifolds. This universal deformation complex takes the form of a differential graded Lie algebra of graphs, denoted $\mathsf{fGC}_2$, together with an injective morphism towards the Chevalley--Eilenberg complex associated with the Schouten algebra. The latter morphism is given by explicit local formulas making implicit use of the supergeometric interpretation of the Schouten algebra as the algebra of functions on a graded symplectic manifold of degree $1$. The ambition of the present work is to generalise Kontsevich's construction to graded symplectic manifolds of arbitrary degree $n\geq1$. The corresponding graph model is given by the full Kontsevich graph complex $\mathsf{fGC}_d$ where $d=n+1$ stands for the dimension of the associated AKSZ type $\sigma$-model.  This generalisation is instrumental to classify universal structures on graded symplectic manifolds. In particular, the zeroth cohomology of the full graph complex $\mathsf{fGC}_{d}$ is shown to act via $\mathsf{Lie}_\infty$-automorphisms on the algebra of functions on graded symplectic manifolds of degree $n$. This generalises the known action of the Grothendieck--Teichm\"{u}ller algebra $\mathfrak{grt}_1\simeq H^0(\mathsf{fGC}_2)$ on the space of polyvector fields. This extended action can in turn be used to generate new universal deformations of Hamiltonian functions, generalising Kontsevich flows on the space of Poisson manifolds to differential graded manifolds of higher degrees. As an application of the general formalism, universal deformations of Courant algebroids via trivalent graphs are presented. 
\end{minipage}
\end{abstract}
\pagenumbering{gobble}
\tableofcontents

\newpage
\setcounter{secnumdepth}{0}
\section{Introduction}
\label{sectionintro}
\setcounter{secnumdepth}{3}

\pagenumbering{arabic}

In a seminal 97' preprint \cite{Kontsevich:1997vb}, M. Kontsevich proved his celebrated formality theorem by constructing an explicit \Lieinf quasi-isomorphism 
\begin{eqnarray}
\label{QI}
\mathcal U_\Phi:\Tpoly\iso \Dpoly
\end{eqnarray}
between $\Tpoly$, the  graded Lie algebra  of polyvector fields on the affine space $\mathbb R^m$,  and $\Dpoly$, the Hochschild differential graded Lie algebra (dg Lie algebra) of multidifferential operators on $\mathbb R^m$, and such that the first Taylor coefficient coincides with the Hochschild--Kostant--Rosenberg (HKR) quasi-isomorphism of complexes\footnote{The subscript $\Phi$ in \eqref{QI} will be hereafter interpreted as denoting a Drinfel'd associator. }. An important corollary of the formality theorem is that it provides an explicit bijective map\footnote{\label{footDGnew}The proof that a \Lieinf quasi-isomorphism between two dg Lie algebras induces a bijection between the associated Deligne groupoids \cite{Goldman1988} can be found in \cite{Kontsevich:1997vb,CattaneoKellerTorossianBrugieres} for the nilpotent case and in \cite{Yekutieli2011} for the pro-nilpotent case. }:
\begin{eqnarray}
\label{bijmap}
\hat{\mathcal U}_\Phi:\mathsf{FPoiss}\iso\mathsf{Star}
\end{eqnarray}
between the set $\mathsf{FPoiss}$ of (equivalence classes of) formal Poisson structures on $\mathbb R^m$ and the set $\mathsf{Star}$ of (equivalence classes of) formal associative deformations of the algebra of functions on $\mathbb R^m$ (also called {\it star products}). The bijection \eqref{bijmap} straightforwardly\footnote{\label{footcomp}Via composition of the bijective map \eqref{bijmap} with the canonical ``formalisation map'' $\mathsf{Poiss}\to\mathsf{FPoiss}:\pi\mapsto \epsilon\,  \pi$ where $\epsilon$ is a formal parameter.} induces a quantization map $\mathsf{Poiss}\to\mathsf{Star}$ assigning to any Poisson bivector $\pi\in\mathsf{Poiss}$ on $\mathbb R^m$ an equivalence class of star products $[*]\in\mathsf{Star}$ quantizing $\pi$.

An important characteristic of Kontsevich's formality morphism is that it is given by {\it universal} formulas \ie formulas applying without distinction to all affine spaces of all finite dimensions and which are defined ``graphically'' via grafting of existing structures on $\Tpoly$ without resorting to additional data. Such formality morphisms were called {\it stable} in \cite{Dolgushev2011}. Informally, these are \Lieinf quasi-isomorphisms whose Taylor coefficients can be written as a sum over Kontsevich's {\it admissible graphs} \cite{Kontsevich:1997vb} where the coefficient in front of each graph is given by a weight function, \cf \eg \cite{Willwacher2015a}. The master equation ensuring that the Taylor maps assemble to a \Lieinf-morphism thus boils down to a series of identities on the weights. Although these equations are algebraic, the only known explicit solutions make use of transcendental methods\footnote{See \cite{Dolgushev2013,Dolgushev2016} for a recursive construction of formality morphisms over rationals. } involving integrals over (compactifications of) configuration spaces of points. 

Kontsevich's formality theorem indisputably constitutes the most remarkable result in the field of deformation quantization, providing a complete solution to the quantization problem formulated in   \cite{Berezin1975,BayenFlatoFronsdalEtAl1978}. However, the transcendental methods involved in the construction are generically difficult to handle thus calling for more algebraic tools allowing to address issues in formality theory and deformation quantization. Such algebraic methods have in fact been introduced by M. Kontsevich  prior to \cite{Kontsevich:1997vb} in the formulation of his {\it Formality conjecture} around 93'-94' \cite{Kontsevich1997} (\cf also \cite{Voronov1997}). More precisely, M. Kontsevich defined in \cite{Kontsevich1997}  a universal version of the deformation theory for formality morphisms. Recall that, on very general grounds, any dg Lie algebra $\alg$ is quasi-isomorphic (as a \Lieinf-algebra) to its cohomology $H(\alg)$ endowed with a certain \Lieinf-structure obtained from the dg Lie algebra structure on $\alg$ via the homotopy transfer theorem. This allows in particular to address formality questions by studying the space of \Lieinf-structures on $H(\alg)$. Going back to the case at hand, the relevant deformation theory is therefore controlled by the Chevalley--Eilenberg dg Lie algebra $\CE(\Tpoly)$ associated with the Schouten algebra of polyvector fields. In \cite{Kontsevich1997}, M. Kontsevich introduced a universal version of $\CE(\Tpoly)$ in the guise of a dg Lie algebra of graphs, denoted $\fGC_2$, together with an injective morphism 
\begin{eqnarray}
\label{morphismfGCCE}
\fGC_2\hookrightarrow \CE(\Tpoly)
\end{eqnarray}
given by local formulas. The morphism \eqref{morphismfGCCE} allows to reformulate questions regarding formality morphisms on affine manifolds (in the stable setting) into purely algebraic questions on the cohomology of the graph complex $\fGC_2$. In particular, obstructions to the existence of a stable formality morphism\footnote{Or equivalently non-trivial stable deformations of the Schouten graded Lie algebra as a \Lieinf-algebra. } live in $H^1(\fGC_2)$ while such morphisms can be shown to be classified by $H^0(\fGC_2)$. More precisely, it was shown by V. A. Dolgushev in \cite{Dolgushev2011} that the exponentiation of the (pro-nilpotent) graded Lie algebra $H^0(\fGC_2)$ acts regularly on the space $\SQI$ of (homotopy classes of) {stable} \Lieinf quasi-isomorphisms of the form \eqref{QI} so that $\SQI$ is a {\it torsor} (or principal homogeneous space) for the pro-unipotent group $\exp\big(H^0(\fGC_{2})\big)$. Furthermore, T. Willwacher constructed in \cite{Willwacher2015} an explicit isomorphism of Lie algebras $H^{0}(\GC_2)\simeq \mathfrak{grt}_1$ where $\mathfrak{grt}_1$ stands for the Grothendieck--Teichm\"uller Lie algebra. Combining these two results leads to a full characterisation of the set $\SQI$ of stable formality maps as a $\GRT_1$-torsor\footnote{This fact was conjectured by M.~Kontsevich in \cite{KontsevichLett.Math.Phys.48:35-721999} based on the relations between the transcendental formulas involved in his formality morphism and the theory of mixed Tate motives. The Grothendieck--Teichm\"uller group and Drinfel'd associators also appear in D. Tamarkin's approach to formality \cite{Tamarkin2007,Hinich2000} via either the use of the Etingof--Kazdhan quantization of Lie bialgebras or the formality of little disks operad, \cf Section \ref{section:Poisson manifolds} for details. } where $\GRT_1$ stands for the (pro-unipotent\footnote{There are different versions of the Grothendieck--Teichm\"uller group, the most important ones being a profinite version $\widehat{\mathsf{GT}}$, a pro-$l$ version $\mathsf{GT}_l$ and a pro-unipotent version $\mathsf{GT}$. The latter is isomorphic to a graded version of the group, denoted $\mathsf{GRT}$. We will only be concerned with the exponentiation $\GRT_1=\exp(\grt_1)$ related to $\mathsf{GRT}$ via $\mathsf{GRT}=\corps^{\times}\ltimes\GRT_1$ where the action of the multiplicative group is via rescaling, \cf \cite{WillwachernotesGRT} for details.  }) Grothendieck--Teichm\"uller group $\GRT_1=\exp(\grt_1)$. The Grothendieck--Teichm\"uller group was introduced by V. Drinfel'd\footnote{Inspired by A. Grothendieck's suggestion in his {\it Esquisse d'un Programme} \cite{Grothendieck} of studying the combinatorial properties of $\text{Gal}(\bar{\mathbb Q}/\mathbb Q)$ via its natural action on the {\it tour} of Teichm\"uller groupoids.} in \cite{Drinfeld1991} in virtue of its relation to the absolute Galois group $\text{Gal}(\bar{\mathbb Q}/\mathbb Q)$ and the theory of quasi-Hopf algebras. Since then, the Grothendieck--Teichm\"uller group (together with the $\GRT_1$-torsor  of Drinfel'd associators) have found a number of applications in various areas of mathematics including the Kashiwara--Vergne conjecture in Lie theory \cite{Alekseev2005,Alekseev2008,Alekseev2010a,Schneps2012}, quantization of Lie bialgebras \cite{Etingof1995,Tamarkin2007b}, the study of multiple zeta values \cite{Le1996, Brown2011,Furusho2008}, rational homotopy automorphisms of the $\textbf{E}_2$-operad \cite{FresseI,Willwacher2015}, \etc

The action of the Grothendieck--Teichm\"uller group on formality morphisms can be traced back to an action of $\GRT_1\simeq\exp\big(H^0(\fGC_{2})\big)$ on $\Tpoly$ via \Lieinf-automorphisms\footnote{We refer to \cite{Willwacher2015} (see also \cite{Merkulov2008}) for the affine space case, \cite{Jost2012}  for a globalisation to any smooth manifolds and \cite{Dolgushev2012b} for a generalisation to the sheaf of polyvector fields on any smooth algebraic variety.}. Explicitly, to any graph cocycle $\Gamgraph\in H^0(\fGC_{2})$ one associates a (homotopy class of) \Lieinf-automorphisms $\mathcal U^\Gamgraphgroup:\Tpoly\iso \Tpoly$ where $\Gamgraphgroup:=\exp(\Gamgraph)\in\exp\big(H^0(\fGC_{2})\big)$. Composition with \eqref{QI} leads to a new (homotopy class of) formality morphisms $\mathcal U_{\Phi\cdot \Gamgraphgroup}:=\mathcal U_\Phi\circ\mathcal U^\Gamgraphgroup:\Tpoly\iso \Dpoly$. Furthermore, the bijection between Deligne groupoids derived from $\mathcal U^\Gamgraphgroup$ (\cf footnote \ref{footDGnew}) induces a stable deformation map $\hat{\mathcal U}^\Gamgraphgroup:\mathsf{FPoiss}\iso\mathsf{FPoiss}$. In particular, the latter can be used to map Poisson bivectors $\pi$ (\cf footnote \ref{footcomp}) towards stable formal Poisson structures deforming $\pi$. At first order, such stable deformations can be interpreted as maps from cocycles in  $H^{0}(\GC_2)$ to stable flows on the space of Poisson bivectors. The first example\footnote{Various examples of $d=2$ flows on the space of Poisson bivectors have recently been systematically investigated in a series of works by A.~V.~{Kiselev} and collaborators, see \cite{Bouisaghouane2016a,Bouisaghouane2016,Buring2017a,Buring2017b,Buring2017c,Buring2018,Kiselev2019}. } of such flows is the so-called {\it tetrahedral} flow introduced by M. Kontsevich in \cite[Section 5.3]{Kontsevich1997}, \cf Section \ref{section:Poisson manifolds}. 

Remarkably, Kontsevich's solution to the quantization problem of \cite{Berezin1975,BayenFlatoFronsdalEtAl1978} is inspired by ideas coming from string theory. Explicitly, Kontsevich's quantization formula can be interpreted \cite{Kontsevich:1997vb,Cattaneo2000} as the Feynman diagram expansion of a 2-dimensional topological field theory -- the {\it Poisson $\sigma$-model} -- introduced in \cite{Ikeda1993a,Ikeda1993,Schaller1994}. As noted in \cite{CattaneoFelder2001a}, the quantization of the Poisson $\sigma$-model can be best interpreted within the AKSZ formalism \cite{AlexandrovKontsevichSchwarzEtAl1997}. The latter deals with theories living on the space of maps between a {\it source} manifold of dimension $d$ and a {\it target} manifold classically endowed with a structure of differential graded symplectic manifold\footnote{Also referred to as a $\NPQ$-manifold of degree $n$, \cf Section \ref{section:NPQmanifolds} below. } of degree $n$ and such that $d=n+1$. The first and simplest example of such construction is provided by the Poisson $\sigma$-model where the source is of dimension $d=2$ and the target is the (shifted) cotangent bundle of a (finite dimensional) Poisson manifold. More generally, we will refer to the geometrical structure necessary to define a AKSZ $\sigma$-model in dimension $d$ as a {\it symplectic Lie $n$-algebroid}, with $d=n+1$. While symplectic Lie $1$-algebroids identify with Poisson manifolds, symplectic Lie $2$-algebroids correspond to {\it Courant algebroids}. The latter first appeared implicitly in the study of integrable Dirac structures \cite{Dorfman1987,Courant1986,Courant1990} before their precise geometric structure was abstracted and explicitly stated by the authors of \cite{Liu1997} in their study of double of Lie bialgebroids.  Courant algebroids play also a central r\^ole in the context of generalised complex geometry, see \cite{Hitchin2002,Gualtieri2004}. 
Their graded geometrical interpretation was put forward by D.~Roytenberg in \cite{Roytenberg2007,Roytenberged.Contemp.Math.Vol.315Amer.Math.Soc.ProvidenceRI2002} and the corresponding {\it Courant $\sigma$-model} was constructed in \cite{Ikeda2003,Roytenberg2007b}. Higher examples of symplectic Lie $n$-algebroids can be found \eg in \cite{Ikeda2010,Liu2016a,Grutzmann2010}.

An interesting open problem that arises from what precedes concerns the possibility of generalising the interplay between deformation quantization results (on the algebraic side) and quantization of AKSZ-type of models (on the field theoretic side). Motivated by this problem, the ambition of the present paper is to generalise some of the algebraic methods introduced by M. Kontsevich in \cite{Kontsevich1997} for Poisson manifolds to the case of higher symplectic Lie $n$-algebroids. Our main tool in this endeavour is given by a stable version of the Chevalley--Eilenberg dg Lie algebra associated with the deformation complex of symplectic Lie $n$-algebroids for arbitrary values of $n\geq1$. Explicitly, this graph model takes the form of an injective morphism of dg Lie algebras:
\begin{eqnarray}
\label{morphismfGCCEn}
\fGC_{d}\hookrightarrow \CE(\Tpolyn)
\end{eqnarray}
where $d=n+1$, thus generalising \eqref{morphismfGCCE} to any $d\geq2$. Here, $\fGC_{d}$ stands for the generalisation of Kontsevich's graph complex to arbitrary dimension $d$ (\cf \eg \cite{Willwacher2015}) and the dg Lie algebra $\Tpolyn$ -- referred to hereafter as the {\it $n$-Schouten algebra} -- controls the deformation theory of symplectic Lie $n$-algebroids. The morphism \eqref{morphismfGCCEn} will allow us to take advantage of the available results regarding the cohomology of $\fGC_{d}$ in order to provide a classification of the stable structures on graded symplectic manifolds of arbitrary (positive) degree. In particular, we propose a classification of \Lieinf-algebra structures deforming the $n$-Schouten algebra in a non-trivial way as well as of  \Lieinf-automorphisms of the $n$-Schouten algebra $\Tpolyn$. The latter yield in particular new stable flows on the space of symplectic Lie $n$-algebroids. 

The present paper will focus on universal structures in the {\it stable} setting [see \Defi{stability} below] \ie we consider cochains of the Chevalley--Eilenberg algebra obtained from graphs belonging to the Kontsevich graph complex of undirected graphs $\fGC_d$ (or its directed analogue $\dfGC_d$). A direct consequence of this choice is that the only incarnation of the Grothendieck--Teichm\"uller Lie algebra as a universal structure occurs in dimension $d=2$ where we recover the above mentioned action of $\GRT_1$ on $\Tpoly$ via \Lieinf-automorphisms. In higher dimensions, the universal structures are insensitive to $\grt_1$ and are in fact classified by loop cocycles. In order to obtain universal structures induced from the Grothendieck--Teichm\"uller Lie algebra in dimensions $d>2$, one needs to move away from the stable setting to enter the {\it (multi)-oriented} regime. [We refer to \cite{Morand2021} for a discussion of universal structures induced by (multi)-oriented graphs \cite{Willwacher2015c,Zivkovic2017a,Zivkovic2017,Merkulov2017,Merkulov2019}  allowing in particular to provide incarnations of the Grothendieck--Teichm\"uller algebra into the deformation theory of (quasi)-Lie bialgebroids.] 
\\ \\
\noindent {\bf Summary and main results}.

After displaying our conventions and notations in Section \ref{Section:Conventions and notations}, we dedicate Sections \ref{Section:Symplectic Lie $n$-algebroids}  and  \ref{section:Graph complexes} to a review -- aimed at non-experts -- of the principal tools and notions involved in the rest of the paper. In Section \ref{Section:Symplectic Lie $n$-algebroids}, we recall the basic concepts of graded geometry, detail the hierarchy of structures endowing graded manifolds (namely $\mathsf{N}$, $\NP$ and $\NPQ$-manifolds) and discuss their associated (non-graded) geometric counterparts. In Section \ref{section:Graph complexes}, we depart from the geometric to the algebraic realm and review the construction of the Kontsevich's full graph complex $\fGC_d$ generalising $\fGC_2$ to arbitrary dimension $d$. The differential graded Lie algebra structure on $\fGC_d$ is best introduced as a convolution Lie algebra from the graph operad $\Gra_d$ whose construction we review. We conclude the section by recalling some known facts regarding the (even and odd) cohomology of $\fGC_d$. 

Building on the last two sections, we introduce our main results in Section \ref{section:Stable structures on graded manifolds}. We start by displaying a tower of representations $\Gra_d\hookrightarrow \End_{\fonc{\V}}$ for all $d>0$ where $\V$ stands for an arbitrary $\NP$-manifold of degree $n$, such that $d=n+1$. This tower of morphism of operads will in turn induce a tower of injective morphisms of dg Lie algebras $\fGC_d\hookrightarrow\CE(\Tpolyn)$ thus providing a stable version of the Chevalley--Eilenberg complex for the $n$-Schouten algebra\footnote{Or equivalently for the graded Poisson algebra of functions $\fonc{\V}$, being isomorphic to $\Tpolyn$ through degree suspension.} $\Tpolyn$. Using this stable model, we show in particular [Corollary \ref{thmcoho}] that the pro-unipotent group $\exp\big(H^0(\fGC_{d})\big)$ acts via \Lieinf-automorphisms on $\Tpolyn$. More generally, stable structures on graded symplectic manifolds are classified in \Prop{propapp}. We conclude the section by discussing Hamiltonian deformations and their linearisation, referred to as Hamiltonian flows. In particular, we present a canonical map from the zeroth cohomology $H^{0}(\fGC_{d})$ to stable Hamiltonian deformations [\Prop{propunivhamdef}] and flows [Corollary \ref{propmapcoHd=0}]  on the space of Hamiltonian functions thus generalising Kontsevich's construction from Poisson bivectors to higher symplectic Lie $n$-algebroids. Furthermore, we describe a novel class of Hamiltonian deformations generated by {\it Weyl factors} induced by elements in $H^{-d}(\fGC_{d})$ [Corollary \ref{cor:WeylHamiltonian}]. 

Finally, Section \ref{section:Applications} is devoted to illustrate some of the machinery developed in Section \ref{section:Stable structures on graded manifolds} to the case of $\NPQ$-manifolds of degrees 1 and 2, respectively. After reviewing some known applications in the case $n=1$ (corresponding to Poisson manifolds), we turn to the case $n=2$ and present new results concerning deformations of Courant algebroids. In particular, we obtain an explicit expression for the unique deformation map for Courant algebroids induced by a loop graph and display a large class of Weyl deformations induced by trivalent graphs (modulo IHX relations). We conclude by a discussion regarding the implications of our results to the deformation quantization problem for Courant algebroids.

\section{Conventions and notations}
\label{Section:Conventions and notations}
\paragraph{Suspension.}

We will work over a ground field $\corps$ of characteristic zero. Let $V=\bigoplus_{k\in\mathbb Z}V^k$ be a graded vector space over $\corps$. The suspended graded vector space $V[k]$ is defined as $V[k]^n=V^{k+n}$ so that the suspension map $s:V[k]\to V$ is of degree $k$.

\paragraph{Invariants and coinvariants.}
Let $G$ be a group and denote $\corps\dl G\dr$ the associated group ring over $\corps$. A (right) representation of $G$ is a (right) module $M$ over the group ring  $\corps\dl G\dr$. 
Letting $M$ be a right $\corps\dl G\dr$-module, we define the two following spaces:
\begin{itemize}
\item \textbf{Invariants}: \hspace{0.43cm}$M^G:=\pset{m\in M\, |\, m\cdot g=m \text{ for all } g\in G}$ 
\item \textbf{Coinvariants}: $M_G:=M/\pset{m\cdot g-m\, |\,  g\in G \text{ and }m\in M}$
\end{itemize}
Note that while the space of invariants is a subspace of $M$, the space of coinvariants (or space of orbits) is defined as a quotient of $M$ by the group action. In other words, there are natural maps $M^G\overset{i}{\longhookrightarrow} M\overset{p}\longtwoheadrightarrow M_G$ where $i$ is injective and $p$ surjective.
If $M$ is a right $\corps\dl G\dr$-module and $N$ a left $\corps\dl G\dr$-module, then $M\otimes N$ is a right $\corps\dl G\dr$-module under the \textbf{diagonal right action} $(M\otimes N)\times G\to M\otimes N:(a,b) \times g\mapsto (a\cdot g,g\un\cdot b)$.

The associated space of coinvariants is then denoted $M\otimes_G N$. Letting $M,N$ be two right $\corps\dl G\dr$-modules, a linear map $f:M\to N$ will be said $G$-\textbf{equivariant} if it is a morphism in the category of  $\corps\dl G\dr$-modules \ie if $f(x\cdot g)=f(x)\cdot g$ for all $x\in M$ and $g\in G$. The space of $G$-equivariant maps will be denoted $\Hom_G(M,N)$.

\paragraph{Symmetric group $\mathbb S_N$.}

The symmetric group $\mathbb S_N$ is defined as the group of automorphisms of the set $\pset{1,2,\ldots,N}$. An element $\sigma\in\mathbb S_N$ is called a \textbf{permutation} and is defined by its image  $\pset{\sigma(1),\sigma(2),\ldots,\sigma(N)}$. The composition $\sigma\cdot\tau$ of two permutations $\sigma,\tau\in\mathbb S_N$ is given by $\pset{1,2,\ldots,N}\overset{\tau}{\mapsto}\pset{\tau(1),\tau(2),\ldots,\tau(N)}\overset{\sigma}{\mapsto}\pset{\sigma\big(\tau(1)\big),\sigma\big(\tau(2)\big),\ldots,\sigma\big(\tau(N)\big)}$.  In the following, we will often represent a permutation $\sigma$ by the $2\times N$ matrix
\[
\sigma=\begin{pmatrix}
1 &2&\cdots&N\\
\sigma(1) &\sigma(2)&\cdots&\sigma(N)
\end{pmatrix}.
\]
We define right and left actions of the symmetric group $\mathbb S_N$ on $V^{\otimes N}$:
\[ \begin{array}{lll}
{V^{\otimes N}\times \mathbb S_N\to V^{\otimes N}} &&
{(v_1,\ldots,v_N)\cdot \sigma=(v_{\sigma(1)},\ldots, v_{\sigma(N)})} \\
{\mathbb S_N \times V^{\otimes N}\to V^{\otimes N}} &&
{\sigma\cdot(v_1,\ldots,v_N)=(v_{\sigma\un(1)},\ldots, v_{\sigma\un(N)})}
\end{array} \]
\begin{Example}
\label{exasym}
Let $\sigma,\tau\in \mS_3$ be defined as follows
\[
{\sigma:=\begin{pmatrix}
1 &2&3\\
1 &3&2
\end{pmatrix}} \hspace{1.5em}
{\tau:=\begin{pmatrix}
1 &2&3\\
3 &1&2
\end{pmatrix}} \hspace{1.5em}
{\sigma\un=\begin{pmatrix}
1 &2&3\\
1 &3&2
\end{pmatrix}} \hspace{1.5em}
{\tau\un=\begin{pmatrix}
1 &2&3\\
2 &3&1
\end{pmatrix}}
\]
admitting inverses as shown. We compute the following compositions:
\[
\sigma\cdot \tau=\begin{pmatrix}
1 &2&3\\
2 &1&3
\end{pmatrix}
\quad,\quad 
\tau\cdot \sigma=\begin{pmatrix}
1 &2&3\\
3 &2&1
\end{pmatrix}.
\]
Now, denoting $v:=(v_1,v_2,v_3)$, one can check that:
\[
(v\cdot \sigma)\cdot \tau=(v_1,v_3,v_2)\cdot\tau=(v_2,v_1,v_3)=v\cdot(\sigma\cdot\tau)\\
\sigma\cdot(\tau\cdot v)=\sigma\cdot(v_2,v_3,v_1)=(v_2,v_1,v_3)=(\sigma\cdot \tau)\cdot v\, .
\]
\end{Example}
The previous actions on $V^{\otimes N}$ induce dual right and left actions of the symmetric group $\mathbb S_N$ on $\Hom(V^{\otimes N},V)$:
\[ \begin{array}{lll}
{\Hom(V^{\otimes N},V)\times \mathbb S_N\to \Hom(V^{\otimes N},V)} &&
{(f\cdot \sigma)(v_1,\ldots,v_N)=f(v_{\sigma\un(1)},\ldots, v_{\sigma\un(N)})} \\
{\mathbb S_N \times \Hom(V^{\otimes N},V)\to \Hom(V^{\otimes N},V)} &&
{(\sigma \cdot f)(v_1,\ldots,v_N)=f(v_{\sigma(1)},\ldots, v_{\sigma(N)})}
\end{array} \]
\begin{Example}
Let  $\sigma,\tau\in \mS_3$ as in Example \ref{exasym} and $f\in\Hom(V^{\otimes3},V)$. Denoting $v:=(v_1,v_2,v_3)$, one can check that:
\begin{eqnarray}
\nn
&&\big((f\cdot \sigma)\cdot \tau\big)v=(f\cdot \sigma)(v_2,v_3,v_1)=(v_2,v_1,v_3)=\big(f\cdot(\sigma\cdot \tau)\big)v\\
&&\big(\sigma\cdot(\tau\cdot f)\big)v=(\tau\cdot f)(v_1,v_3,v_2)=f(v_2,v_1,v_3)=\big((\sigma\cdot\tau)\cdot f\big)v.\nn
\end{eqnarray}
\end{Example}
In the following, we will denote $\sgn_N$ the \textbf{signature representation} of $\mS_N$ \ie the one-dimensional $\corps\dl\mathbb S_N\dr$-module associating to each permutation $\sigma\in\mathbb S_N$ its signature $|\sigma|\in\pset{-1,1}$. A collection of right $\corps\dl\mathbb S_N\dr$-modules $M(N)$ for $N\geq1$ will be referred to as a $\mathbb S$-\textbf{module}. 
\paragraph{(Un)shuffles.}

Let $p,q\in\mathbb N$. A \textbf{shuffle} of type $(p,q)$ is a permutation $\sigma\in \mS_{p+q}$ such that $\sigma$ sends $\pset{1,\ldots, p+q}$ to $\pset{i_1,\ldots, i_p\, |\, j_1,\ldots ,j_q}$ where $i_1<\cdots<i_p$ and $j_1<\cdots<j_q$. 
\begin{Example}[Shuffles]
\hfill
\begin{itemize}
\item $\text{\rm Sh}(1,1)=\pset{(1|2),(2|1)}$
\item $\text{\rm Sh}(1,2)=\pset{(1|23),(2|13),(3|12)}$
\item $\text{\rm Sh}(2,1)=\pset{(12|3),(13|2),(23|1)}$
\item $\text{\rm Sh}(1,3)=\pset{(1|234),(2|134),(3|124),(4|123)}$
\item $\text{\rm Sh}(2,2)=\pset{(12|34),(13|24),(14|23),(23|14),(24|13),(34|12)}$
\item $\text{\rm Sh}(3,1)=\pset{(123|4),(124|3),(134|2),(234|1)}$
\end{itemize}
\end{Example}
The set of shuffles of type $(p,q)$ is denoted $\text{\rm Sh}(p,q)$. Since a shuffle $\sigma\in\text{\rm Sh}(p,q)$ is completely determined by the set $\pset{i_1,\ldots, i_p}$, there are $\binom{p+q}{p}$ shuffles of type $(p,q)$. A \textbf{unshuffle} of type $(p,q)$ is a permutation $\sigma\in \mS_{p+q}$ such that the inverse permutation $\sigma\un$ is a shuffle of type $(p,q)$. 
The set of unshuffles of type $(p,q)$ is denoted $\text{\rm Sh}\un(p,q)$. 

\begin{Example}[Unshuffles]
\hfill
\begin{itemize}
\item $\text{\rm Sh}\un(1,1)=\pset{(12),(21)}$
\item $\text{\rm Sh}\un(1,2)=\pset{(123),(213),(231)}$
\item $\text{\rm Sh}\un(2,1)=\pset{(123),(132),(312)}$
\item $\text{\rm Sh}\un(1,3)=\pset{(1234),(2134),(2314),(2341)}$
\item $\text{\rm Sh}\un(2,2)=\pset{(1234),(1324),(1342),(3124),(3142),(3412)}$
\item $\text{\rm Sh}\un(3,1)=\pset{(1234),(1243),(1423),(4123)}$
\end{itemize}
\end{Example}
\paragraph{Operads.}
We will consider operads in the category of (graded) vector spaces over $\corps$. Our conventions will mostly follow the ones of the book \cite{Loday2012}. We will denote $\mAss$, $\mLie$ and $\mCom$ the operads of (graded) vector spaces encoding (graded) associative, Lie and commutative associative algebras without unit, respectively. The cooperad governing cocommutative coassociative algebras without counit will be denoted $\mathsf{coCom}$. The latter is defined explicitly as:
\[
\mathsf{coCom}=
\begin{cases}
\mathbf{0}\text{ for }n=0\\
\corps \text{ for all }n>0
\end{cases}
\]
where $\corps$ stands for the trivial representation of $\mathbb S_n$. 

\noindent Letting $\mathcal O$ be an operad in the category of graded vector spaces, the set of graded vector spaces:
\begin{eqnarray}
\nn
\mathcal O\{\dime\}(N):=\begin{cases}
\mathcal O(N) [\dime(1-N)]\text{ for }\dime\text{ even}\\
\mathcal O(N)\otimes \sgn_N[\dime(1-N)]\text{ for }\dime\text{ odd}
\end{cases}
\end{eqnarray}
assemble to a $\mathbb S$-module. Endowing this $\mathbb S$-module with the partial composition maps, identity and right-actions of $\mathcal O$ defines the $\dime$-\textbf{suspended operad} $\mathcal O\{\dime\}$. Alternatively, the $\dime$-suspended operad $\mathcal O\{\dime\}$ can be characterised as the unique operad for which the set of algebras of the operad $\mathcal O$ on a graded vector space $V$ are in one-to-one correspondence with the set of algebras of $\mathcal O\{\dime\}$  on the suspended graded vector space $V[\dime]$. In particular, $\End_{V}\{d\}=\End_{V[d]}$ where $\End_V$ denotes the endomorphism operad associated with the graded vector space $V$.

\section{Graded geometry}
\label{Section:Symplectic Lie $n$-algebroids}
The aim of the present section is to provide a short introduction to graded manifolds as well as their (non-graded) geometric counterparts.
The latter objects are defined as the geometrical data associated with graded\footnote{We will only deal with $\mathbb N$-\textbf{graded manifolds} for which the corresponding degree assigned to each local coordinate is a non-negative integer. } manifolds
 -- understood as manifolds endowed with a grading of the corresponding structure sheaf [\cf \cite{Mehta2007} for precise definitions] -- supplemented with some additional graded structures. 
\subsection{$\mathsf{N}$-manifolds}

Letting $\V$ be a $\mathbb N$-graded manifold, or $\mathsf{N}$-manifold, of degree $n$ (\ie concentrated in degrees $0,\ldots ,n$) we will denote $\fonc{\V}$ the associated algebra of functions. The subvector space of homogeneous functions of degree $k$ will be denoted $\foncg{k}\subset\fonc{\V}$ so that $\fonc{\V}=\bigoplus_{k\geq0}\foncg{k}$ is a graded algebra. Moreover, $\fonc{\V}$ is a filtered algebra. Letting $\mathcal A_k$ denote the (graded) subalgebra of $\fonc{\V}$ locally generated by functions of degree $\leq k$, there is an increasing sequence: 
\begin{eqnarray}
\label{filtration}
\fonc{\M}=\mathcal A_0\subset\mathcal A_1\subset\cdots\subset\mathcal A_n=\fonc{\V}
\end{eqnarray}
where we have $\foncg{k}=\mathcal A_k / \mathcal A_{k-1}$ so that $\fonc{\V}=\bigoplus_{k\geq0}\foncg{k}$ is the graded algebra associated with the filtration \eqref{filtration}.
Corresponding to this filtration, there is a tower of fibrations:
\begin{eqnarray}
\label{fibration}
\M=\M_0\leftarrow\M_1\leftarrow\cdots\leftarrow\M_n=\V
\end{eqnarray}
where $\M$ is an ordinary smooth manifold -- referred to as the \textbf{base} -- and such that $\fonc{\M}=\foncg{0}=\mathcal A_0$. Furthermore, $\M_1$ is a vector bundle over $\M$ and for all $k\geq1$, $\M_k\leftarrow \M_{k+1}$ is an affine fibration, \cf \cite{Roytenberged.Contemp.Math.Vol.315Amer.Math.Soc.ProvidenceRI2002} for details.

The geometry of the fibration underlying graded manifolds can be enriched by introducing additional (hierarchised) data on $\V$ (\cf \cite{Roytenberged.Contemp.Math.Vol.315Amer.Math.Soc.ProvidenceRI2002,CattaneoSchaetz2010}):
\begin{itemize}
\item A $\NP$-\textbf{manifold $(\V,\om)$ of degree $n$} is a $\mathbb N$-graded manifold $\V$ endowed with a symplectic 2-form $\om$ of intrinsic degree $n$.
\item A $\NPQ$-\textbf{manifold $(\V,\om,\Q)$ of degree $n$} is a $\mathbb N$-graded manifold $\V$ endowed with a symplectic 2-form $\om$ of intrinsic degree $n$ and a homological vector field $\Q$ (\ie $\Q$ is of degree $1$ and satisfies $\Q^2=0$) such that $\Lag_Q\om=0$.
\end{itemize}
These additional data induce some extra geometric structures on the fibration \eqref{fibration}. 
We will refer to the (non-graded) geometrical data associated with {$\NPQ$-manifolds} of degree $n$ as \textbf{symplectic Lie $n$-algebroids}.
\subsection{$\NP$-manifolds}
\label{section:NPmanifolds}

Endowing a graded manifold with a symplectic (\ie non-degenerate and closed) 2-form has a number of consequences. 
First of all, the existence of a symplectic $2$-form of degree $n$ on a $\mathbb N$-graded manifold $\V$ constrains the degree of $\V$ to not exceed $n$ [\cf \cite{Roytenberged.Contemp.Math.Vol.315Amer.Math.Soc.ProvidenceRI2002} Lemma 2.4]. Secondly, it can be shown than any homogeneous symplectic 2-form of degree $n\geq1$ is exact [\cf \cite{Roytenberged.Contemp.Math.Vol.315Amer.Math.Soc.ProvidenceRI2002} Lemma 2.2]. These two properties can be used in order to provide a local presentation of $\NP$-manifolds \ala Darboux.

We distinguish between odd and even cases as follows. When $n$ is odd, we introduce a set of homogeneous coordinates\footnote{The subscript denotes the corresponding degree.}:
\[
\pset{\underset{0}{x^\mu},\underset{i}{\psi^{\alpha_i}},\underset{n-i}{\chi_{\alpha_i}},\underset{n}{p_\mu}}\quad\text{ where }\quad i\in\pset{1,\dots,\half(n-1)}.
\]
The symplectic 2-form of odd degree $n$ can thus be written as: 
\begin{eqnarray}
\nn
\om=dx^\mu\w dp_\mu+\sum_{i=1}^{\half(n-1)}d\psi^{\alpha_i}\w d\chi_{\alpha_i}\, .
\end{eqnarray}
The associated Poisson bracket of degree $-n$ acts as follows:
\begin{equation}\label{PBodd}
\begin{split}
\pb{f}{g}_\om &= (-1)^k\frac{\p f}{\p x^\mu}\, \frac{\p g}{\p p_\mu}+\frac{\p f}{\p p_\mu}\, \frac{\p g}{\p x^\mu} \\
&+ \sum_{i=1}^{\half(n-1)}\Bigg\{(-1)^{k(i+1)}\frac{\p f}{\p \psi^{\alpha_i}}\frac{\p g}{\p \chi_{\alpha_i}}+(-1)^{ik}\frac{\p f}{\p \chi_{\alpha_i}}\frac{\p g}{\p \psi^{\alpha_i}}\Bigg\}
\end{split}
\end{equation}
on homogeneous functions $f\in \foncg{k}$ and $g\in \foncg{l}$ of degree $k$ and $l$, respectively.

When $n$ is even, the corresponding set of homogeneous coordinates reads 
\[
\pset{\underset{0}{x^\mu},\underset{i}{\psi^{\alpha_i}},\underset{n/2}{\xi^a},\underset{n-i}{\chi_{\alpha_i}},\underset{n}{p_\mu}}\quad\text{ where }\quad i\in\pset{1,\dots,\half n-1}\, .
\]
The symplectic 2-form of even degree $n$ is written as: 
\begin{eqnarray}
\nn
\om=dx^\mu\w dp_\mu+\sum_{i=1}^{\half n-1}d\psi^{\alpha_i}\w d\chi_{\alpha_i}+\half \kappa_{ab}\,  d\xi^a\w d\xi^b
\end{eqnarray}
where the bilinear form $\kappa$ is non-degenerate and symmetric (resp. skewsymmetric\footnote{Note that, whenever $n=4\, k$ (for some integer $k=0,1,2,\ldots$) the indices of type $a,b,\ldots$ should run over an even number of dimensions in order to ensure the existence of a skewsymmetric invertible bilinear form $\kappa$.}) for $n/2$ odd (resp. even) \ie $\kappa_{ab}=-(-1)^{n/2}\kappa_{ba}$. The associated Poisson bracket thus takes the form:
\begin{equation}\label{PBeven}
\begin{split}
\pb{f}{g}_\om &= \frac{\p f}{\p x^\mu}\, \frac{\p g}{\p p_\mu}-\frac{\p f}{\p p_\mu}\, \frac{\p g}{\p x^\mu}\\
&+\sum_{i=1}^{\half n-1}\Bigg\{(-1)^{ik}\frac{\p f}{\p \psi^{\alpha_i}}\frac{\p g}{\p \chi_{\alpha_i}}-(-1)^{i(k+1)}\frac{\p f}{\p \chi_{\alpha_i}}\frac{\p g}{\p \psi^{\alpha_i}}\Bigg\}+(-1)^{k{n}/{2}}\frac{\p f}{\p \xi^a}\, \kappa^{ab}\, \frac{\p g}{\p \xi^b}.
\end{split}
\end{equation}
It can be checked that the Poisson brackets \eqref{PBodd} and \eqref{PBeven} satisfy the following properties:
\begin{enumerate}
\item $\pb{f}{g}_{\om}=-(-1)^n(-1)^{kl}\pb{g}{f}_{\om}$
\item $\pb{f}{g\cdot h}_{\om}=\pb{f}{g}_{\om}\cdot h+(-1)^{l(k-n)}g\cdot \pb{f}{h}_{\om}$
\item $\pb{\pb{f}{g}_{\om}}{h}_{\om}+(-1)^{k(l+m)}\pb{\pb{g}{h}_{\om}}{f}_{\om}+(-1)^{m(k+l)}\pb{\pb{h}{f}_{\om}}{g}_{\om}=0$
\end{enumerate}
for all homogeneous functions $f\in \foncg{k}$, $g\in \foncg{l}$ and $h\in \foncg{m}$ of degree $k$, $l$ and $m$ respectively.\footnote{In other words, the triplet $\big(\fonc{\V},\cdot,\pbdot_{\om}\big)$ is a $\mGer_{n+1}$-algebra, \cf Section \ref{Stable structures}.}

We conclude this quick survey of $\NP$-manifolds by discussing the notion of gauge transformations.

\begin{Definition}[Gauge transformation]
\label{defi:gauge transfo}
Let $(\V,\om)$ be a $\NP$-manifold. A gauge transformation on $\V$ is a diffeomorphism $\phi:\V\to\V$ being 
\begin{enumerate}
\item degree preserving
\item symplectomorphic \ie $\phi^*\om=\om$.
\end{enumerate}
\end{Definition}
\noindent In other words, the pullback map $\phi^*:\fonc{\V}\to\fonc{\V}$ is:
\begin{enumerate}
\item of degree 0 \ie $|\phi^*(f)|=|f|$
\item a morphism of Poisson algebras \ie $\phi^*(f\cdot g)=\phi^*(f)\cdot \phi^*(g)$ and $\phi^*(\pb{f}{g}_\om)=\big\{\phi^*(f),$ $\phi^*(g)\big\}_\om$
for all homogeneous functions $f,g\in\fonc{\V}$.
\end{enumerate}
Two functions $\mathcal F, \mathcal F'\in\fonc{\V}$ will be said to be equivalent if there exists a gauge transformation $\phi:\V\to\V$ such that $\mathcal F'=\phi^*(\mathcal F)$. Equivalence classes thereof will be denoted $[\mathcal F]$. 
Infinitesimal gauge transformations are symplectic vector fields $\mathscr{X}=\pb{f}{\cdot}_\om$ of degree 0 generated by arbitrary functions ${f}\in\foncg{n}$ of degree $n$. 
Two equivalent functions $\mathcal F, \mathcal F'\in[\mathcal F]$ differ infinitesimally by a term of the form $\pb{\mathcal F}{f}_\om$, with ${f}\in\foncg{n}$.

\subsection{$\NPQ$-manifolds}
\label{section:NPQmanifolds}

We now turn to $\NPQ$-manifolds and start by pointing out that the latter can be equivalently described in terms of a Poisson bracket together with a Hamiltonian function\footnote{Indeed, it follows from Cartan's homotopy formula that the compatibility relation between the symplectic $2$-form and the vector field ensures that the latter is Hamiltonian, \cf Lemma 2.2 in \cite{Roytenberged.Contemp.Math.Vol.315Amer.Math.Soc.ProvidenceRI2002}. } \ie as a triplet $(\mathcal V,\pbdot_\om,\cH)$ where:
\begin{enumerate}
\item $\mathcal V$ is a $\mathbb N$-graded manifold. 
\item $\pbdot_\om$ is a non-degenerate Poisson bracket  of degree $-n$ acting on the graded algebra of functions on $\V$. 
\item $\cH$ is a Hamiltonian function \ie a homogeneous function of degree $n+1$ being nilpotent with respect to the graded Poisson bracket \ie $\pb{\cH}{\cH}_\om=0$. The set of Hamiltonian functions will be denoted $\Ham$.
\end{enumerate}
Equivalence between the homological and Hamiltonian presentations of $\NPQ$-manifolds of degree $n$ is realised by identifying $\om$ as the symplectic $2$-form of degree $n$ dual to $\pbdot_\om$ and defining the privileged vector field $\Q\in\Gamma^1(T\V)$ of degree 1 on $\V$ as $\Q=\pb{\cH}{\cdot}_\om$. The nilpotency of $\cH$ ensures that $\Q$ is homological \ie $\br{\Q}{\Q}=0$, with $\brdot$ the graded Lie bracket on $\V$.

The importance of $\NPQ$-manifolds (or equivalently symplectic Lie $n$-algebroids) stems from the fact that these naturally form the target space of the classical action associated with AKSZ-type $\sigma$-models \cite{AlexandrovKontsevichSchwarzEtAl1997} for which the source manifold has dimension $d=n+1$.

Two Hamiltonian functions $\cH, \cH'\in\Ham$ will be said to be equivalent if there exists a gauge transformation [see \Defi{defi:gauge transfo}] denoted $\phi:\V\to\V$ such that $\cH'=\phi^*(\cH)$. Two equivalent Hamiltonians differ infinitesimally by a coboundary $\pb{\cH}{f}_\om$, with ${f}\in\foncg{n}$.

We conclude this brief survey by displaying examples of symplectic Lie $n$-algebroids in low degrees.

\begin{Example}[Symplectic Lie $n$-algebroids]
\label{exasymp}
\hfill
\begin{itemize}
\item $n=0$ (Symplectic manifolds)

The manifold is coordinatised by a unique set of homogeneous coordinates $\xi^a$ of degree 0, with $a\in\pset{1,\ldots,D}$ and $D$ the (even) dimension of the manifold. The manifold is thus non-graded (or bosonic) \ie $\V$ identifies with its base $\M$. The symplectic 2-form of degree $0$ takes the usual form $\om=\half \kappa_{ab}\,  d\xi^a\w d\xi^b$ where the bilinear form $\kappa$ is non-degenerate and skewsymmetric \ie $\kappa_{ab}=-\kappa_{ba}$. The associated Poisson bracket thus takes the form:
\begin{eqnarray}
\nn
\pb{f}{g}_\om=\frac{\p f}{\p \xi^a}\, \kappa^{ab}\, \frac{\p g}{\p \xi^b}\, .
\end{eqnarray}
The absence of degree 1 coordinates prevents the existence of a Hamiltonian function $\cH$ (of would-be degree 1) in this case. Symplectic Lie $0$-algebroids are thus in one-to-one correspondence with (ordinary) symplectic manifolds. Gauge transformations identify in this case with usual symplectomorphisms. 

\item $n=1$ (Poisson manifolds)

The set of homogeneous coordinates takes the form $\pset{\underset{0}{x^\mu},\underset{1}{p_\mu}}$. The symplectic 2-form of odd degree $1$ can thus be written as $\om=dx^\mu\w dp_\mu$
while the associated Poisson bracket of degree $-1$ acts as follows:
\begin{eqnarray}
\nn
\pb{f}{g}_\om=(-1)^k\frac{\p f}{\p x^\mu}\, \frac{\p g}{\p p_\mu}+\frac{\p f}{\p p_\mu}\, \frac{\p g}{\p x^\mu}
\end{eqnarray}
on homogeneous functions $f\in \foncg{k}$ and $g\in \foncg{l}$ of degree $k$ and $l$ respectively. Up to degree suspension, $\pbdot_\om$ identifies with the Schouten bracket acting on polyvector fields. The most general function of degree $2$ reads $\cH=\frac{1}{2}\pi^{\mu\nu}(x)p_\mu\,  p_\nu$ with $\pi$ a bivector, \ie $\pi^{\mu \nu}=-\pi^{\nu\mu}$. It can be checked that\footnote{Here and in the following, round (resp.\ square) brackets of indices will denote (skew)symmetrisation.} $\pb{\cH}{\cH}_\om=0\Leftrightarrow \pi^{\rho[\lambda}\p_\rho\pi^{\mu\nu]}=0$ \ie $\cH$ is Hamiltonian if and only if $\pi$ is a Poisson bivector. It follows that symplectic Lie $1$-algebroids are in one-to-one correspondence with Poisson manifolds. Gauge transformations identify in this case with usual diffeomorphisms on the base manifold. 
\item $n=2$ (Courant algebroids)

The set of homogeneous coordinates can be decomposed as $\pset{\underset{0}{x^\mu},\underset{1}{\xi^a},\underset{2}{p_\mu}}$. The symplectic 2-form of even degree $2$ can thus be written as $\om=dx^\mu\w dp_\mu+\half \kappa_{ab}\,  d\xi^a\w d\xi^b$
where the bilinear form $\kappa$ is non-degenerate and symmetric \ie $\kappa_{ab}=\kappa_{ba}$.
The associated Poisson bracket of degree $-2$ acts as follows:
\begin{eqnarray}
\nn
\pb{f}{g}_\om=\frac{\p f}{\p x^\mu}\, \frac{\p g}{\p p_\mu}-\frac{\p f}{\p p_\mu}\, \frac{\p g}{\p x^\mu}+(-1)^{k}\frac{\p f}{\p \xi^a}\, \kappa^{ab}\, \frac{\p g}{\p \xi^b}
\end{eqnarray}
on homogeneous functions $f\in \foncg{k}$ and $g\in \foncg{l}$ of degree $k$ and $l$ respectively. The most general function of degree $3$ reads $\cH=\rho_a{}^\mu\, \xi^a p_\mu+\frac16\, T_{abc}\, \xi^a\xi^b\xi^c$ where $T_{abc}$ is totally skewsymmetric. 
It can be checked that the nilpotency condition $\pb{\cH}{\cH}_\om=0$ is equivalent to the three following constraints:
\begin{enumerate}
\item ${\mathcal C_1}^{\mu \nu}:=\rho_{a}{}^{\mu}\kappa^{ab}\rho_{b}{}^{\nu}=0$\itemnum\label{ConstraintCA1}
\item ${\mathcal C_2}^{\mu}_{a b}:=\rho_{ c}{}^{\mu}\kappa^{cd}T_{d a b}+2\, \rho_{[a}{}^{\lambda}\, \p_{\lambda}\rho_{b]}{}^{\mu}=0$\itemnum\label{ConstraintCA2}
\item ${\mathcal C_3}_{a b c d}:=\frac{1}{4}T_{e [a b}\kappa^{ef}T_{c d] f}+\frac{1}{3}\rho_{[a}{}^{\mu}\, \p_\mu T_{bcd]}=0$.\itemnum\label{ConstraintCA3}
\end{enumerate}
As will be reviewed in Section \ref{section:Courant algebroids}, symplectic Lie $2$-algebroids are in one-to-one correspondence with Courant algebroids. Infinitesimal gauge transformations are generated by functions $f=X^\mu p_\mu+ \Lambda_{ab}\, \xi^a\xi^b$ where the two terms on the right-hand side correspond to infinitesimal diffeomorphisms on the base manifold and infinitesimal rotations on the fiber, respectively. 
\end{itemize}
\end{Example}

\section{Graph complexes}
\label{section:Graph complexes}
The aim of the present section is to review a particular family of graph complexes\footnote{Graph complexes come in many variants. As shown in \cite{Getzler1994b,Markl1998,Conant2002}, to any cyclic operad $\mathcal O$ one can associate a class of $\mathcal O$-graph complexes. In particular, $\mathcal O=\mAss$ corresponds to the class of ribbon graphs computing cohomology of moduli spaces of curves \cite{Penner1987,Penner1988} while the graph complex for $\mathcal O=\mLie$ computes cohomology of outer automorphisms of free groups \cite{Culler1986}. We will solely be interested in the case $\mathcal O=\mCom$. Also, graph complexes come in two dual versions: a \textit{homological} version in which the boundary operator acts via ``collapsing'' of edges \cite{Kontsevich1993,Kontsevich1994} and a \textit{cohomological} one in which the coboundary operator acts by ``blowing up'' edges \cite{Kontsevich1997,Willwacher2015}. We will hereafter focus on the cohomological version. } introduced by M. Kontsevich in  \cite{Kontsevich1993,Kontsevich1994,Kontsevich1997}. The former is most clearly defined in terms of the convolution Lie algebra constructed from a suitable graph operad. We start by reviewing the construction of this graph operad -- denoted $\Gra_d$ hereafter\footnote{In Section \ref{section:Stable structures on graded manifolds}, we will relate the integer $d$ (in the case when $d\geq1$) to the dimension of the source of the relevant AKSZ $\sigma$-model on which $\Gra_d$ will be shown to act. In other words, we will consider $d=n+1$ where $n$ is the degree of the corresponding $\NPQ$-manifold, \cf Section \ref{section:NPQmanifolds}. } -- 
 from a combinatorial point of view before turning to the definition of the so-called full graph complex $\fGC_d$. After reviewing results regarding the cohomology of $\fGC_d$, we conclude by presenting a variant of the full graph complex whose elements are directed graphs. 
  The material covered in this section is standard and can be found for example in \cite{Willwacher2015,Dolgushev2012a,Dolgushev2012b}.

\subsection{The graph operad $\Gra_\dime$}
\label{section:graphoperad}
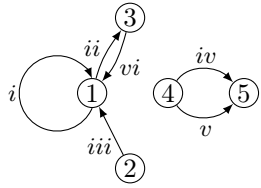
\begin{figure}[h]
  \begin{center}
\begin{tikzpicture}[scale=0.5, >=stealth']
\tikzstyle{w}=[circle, draw, minimum size=4, inner sep=1]
\tikzstyle{b}=[circle, draw, fill, minimum size=4, inner sep=1]
\node [ext] (b1) at (0,0) {1};
\node [ext] (b2) at (1,2) {3};
\node [ext] (b3) at (1,-2) {2};
\node [ext] (b4) at (2,0) {4};
\node [ext] (b5) at (4,0) {5};
\draw (0,1.2) node[anchor=center] {{\small $ii$}};
\draw (b1) edge[>=latex,<-] (b3);
\draw (0.2,-1.4) node[anchor=center] {{\small $iii$}};
\draw (-2.1,0) node[anchor=center] {{\small $i$}};
    \draw [>=latex,->]  (0,-0.4)arc(157:-157:-1);
        \draw (3,1) node[anchor=center] {{\small $iv$}};
                \draw (3,-1) node[anchor=center] {{\small $v$}};
\draw[black,->,>=latex] (b4) to[out=55, in=135, looseness=1.] (b5);
\draw[black,->,>=latex] (b4) to[out=-55, in=-135, looseness=1.] (b5);
\draw[black,->,>=latex] (b1) to[out=75, in=-125, looseness=1.] (b2);
\draw[black,<-,>=latex] (b1) to[out=50, in=-105, looseness=1.] (b2);
\draw (1,0.7) node[anchor=center] {{\small $vi$}};
\end{tikzpicture}
  \end{center}
  \caption{Example of graph in $\gra_{5,6}$}
  \label{figgraphex}
\end{figure}

Our starting point towards a definition of the graph operad $\Gra_\dime$ will be the set of \textbf{multidigraphs} (or quivers) \ie directed graphs which are allowed to contain multiple edges and tadpoles\footnote{Formally, a multidigraph is defined as a four-tuple $\Gamgraph=(V_\Gamgraph,E_\Gamgraph,s,t)$ where:
    \begin{itemize}
  \item $V_\Gamgraph$ is a set whose elements are called \textbf{vertices}. 
  \item $E_\Gamgraph$ is a set whose elements are called \textbf{edges}. 
  \item The map $s:E_\Gamgraph\to V_\Gamgraph$ assigns to each edge its \textbf{source}. 
    \item The map $t:E_\Gamgraph\to V_\Gamgraph$ assigns to each edge its \textbf{target}. 
      \end{itemize}
   An edge $e\in E_\Gamgraph$ such that $s(e)=t(e)$ is called a \textbf{tadpole} while pairs of edges $e_1,e_2\in E_\Gamgraph$ such that $s(e_1)=s(e_2)$ and $t(e_1)=t(e_2)$ are called \textbf{double edges}. The set of edges connecting a given vertex $v\in V_\Gamgraph$ will be denoted $E_\Gamgraph(v)$. We will mostly deal with {\it labeled} multidigraphs \ie multidigraphs endowed with two bijective maps $l_V:V_\Gamgraph\to[|V_\Gamgraph|]$ and $l_E:E_\Gamgraph\to[|E_\Gamgraph|]$ where $|V_\Gamgraph|$ (resp. $|E_\Gamgraph|$) denotes the number of vertices (resp. edges) of $\Gamgraph$ and $[n]:=\pset{1,2,\ldots,n}$. While depicting multidigraphs pictorially, we will represent edges by arrows from source to target vertices. To avoid ambiguity, labelling will be performed using Hindu-Arabic numerals for vertices and Roman numerals for edges. Note that we do not assume any compatibility between the labelling of vertices and edges {\it a priori}.}. 
      The set of multidigraphs with $N$ vertices and $k$ directed edges will be denoted $\gra_{N,k}$. A typical example\footnote{Note that the definition of a multidigraph does not assume connectedness. }of multidigraph is given in Figure \ref{figgraphex}. 
There is a natural right-action of the semi-direct product $\mS_k\ltimes \mS_2^{\times k}$ on elements of $\gra_{N,k}$ by permutation of the ordering $(\mathbb S_k)$ and flipping of the directions of the edges $(\mS_2^{\times k})$. We will consider the 1-dimensional signature representation $\sgn_k$ (resp. $\sgn_{2}^{\otimes k}$) as a left $\corps\dl\mS_k\ltimes \mS_2^{\times k}\dr$-module with trivial action of $\mS_2^{\times k}$ (resp. $\mS_k$).
For all $N\geq1$ and $d\in\mathbb N$, we define the collection of graded vector spaces $\Gra_\dime(N)$ as:
\begin{itemize}
\item $\dime$ even: \hspace{0.15cm}$\displaystyle\Gra_\dime(N):=\prod_{k\geq0}\big(\corps\dl\gra_{N,k}\dr\otimes_{\mS_k\ltimes \mS_2^{\times k}}\sgn_k\big)[k(\dime-1)]$\itemnum\label{defGraeven}
\item $\dime$ odd: \hspace{0.27cm}$\displaystyle\Gra_\dime(N):=\prod_{k\geq0}\big(\corps\dl\gra_{N,k}\dr\otimes_{\mS_k\ltimes \mS_2^{\times k}}\sgn_{2}^{\otimes k}\big)[k(\dime-1)]$\itemnum\label{defGraodd}
\end{itemize}
where the subscript stands for taking coinvariants with respect to the diagonal right action of ${\mS_k\ltimes \mS_2^{\times k}}$ and the term between brackets denotes degree suspension (\cf Section \ref{Section:Conventions and notations} for conventions).  Elements of $\Gra_\dime(N)$ are linear combinations of equivalence classes of graphs in $\gra_{N,k}$, for arbitrary $k\geq0$. 
Two graphs $\Gamgraph,\Gamgraph'\in\gra_{N,k}$ will be said equivalent (\ie $\Gamgraph\sim\Gamgraph'$) if one of the two following condition holds: 
\begin{enumerate}
\item There exists an element $\sigma\in \mS_2^{\times k}$ such that $\Phi^{\text{\rm dir}}_\sigma(\Gamgraph)=(-1)^{\dime|\sigma|}\Gamgraph'$ where $\Phi^{\text{\rm dir}}_\sigma$ stands for the automorphism of $\gra_{N,k}$ that flips the direction of the edges according to $\sigma$,
\begin{eqnarray}
\nn
\text{\eg}\quad\raisebox{-1ex}{\hbox{
\begin{tikzpicture}[scale=0.5, >=stealth']
\tikzstyle{w}=[circle, draw, minimum size=4, inner sep=1]
\tikzstyle{b}=[circle, draw, fill, minimum size=4, inner sep=1]
\node [ext] (b4) at (2,0) {1};
\node [ext] (b5) at (4,0) {2};
\draw (3,0.5) node[anchor=center] {{\small $i$}};
\draw[black,->,>=latex]  (b4) to (b5);
\end{tikzpicture}}}
\quad\sim\hspace{0.15cm}(-1)^d
\raisebox{-1ex}{\hbox{
\begin{tikzpicture}[scale=0.5, >=stealth']
\tikzstyle{w}=[circle, draw, minimum size=4, inner sep=1]
\tikzstyle{b}=[circle, draw, fill, minimum size=4, inner sep=1]
\node [ext] (b4) at (2,0) {1};
\node [ext] (b5) at (4,0) {2};
\draw (3,0.5) node[anchor=center] {{\small $i$}};
\draw[black,<-,>=latex]  (b4) to (b5);
\end{tikzpicture}}}\quad.
\end{eqnarray}
\item There exists an element $\sigma\in \mS_k$ such that $\Phi^{\text{\rm order}}_\sigma(\Gamgraph)=(-1)^{(\dime+1)|\sigma|}\Gamgraph'$ where $\Phi^{\text{\rm order}}_\sigma$ stands for the automorphism of $\gra_{N,k}$ that permutes the order of the edges according to $\sigma$,
\begin{eqnarray}
\nn
\text{\eg}\quad\raisebox{-1ex}{\hbox{
\begin{tikzpicture}[scale=0.5, >=stealth']
\tikzstyle{w}=[circle, draw, minimum size=4, inner sep=1]
\tikzstyle{b}=[circle, draw, fill, minimum size=4, inner sep=1]
\node [ext] (b4) at (2,0) {1};
\node [ext] (b5) at (4,0) {2};
\node [ext] (b6) at (6,0) {3};
\draw (3,0.5) node[anchor=center] {{\small $i$}};
\draw (5,0.5) node[anchor=center] {{\small $ii$}};
\draw[black,->,>=latex]  (b4) to (b5);
\draw[black,->,>=latex]  (b5) to (b6);
\end{tikzpicture}}}
\quad\sim\hspace{0.15cm}(-1)^{d+1}
\raisebox{-1ex}{\hbox{
\begin{tikzpicture}[scale=0.5, >=stealth']
\tikzstyle{w}=[circle, draw, minimum size=4, inner sep=1]
\tikzstyle{b}=[circle, draw, fill, minimum size=4, inner sep=1]
\node [ext] (b4) at (2,0) {1};
\node [ext] (b5) at (4,0) {2};
\node [ext] (b6) at (6,0) {3};
\draw (3,0.5) node[anchor=center] {{\small $ii$}};
\draw (5,0.5) node[anchor=center] {{\small $i$}};
\draw[black,->,>=latex]  (b4) to (b5);
\draw[black,->,>=latex]  (b5) to (b6);
\end{tikzpicture}}}\quad.
\end{eqnarray}
\end{enumerate}
According to the degree suspension in \eqref{defGraeven}-\eqref{defGraodd}, each edge is assigned an intrinsic degree $1-d$, so that the degree of an element $\Gamgraph\in\gra_{N,k}$ as seen in $\Gra_\dime(N)$ is given by $|\Gamgraph|=k(1-\dime)$. It is also clear from their definition that graded vector spaces $\Gra_\dime(N)$ for different $\dime$ of same parity only differ by their degree assignment and are thus isomorphic to each other. Following \cite{Rutten2018}, we will call \textbf{zero graph} a graph $\Gamgraph\in\gra_{N,k}$ which equals minus itself in $\Gra_\dime(N)$ and thus belongs to the zero class in $\Gra_\dime(N)$. It follows that a graph admitting an automorphism that permutes the edges ordering by an odd permutation is a zero graph whenever $\dime$ is even. In particular, graphs admitting multiple edges are zero graphs for $\dime$ even\footnote{Whenever $\dime$ is even, the double edges graph $\vcenter{\hbox{\begin{tikzpicture}[scale=0.5, >=stealth']
\tikzstyle{w}=[circle, draw, minimum size=4, inner sep=1]
\tikzstyle{b}=[circle, draw, fill, minimum size=4, inner sep=1]
\node [ext] (b4) at (2,0) {1};
\node [ext] (b5) at (4,0) {2};
        \draw (3,1) node[anchor=center] {{\small $i$}};
                \draw (3,-1) node[anchor=center] {{\small $ii$}};
\draw[black] (b4) to[out=55, in=125, looseness=1.] (b5);
\draw[black] (b4) to[out=-55, in=-125, looseness=1.] (b5);
\end{tikzpicture}}}$
satisfies $\vcenter{\hbox{\begin{tikzpicture}[scale=0.5, >=stealth']
\tikzstyle{w}=[circle, draw, minimum size=4, inner sep=1]
\tikzstyle{b}=[circle, draw, fill, minimum size=4, inner sep=1]
\node [ext] (b4) at (2,0) {1};
\node [ext] (b5) at (4,0) {2};
        \draw (3,1) node[anchor=center] {{\small $i$}};
                \draw (3,-1) node[anchor=center] {{\small $ii$}};
\draw[black] (b4) to[out=55, in=125, looseness=1.] (b5);
\draw[black] (b4) to[out=-55, in=-125, looseness=1.] (b5);
\end{tikzpicture}}}\sim-\vcenter{\hbox{\begin{tikzpicture}[scale=0.5, >=stealth']
\tikzstyle{w}=[circle, draw, minimum size=4, inner sep=1]
\tikzstyle{b}=[circle, draw, fill, minimum size=4, inner sep=1]
\node [ext] (b4) at (2,0) {1};
\node [ext] (b5) at (4,0) {2};
        \draw (3,1) node[anchor=center] {{\small $ii$}};
                \draw (3,-1) node[anchor=center] {{\small $i$}};
\draw[black] (b4) to[out=55, in=125, looseness=1.] (b5);
\draw[black] (b4) to[out=-55, in=-125, looseness=1.] (b5);
\end{tikzpicture}}}$ and is thus a zero graph in $\Gra_d(2)$.}. On the other hand, a graph admitting an automorphism that flips an odd number of edges is automatically a zero graph whenever $\dime$ is odd. In particular, graphs with tadpoles are zero graphs for $\dime$ odd\footnote{ Whenever $\dime$ is odd, the tadpole graph $\vcenter{\hbox{\begin{tikzpicture}[scale=0.5, >=stealth']
\tikzstyle{w}=[circle, draw, minimum size=4, inner sep=1]
\tikzstyle{b}=[circle, draw, fill, minimum size=4, inner sep=1]
\node [ext] (b1) at (0,0) {1};
    \draw [>=latex,->]  (0,-0.4)arc(157:-157:-1);
\end{tikzpicture}}}$
satisfies $\vcenter{\hbox{\begin{tikzpicture}[scale=0.5, >=stealth']
\tikzstyle{w}=[circle, draw, minimum size=4, inner sep=1]
\tikzstyle{b}=[circle, draw, fill, minimum size=4, inner sep=1]
\node [ext] (b1) at (0,0) {1};
    \draw [>=latex,->]  (0,-0.4)arc(157:-157:-1);
\end{tikzpicture}}}\sim-\vcenter{\hbox{\begin{tikzpicture}[scale=0.5, >=stealth']
\tikzstyle{w}=[circle, draw, minimum size=4, inner sep=1]
\tikzstyle{b}=[circle, draw, fill, minimum size=4, inner sep=1]
\node [ext] (b1) at (0,0) {1};
    \draw [>=latex,<-]  (0,-0.4)arc(157:-157:-1);
\end{tikzpicture}}}$ and is thus a zero graph in $\Gra_\dime(1)$.}.

 For all $N\geq1$, the symmetric group $\mS_N$ acts naturally on the right on the graded vector space $\Gra_\dime(N)$ by permuting the label of vertices as $\pset{1,2,\ldots,N}\overset{\sigma}{\mapsto}\pset{\sigma\un(1),\sigma\un(2),\ldots,\sigma\un(N)}$. We will denote $\Sigma_N:\Gra_d(N)\times\mathbb S_N\to\Gra_d(N)$ the corresponding right action\footnote{For example, letting $\sigma\in \mS_3$ be defined as 
\[
\sigma:=\begin{pmatrix}
1 &2&3\\
3 &1&2
\end{pmatrix}\, ,
\] the right-action  of $\sigma$ on the graph $\Gamgraph\in\Gra_\dime(3)$ defined as: 
\[
\Gamgraph:=\raisebox{-0.8ex}{\hbox{
\begin{tikzpicture}[scale=0.5, >=stealth']
\tikzstyle{w}=[circle, draw, minimum size=4, inner sep=1]
\tikzstyle{b}=[circle, draw, fill, minimum size=4, inner sep=1]
\node [ext] (b4) at (2,0) {1};
\node [ext] (b5) at (4,0) {2};
\node [ext] (b6) at (6,0) {3};
\draw (3,0.5) node[anchor=center] {{\small $i$}};
\draw[black,->,>=latex]  (b4) to (b5);
\draw (5,0.5) node[anchor=center] {{\small $ii$}};
\draw[black,->,>=latex]  (b5) to (b6);
\end{tikzpicture}}}
\]
 reads:
 \[
 \Sigma_3(\Gamgraph|\sigma)=\raisebox{-0.8ex}{\hbox{
\begin{tikzpicture}[scale=0.5, >=stealth']
\tikzstyle{w}=[circle, draw, minimum size=4, inner sep=1]
\tikzstyle{b}=[circle, draw, fill, minimum size=4, inner sep=1]
\node [ext] (b4) at (2,0) {2};
\node [ext] (b5) at (4,0) {3};
\node [ext] (b6) at (6,0) {1};
\draw (3,0.5) node[anchor=center] {{\small $i$}};
\draw[black,->,>=latex]  (b4) to (b5);
\draw (5,0.5) node[anchor=center] {{\small $ii$}};
\draw[black,->,>=latex]  (b5) to (b6);
\end{tikzpicture}}}\ .
\]}. In other words, the set of graded vector spaces $\pset{\Gra_\dime(N)}_{N\geq1}$ assemble to a {$\mathbb S$-module} over $\corps$. The $\mathbb S$-module $\pset{\Gra_\dime(N)}_{N\geq1}$ can further be given  the structure of an operad by endowing it with partial composition operations. Explicitly, we define partial composition operations: 
\begin{eqnarray}
\label{parcompGra}
\circ_i:\gra_{M,j}\otimes\gra_{N,k}\to\Gra_d(M+N-1)\text{ for all }1\leq i\leq M\ \text{ as }\ \Gamgraph\circ_i\Gamgraph'\ =\hspace{-0.5cm}\underset{f\, \in\, \Hom(E_{\Gamgraph}(v_i),V_{\Gamgraph'})}{\sum} \hspace{-0.7cm}\Gamgraph\circ_i^f \Gamgraph'
\end{eqnarray}
 where $\gamma\in\gra_{M,j}$ and $\gamma'\in\gra_{N,k}$. In the above formula, we let $v_i$ be the $i^\text{th}$ vertex of $\Gamgraph$ and the sum be performed over homomorphisms of sets between the set $E_{\Gamgraph}(v_i)$ of edges of $\Gamgraph$ connecting $v_i$ and the set $V_{\Gamgraph'}$ of vertices of $\Gamgraph'$. The operation $\circ_i^f$ consists in first inserting the graph $\Gamgraph'$ in place of the vertex $v_i\in\Gamgraph$ and then reconnecting the elements in $E_{\Gamgraph}(v_i)$ to vertices of $\Gamgraph'$ along the map $f$. The output is a sum of graphs with $j+k$ edges in $\Gra_d(M+N-1)$. As for labelling of vertices and edges, we follow the rules:
   \begin{itemize}
   \item The labels of the first $i-1$ vertices of $\Gamgraph$ are left unchanged.
   \item The labels of the vertices of $\Gamgraph'$ are shifted up by $i-1$. 
   \item The last $M-i$ vertices of $\Gamgraph$ are shifted up by $N-1$. 
   \item All edges originating from $\Gamgraph$ are declared smaller than all edges originating from $\Gamgraph'$. 
   \end{itemize}
   The partial composition operations $\circ_i$ can be checked to be equivariant with respect to the right-action of $\mS_k\ltimes \mS_2^{\times k}$ on $\gra_{N,k}$ allowing to define partial composition operations $\circ_i:\Gra_\dime(M)\otimes\Gra_\dime(N)\to\Gra_\dime(M+N-1)$.

\begin{figure}[h]
  \begin{center}
\begin{eqnarray}
\raisebox{-4ex}{\hbox{
\begin{tikzpicture}[scale=0.5, >=stealth']
\tikzstyle{w}=[circle, draw, minimum size=4, inner sep=1]
\tikzstyle{b}=[circle, draw, fill, minimum size=4, inner sep=1]
\node [ext] (b1) at (0,0) {1};
\node [ext] (b2) at (2,0) {2};
\node [ext] (b3) at (1,-1.73) {3};
\draw[black,->,>=latex]  (b1) to (b2);
\draw (1,0.5) node[anchor=center] {{\small $i$}};
\draw (2,-0.85) node[anchor=center] {{\small $ii$}};
\draw (-0.08,-0.9) node[anchor=center] {{\small $iii$}};
\draw[black,->,>=latex]  (b2) to (b3);
\draw[black,->,>=latex]  (b3) to (b1);
\end{tikzpicture}}}
\circ_2
\raisebox{-0.8ex}{\hbox{
\begin{tikzpicture}[scale=0.5, >=stealth']
\tikzstyle{w}=[circle, draw, minimum size=4, inner sep=1]
\tikzstyle{b}=[circle, draw, fill, minimum size=4, inner sep=1]
\node [ext] (b4) at (2,0) {1};
\node [ext] (b5) at (4,0) {2};
\draw (3,0.5) node[anchor=center] {{\small $i$}};
\draw[black,->,>=latex]  (b4) to (b5);
\end{tikzpicture}}}
=
\raisebox{-6ex}{\hbox{
\begin{tikzpicture}[scale=0.5, >=stealth']
\tikzstyle{w}=[circle, draw, minimum size=4, inner sep=1]
\tikzstyle{b}=[circle, draw, fill, minimum size=4, inner sep=1]
\node [ext] (b4) at (0,0) {4};
\node [ext] (b3) at (2,0) {3};
\node [ext] (b2) at (2,2) {2};
\node [ext] (b1) at (0,2) {1};
\draw (1,2.5) node[anchor=center] {{\small $i$}};
\draw (1,-0.5) node[anchor=center] {{\small $ii$}};
\draw (-0.5,1) node[anchor=center] {{\small $iii$}};
\draw (2.5,1) node[anchor=center] {{\small $iv$}};
\draw[black,->,>=latex]  (b1) to (b2);
\draw[black,->,>=latex]  (b2) to (b3);
\draw[black,->,>=latex]  (b3) to (b4);
\draw[black,->,>=latex]  (b4) to (b1);
\end{tikzpicture}}}
+
\raisebox{-6ex}{\hbox{
\begin{tikzpicture}[scale=0.5, >=stealth']
\tikzstyle{w}=[circle, draw, minimum size=4, inner sep=1]
\tikzstyle{b}=[circle, draw, fill, minimum size=4, inner sep=1]
\node [ext] (b4) at (0,0) {4};
\node [ext] (b3) at (2,0) {2};
\node [ext] (b2) at (2,2) {3};
\node [ext] (b1) at (0,2) {1};
\draw (1,2.5) node[anchor=center] {{\small $i$}};
\draw (1,-0.5) node[anchor=center] {{\small $ii$}};
\draw (-0.5,1) node[anchor=center] {{\small $iii$}};
\draw (2.5,1) node[anchor=center] {{\small $iv$}};
\draw[black,->,>=latex]  (b1) to (b2);
\draw[black,<-,>=latex]  (b2) to (b3);
\draw[black,->,>=latex]  (b3) to (b4);
\draw[black,->,>=latex]  (b4) to (b1);
\end{tikzpicture}}}
+
\raisebox{-4ex}{\hbox{
\begin{tikzpicture}[scale=0.5, >=stealth']
\tikzstyle{w}=[circle, draw, minimum size=4, inner sep=1]
\tikzstyle{b}=[circle, draw, fill, minimum size=4, inner sep=1]
\node [ext] (b1) at (0,0) {1};
\node [ext] (b2) at (2,0) {2};
\node [ext] (b3) at (1,-1.73) {4};
\node [ext] (b4) at (3.5,1.3) {3};
\draw[black,->,>=latex]  (b1) to (b2);
\draw (1,0.5) node[anchor=center] {{\small $i$}};
\draw (2,-0.85) node[anchor=center] {{\small $ii$}};
\draw (-0.08,-0.9) node[anchor=center] {{\small $iii$}};
\draw (2.4,1.25) node[anchor=center] {{\small $iv$}};
\draw[black,->,>=latex]  (b2) to (b3);
\draw[black,->,>=latex]  (b3) to (b1);
\draw[black,->,>=latex]  (b2) to (b4);
\end{tikzpicture}}}+
\raisebox{-4ex}{\hbox{
\begin{tikzpicture}[scale=0.5, >=stealth']
\tikzstyle{w}=[circle, draw, minimum size=4, inner sep=1]
\tikzstyle{b}=[circle, draw, fill, minimum size=4, inner sep=1]
\node [ext] (b1) at (0,0) {1};
\node [ext] (b2) at (2,0) {3};
\node [ext] (b3) at (1,-1.73) {4};
\node [ext] (b4) at (3.5,1.3) {2};
\draw[black,->,>=latex]  (b1) to (b2);
\draw (1,0.5) node[anchor=center] {{\small $i$}};
\draw (2,-0.85) node[anchor=center] {{\small $ii$}};
\draw (-0.08,-0.9) node[anchor=center] {{\small $iii$}};
\draw (2.4,1.25) node[anchor=center] {{\small $iv$}};
\draw[black,->,>=latex]  (b2) to (b3);
\draw[black,->,>=latex]  (b3) to (b1);
\draw[black,<-,>=latex]  (b2) to (b4);
\end{tikzpicture}}}\nn
\end{eqnarray}
  \end{center}
  \caption{Example of partial composition $\Gra_d(3)\circ_2\Gra_d(2)\to\Gra_d(4)$}
  \label{figtypicalodd}
\end{figure}
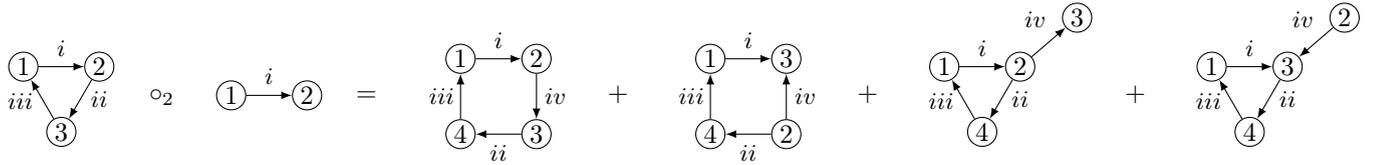

The partial composition operations $\circ_i$ on $\Gra_\dime$ preserve the number of edges and thus have zero intrinsic degree. 
Further, they can be checked to satisfy the following properties for all $\Gamgraph_m\in\Gra_\dime(m)$:
\begin{itemize}
\item \textbf{Sequential composition}:
\begin{eqnarray}
\nn
(\Gamgraph_m\circ_j\Gamgraph_n)\circ_{i}\Gamgraph_p=\Gamgraph_m\circ_j(\Gamgraph_n\circ_{i-j+1}\Gamgraph_p) \text{ for all }j\leqslant i\leqslant j+n-1\, . 
\end{eqnarray}
\item \textbf{Parallel composition}:
\begin{eqnarray}
\nn
(\Gamgraph_m\circ_j\Gamgraph_n)\circ_{i}\Gamgraph_p=(-1)^{|\Gamgraph_n||\Gamgraph_p|}(\Gamgraph_m\circ_{i-n+1}\Gamgraph_p)\circ_j\Gamgraph_n \text{ for all }i\geqslant j+n\, . 
\end{eqnarray}
\end{itemize}
Finally, the partial composition operations are equivariant with respect to the right-action $\Sigma_N$ of $\mathbb S_N$ on $\Gra_d(N)$.

The previous properties ensure that the $\mathbb S$-module $\pset{\Gra_\dime(N)}_{N\geq1}$ is naturally endowed with a structure of operad \cite{Willwacher2015}:
\begin{Proposition}[Operad $\Gra_\dime$]
For all $d\in\mathbb N$, one can define an operad in the category of graded vector spaces as the quadruplet $\big(\Gra_\dime,\Sigma,\circ_i,\Id\big)$ where:
\begin{itemize}
\item The set of graded vector spaces $\pset{\Gra_\dime(N)}_{N\geq1}$ endowed with the set of natural right-actions $\Sigma_N:\Gra_\dime(N)\times \mS_N\to \Gra_\dime(N)$ is a $\mathbb S$-module. 
\item $\circ_i:\Gra_\dime(M)\otimes\Gra_\dime(N)\to\Gra_\dime(M+N-1)$ is the set of equivariant partial composition operations defined in eq.\eqref{parcompGra}.
\item The identity element $\Id\in\Gra_\dime(1)$ is defined as the graph $\Id:=
\vcenter{\hbox{\begin{tikzpicture}[scale=0.5, >=stealth']
\node [ext] (b1) at (0,0) {\rm 1};
\end{tikzpicture}}}$ of degree $0$.
\end{itemize}
\end{Proposition}
As usual, representations of the graded operad $\Gra_d$ (or $\Gra_d$-algebras) are ordered pairs $(V,\rho)$ where $V$ is a graded vector space and $\rho:\Gra_d\to\End_V$ is a morphism of operads, with $\End_V$ the endomorphism operad on $V$, see \cite{Loday2012} for details.

\subsection{Stable structures}
\label{Stable structures}

The notion of universal structures was first introduced in \cite{Kontsevich1997} to characterise a subclass of cochains in the Chevalley--Eilenberg algebra of polyvector fields $\CE(\Tpoly)$. The terminology referred to the fact that such cochains are defined ``graphically'' via grafting of existing structures on $\Tpoly$ without resorting to additional data and thus independently of the dimension of the underlying manifold. Such universal cochains were then argued to constitute natural candidate recipients for the possible obstructions to the existence of a formality morphism. The corresponding class of formality morphisms was then precisely defined in \cite{Dolgushev2011} in terms of the operads $\mathsf{OC}$ and $\mathsf{KGra}$ (\cf definitions therein). Informally, these are \Lieinf quasi-isomorphisms whose Taylor coefficients can be written as a sum over Kontsevich admissible graphs \cite{Kontsevich:1997vb}, independently of the dimension\footnote{Two such morphisms thus only differ by their weight function, the latter depending on the choice of a Drinfel'd associator. }. More generally, universal structures can be loosely defined as originating from ``graph operads''. In the present work, we will focus\footnote{We refer to \cite{Shoikhet2008a} and \cite{Morand2021} for examples of universal structures induced from (multi)-oriented graphs on infinite-dimensional polyvector fields and (quasi)-Lie bialgebroids, respectively. } on {\it stable} structures originating from the Kontsevich's operads $\Gra_{d}$ and $\dGra_{d}$ of (un)directed graphs. The definition of stable structures adopted in the present work is adapted from \cite[Definition 4.4.1.4]{Almthesis}:
\begin{Definition}[Stable structure]
\label{stability}
Let $\mP$ be an operad in the category of graded vector spaces and $V$ a graded vector space. 
A $\mP$-algebra structure on $V$ will be said stable if the action of the operad $\mP$ on $V$ factors through the Kontsevich's graph operad $\Gra_{d}$ (or its directed avatar $\dGra_d$)  as $\mP\longrightarrow\mathsf{(d)Gra}_d\longrightarrow \End_{V}$, for some $d\in\mathbb N$. 
\end{Definition}

As will be recalled in Section \ref{subsection:Stable structures on}, an important example of stable structures on symplectic Lie $n$-algebroids is given by the notion of $\mGerd$-algebras, namely a triplet $\big(\alg,\w,\brdot\big)$ such that:
\begin{enumerate}
\item $\big(\alg,\w\big)$ is a $\mCom$-algebra.
\item $\big(\alg[d-1],\brdot\big)$ is a $\mLie$-algebra.
\item The bracket $\brdot$ is a bi-derivation with respect to the product $\w$.\footnote{Note that, in order to explicitly state the third compatibility relation, one needs first to pullback one of the defining maps along the suspension map $s:\alg[d-1]\to\alg$ of degree $d-1$ so that both products act on the same space. Explicitly, one can define the pushforward $\pbdot$ of the graded Lie bracket $\brdot$ on $\alg$ as $\pbdot:=s\, \circ\, \brdot\circ(s\un\otimes s\un)$ so that $\pb{a}{b}=(-1)^{(d-1)a}\, s \circ\br{s\un(a)}{s\un(b)}$ for all $a,b\in\alg$. The pushforward bracket $\pbdot$ is of degree $1-d$ and satisfies the following properties:
\begin{itemize}
\item graded-(skew)symmetric \ie $\pb{a}{b}=(-1)^d(-1)^{ab}\pb{b}{a}$
\item graded-Jacobi identity \ie $\pb{\pb{a}{b}}{c}+(-1)^{a(b+c)}\pb{\pb{b}{c}}{a}+(-1)^{c(a+b)}\pb{\pb{c}{a}}{b}=0$.
\end{itemize}
The graded Poisson identity on $\alg$ thus reads $\pb{a}{b\w c}=\pb{a}{b}\w c+(-1)^{b(a+1-d)}\, b\w \pb{a}{c}$.
}
\end{enumerate}
We will denote $\mGerd$ the operad whose associated representations are $\mGerd$-algebras. The notion of $\mGer_1$-algebra identifies with the one of Poisson algebra for which both binary operations are of zero degree on $\alg$. The case $d=2$ was first introduced by M. Gerstenhaber in \cite{Gerstenhaber1963} in order to characterise the natural structure living on the Hochschild cohomology of an associative algebra. For this reason, $\mGer_2$-algebras are usually referred to as \textbf{Gerstenhaber algebras}.  Note that the definition of $\mGerd$-algebras coincides with the one of $\mathsf{e}_d$-algebras \big[\cf \eg \cite[Section 13.3.16]{Loday2012}\big] for $d\geq2$ while $\mathsf{e}_1$-algebras are conventionally chosen to be associative algebras. The following Proposition asserts that any $\Gra_d$-algebra (\ie an algebra over Kontsevich's operad of undirected graphs) is endowed with a stable structure of $\mGer_d$-algebra.

\begin{Proposition}[T. Willwacher \cite{Willwacher2015}]
\label{propembedGer}
For all $d\in\mathbb N$, there is a natural embedding of operads 
$i_d:\mGerd{\longhookrightarrow}\, \Gra_{d}$.
\end{Proposition}
Explicitly, the embedding of operads $i_d$ is defined by the following action on generators $a_1\w a_2,\pb{a_1}{a_2}\in\mGerd(2)$:
\begin{itemize}
\item $i_d(a_1\w a_2)=\Gammult$ with $\w$ the graded commutative associative product of degree $0$
\item $i_d(\pb{a_1}{a_2})=\Gambr$ with $\pbdot$ the graded Lie bracket of degree $1-d$
\end{itemize}
where $\Gammult$ and $\Gambr\in\Gra_{d}(2)$ are respectively defined as:
\begin{eqnarray}
\label{gampart}
\Gammult:=\raisebox{-0.8ex}{\hbox{
\begin{tikzpicture}[scale=0.5, >=stealth']
\tikzstyle{w}=[circle, draw, minimum size=4, inner sep=1]
\tikzstyle{b}=[circle, draw, fill, minimum size=4, inner sep=1]
\node [ext] (b4) at (2,0) {\rm 1};
\node [ext] (b5) at (3.5,0) {\rm 2};
\end{tikzpicture}}}\quad,\quad\Gambr:=\raisebox{-0.8ex}{\hbox{
\begin{tikzpicture}[scale=0.5, >=stealth']
\tikzstyle{w}=[circle, draw, minimum size=4, inner sep=1]
\tikzstyle{b}=[circle, draw, fill, minimum size=4, inner sep=1]
\node [ext] (b4) at (2,0) {\rm 1};
\node [ext] (b5) at (4,0) {\rm 2};
\draw (3,0.5) node[anchor=center] {{\small $i$}};
\draw[black,->,>=latex]  (b4) to (b5);
\end{tikzpicture}}}\, .
\end{eqnarray}
In particular, the previous embeddings provide a canonical morphism of operad $\mLie\{1-d\}\to\Gra_{d}$ so that any $\Gra_{d}$-algebra is naturally endowed with a Lie bracket of degree $1-d$. Deformation complexes for these canonical morphisms will be shown to provide the definition of graph complexes in the next section.

We conclude by pointing out that in the case $d=1$, there is a natural embedding of operads $\mAss\, {\longhookrightarrow}\, \Gra_{1}$ mapping the generator $m_2\in\mAss(2)$ (\ie the associative binary product) of the associative operad $\mAss$ to the element of $\Gra_1(2)$ being explicitly defined as the infinite sum of graphs \cite{Khoroshkin2014}:
\begin{eqnarray}
\label{Moyalgraph}
\raisebox{-2.2ex}{\hbox{\begin{tikzpicture}[scale=0.5, >=stealth']
\tikzstyle{b}=[circle, draw, fill, minimum size=2, inner sep=0.02]
\node [ext] (b2) at (0,0) {1};
\node [] (b1) at (1,0.25) {$\vdots$};
\node [ext] (b3) at (2,0) {2};
\draw[black]  (b2) to[out=55, in=125, looseness=1.] (b3);
\draw[black]  (b2) to[out=-55, in=-125, looseness=1.] (b3);
\end{tikzpicture}}}:=\sum_{j\geq0}
\, \frac{1}{j!}\, 
\raisebox{-5ex}{\hbox{\begin{tikzpicture}[scale=0.5, >=stealth']
\tikzstyle{b}=[circle, draw, fill, minimum size=2, inner sep=0.02]
\node [ext] (b2) at (0,0) {1};
\node [] (b1) at (1,0.25) {$\vdots$};
\node [] (b4) at (1,-1.2) {\footnotesize{$j$ edges}};
\node [ext] (b3) at (2,0) {2};
\draw[black]  (b2) to[out=55, in=125, looseness=1.] (b3);
\draw[black]  (b2) to[out=-55, in=-125, looseness=1.] (b3);
\end{tikzpicture}}}\, .
\end{eqnarray}
As a result, $\Gra_1$-algebras are naturally endowed with a stable associative product [see eq.\eqref{GMproduct} below for an example].
\subsection{The full graph complex $\fGC_\dime$}
\label{section:The full graph complex}
We now turn to the definition of the full graph complex, denoted $\fGC_\dime$ hereafter. The differential on $\fGC_d$ stems from
a richer structure -- namely a pre-Lie structure -- defined in terms of the graph operad $\Gra_d$ using one of the following equivalent constructions:
\begin{enumerate}
\item The pre-Lie algebra associated with the suspended operad $\Gra_d\{d\}$.
\item The convolution pre-Lie algebra $\Hom_{\mathbb S}(\mathsf{coCom},\Gra_d\{d\})$. 
\item The deformation complex of the trivial operad morphism $0:\mLie\{1-d\}\to\Gra_{d}$.
\end{enumerate}

We pass on the explicit unfolding of these definitions\footnote{We refer to \cite{Loday2012,Merkulov2007} for generic constructions and to \cite{Willwacher2015,Dolgushev2012a} for applications to the case at hand. } and merely present the final result:

\begin{Proposition}[Pre-Lie structure on $\fGC_\dime$]
\label{propprelie}
For all $\dime\in\mathbb N$, the couple $\big(\fGC_\dime,\circ\big)$ where:
\begin{itemize}
\item  The graded vector space $\fGC_\dime$ is defined as\footnote{The sign conventions used relatively to the action of the various symmetry groups are summed up in Table \ref{figgraphsym}.}:

\begin{itemize}
\item $\dime$ even: $\displaystyle\fGC_\dime:=\prod_{N\geq1}\big(\Gra_\dime(N)[\dime(1-N)]\big)^{\mS_N}$
\item $\dime$ odd: $\displaystyle\fGC_\dime:=\prod_{N\geq1}\big(\Gra_\dime(N)\otimes \sgn_N[\dime(1-N)]\big)^{\mS_N}$
\end{itemize}
where the superscript stands for taking invariants with respect to the right action of $\mS_N$ with $\sgn_N$ the 1-dimensional signature representation of $\mS_N$. The terms between brackets denote degree suspension\footnote{According to the suspension, the degree of an element $\Gamgraph\in\fGC_\dime$ with $N$ vertices and $k$ edges is given by $|\Gamgraph|=\dime(N-1)+k(1-\dime)$. }.

\item The binary operation $\circ:\fGC_\dime\otimes \fGC_\dime\to\fGC_\dime$ is of degree $0$ and defined via the formula
 \begin{eqnarray}
 \nn
 \Gamgraph\circ \Gamgraph'\ =\sum_{\sigma\in \text{\rm Sh}\un(N',N-1)}(-1)^{\dime|\sigma|}\Sigma_{N+N'-1}\big(\Gamgraph\circ_1\Gamgraph'\big|\sigma)
 \end{eqnarray}
 where $\Sigma_N:\Gra_\dime(N)\times \mS_N\to \Gra_\dime(N)$ denotes the right action defined previously while $N,N'$ stand for the number of vertices in the homogeneous graphs $\Gamgraph, \Gamgraph'$, respectively. The sum is performed over the unshuffles of type $(N',N-1)$ and $|\sigma|$ denotes the signature of the permutation $\sigma\in \mS_{N+N'-1}$. 
\end{itemize}
is a graded pre-Lie algebra \ie    for all $\Gamgraph_m\in\fGC_\dime$, the following relation holds:
\begin{eqnarray}
\nn
( \Gamgraph_1\circ\Gamgraph_2)\circ \Gamgraph_3-\Gamgraph_1\circ(\Gamgraph_2\circ \Gamgraph_3)=(-1)^{|\Gamgraph_2||\Gamgraph_3|}\big(( \Gamgraph_1\circ\Gamgraph_3)\circ \Gamgraph_2-\Gamgraph_1\circ(\Gamgraph_3\circ \Gamgraph_2)\big)\, .
 \end{eqnarray}
\end{Proposition}
\Prop{propprelie} can be reformulated as the existence of a morphism of operads $\mathsf{preLie}\to\End_{\, \fGC_d}$. Composing with the morphism of operads  $\mLie\to\mathsf{preLie}$ allows to endow $\fGC_d$ with a structure of graded Lie algebra through the commutator (graded) Lie bracket 
$\brdot$ defined as:
 \begin{eqnarray}
\nn
 \br{\Gamgraph_1}{\Gamgraph_2}=\Gamgraph_1\circ\Gamgraph_2-(-1)^{|\Gamgraph_1||\Gamgraph_2|}\Gamgraph_2\circ\Gamgraph_1\, .
 \end{eqnarray}
For all $\dime\in\mathbb N$, it can be checked that the element $\Gambr\in\fGC^1_\dime$ [\cf \eqref{gampart}] is a Maurer--Cartan element for the graded Lie algebra $(\fGC_\dime,\brdot)$ \ie $  \br{\Gambr}{\Gambr}=0$.
 This property allows to define the differential operator $\delta:=\br{\Gambr}{\cdot}$ acting through the adjoint action associated with the Maurer--Cartan element. The latter can be shown to square to zero\footnote{In retrospect, it can be checked that the choices made in Table \ref{figgraphsym} are the only ones ensuring that $\delta^2\equiv0$ \cite{Willwacher2015e}. } as well as to be a derivation of the graded Lie bracket. 
 
   We sum up the previous discussion by the following proposition:
 \begin{Proposition}
 The triplet $(\fGC_\dime,\deltabr,\brdot)$ is a dg Lie algebra\footnote{Note that the dg Lie algebra $(\fGC_\dime,\deltabr,\brdot)$ can be defined from the onset as the deformation complex of the (non-trivial) operad morphism $\mLie\{1-d\}\to\Gra_{d}$ defined in Section \ref{section:graphoperad}, \cf \cite{Willwacher2015}. }. 
\end{Proposition}
  Forgetting the Lie bracket, we refer to the couple $(\fGC_\dime,\deltabr)$ as the \textbf{full graph complex}. 

\pagebreak

\begin{table}
\begin{center}
\begin{tabular}{cccc}
    \toprule
           &$\mathbb S_2^{\times k}$&$\mathbb S_k$&$\mathbb S_N$ \\ \midrule
    $d$ even&$+$&$-$&$+$  \\
    $d$ odd&$-$&$+$&$-$  \\
    \bottomrule
  \end{tabular}
  \end{center}
  \caption{Symmetries of graphs in $\fGC_d$}
  \label{figgraphsym}
\end{table}

\noindent We conclude by displaying\footnote{As is customary, we will represent a given element of $\fGC_d$ as a linear combination of undirected graphs with black vertices since taking invariants with respect to $\mathbb S_N$ makes the vertices undistinguishable. In order to obtain an explicit element of $\fGC_d$ from such a graph, one needs to go through the following steps (\cf Figure \ref{graphex} for an example):
\begin{enumerate}
\item Choose an ordering of the edges.
\item Choose an orientation of the edges. 
\item Sum over all possible ways of assigning labels to the vertices.
\item Divide by the order of the symmetry of the given graph. 
\end{enumerate}
Note that the overall sign is left ambiguous. 
}
 distinguished examples of graphs in $\fGC_d$:
 \begin{figure}[h]
  \begin{center}
  \begin{eqnarray}
\nn
 \raisebox{-0.5ex}{\hbox{\begin{tikzpicture}[scale=0.5, >=stealth']
\tikzstyle{b}=[circle, draw, fill, minimum size=2, inner sep=0.02]
\node [b] (b1) at (-1.5,0) {1};
\node [b] (b2) at (0,0) {2};
\node [b] (b3) at (1.5,0) {3};
\draw[black]  (b2) to (b3);
\end{tikzpicture}}}
\quad =\quad
\frac{1}{3}\Big(
\raisebox{-0.8ex}{\hbox{
\begin{tikzpicture}[scale=0.5, >=stealth']
\tikzstyle{w}=[circle, draw, minimum size=4, inner sep=1]
\tikzstyle{b}=[circle, draw, fill, minimum size=4, inner sep=1]
\node [ext] (b1) at (0,0) {1};
\node [ext] (b4) at (2,0) {\rm 2};
\node [ext] (b5) at (4,0) {\rm 3};
\draw (3,0.5) node[anchor=center] {{\small $i$}};
\draw[black,->,>=latex]  (b4) to (b5);
\end{tikzpicture}}}
\quad +\quad
\raisebox{-0.8ex}{\hbox{
\begin{tikzpicture}[scale=0.5, >=stealth']
\tikzstyle{w}=[circle, draw, minimum size=4, inner sep=1]
\tikzstyle{b}=[circle, draw, fill, minimum size=4, inner sep=1]
\node [ext] (b1) at (0,0) {2};
\node [ext] (b4) at (2,0) {\rm 3};
\node [ext] (b5) at (4,0) {\rm 1};
\draw (3,0.5) node[anchor=center] {{\small $i$}};
\draw[black,->,>=latex]  (b4) to (b5);
\end{tikzpicture}}}
\quad +\quad
\raisebox{-0.8ex}{\hbox{
\begin{tikzpicture}[scale=0.5, >=stealth']
\tikzstyle{w}=[circle, draw, minimum size=4, inner sep=1]
\tikzstyle{b}=[circle, draw, fill, minimum size=4, inner sep=1]
\node [ext] (b1) at (0,0) {3};
\node [ext] (b4) at (2,0) {\rm 1};
\node [ext] (b5) at (4,0) {\rm 2};
\draw (3,0.5) node[anchor=center] {{\small $i$}};
\draw[black,->,>=latex]  (b4) to (b5);
\end{tikzpicture}}}\nn
\Big)
\end{eqnarray}
  \end{center}
  \caption{Example of graph in $\fGC_d$}
  \label{graphex}
\end{figure}
\begin{Example}
\hfill
\begin{itemize}
\item The graph $\raisebox{-0.5ex}{\hbox{\begin{tikzpicture}[scale=0.5, >=stealth']
\tikzstyle{b}=[circle, draw, fill, minimum size=2, inner sep=0.02]
\node [b] (b2) at (0,0) {2};
\node [b] (b3) at (1.5,0) {3};
\draw[black]  (b2) to (b3);
\end{tikzpicture}}}$ is a cocycle in the even and odd graph complexes.

\item The tadpole graph $\vcenter{\hbox{\begin{tikzpicture}[scale=0.5, >=stealth']
\tikzstyle{w}=[circle, draw, minimum size=4, inner sep=1]
\tikzstyle{b}=[circle, draw, fill, minimum size=4, inner sep=1]
\node [b] (b1) at (0,0) {1};
    \draw [>=latex]  (0,-0.4)arc(157:-157:-1);
\end{tikzpicture}}}$ is a cocycle in the even graph complex and a zero graph in the odd graph complex.

\item The multi-arrows graph $\raisebox{-2ex}{\hbox{\begin{tikzpicture}[scale=0.5, >=stealth']
\tikzstyle{b}=[circle, draw, fill, minimum size=2, inner sep=0.02]
\node [b] (b2) at (0,0) {2};
\node [b] (b3) at (2,0) {3};
\draw[black]  (b2) to[out=55, in=125, looseness=1.] (b3);
\draw[black]  (b2) to (b3);
\draw[black]  (b2) to[out=-55, in=-125, looseness=1.] (b3);
\end{tikzpicture}}}$ -- sometimes referred to as the ``$\Theta$-graph'' -- is a cocycle in the odd graph complex and a zero graph in the even graph complex.
\item The $\Theta$-graph cocycle can be promoted to a Maurer--Cartan element\footnote{The obstruction to the prolongation of the $\Theta$-graph to a full Maurer--Cartan element lies in $H^2(\fGCconn_{1})\simeq\corps\dl L_{3}\dr$, \cf Section \ref{section:Cohomology of the full graph complex}. Since the obstruction to the prolongation of the $\Theta$-graph at order $k\geq2$ has Betti number $k+2$, it never hits the loop graph $L_3$ of Betti number $1$. The prolongation of the $\Theta$-graph to a Maurer--Cartan element in $\fGC_1$ is thus unobstructed at all orders. } in $(\fGC_1,\deltabr,\brdot)$ as the sum of multi-arrowed graphs \cite{Khoroshkin2014}: 
\begin{eqnarray}
\label{MoyalgraphMC}
\raisebox{-2.2ex}{\hbox{\begin{tikzpicture}[scale=0.5, >=stealth']
\tikzstyle{b}=[circle, draw, fill, minimum size=2, inner sep=0.02]
\node [b] (b2) at (0,0) {1};
\node [] (b1) at (1,0.25) {$\vdots$};
\node [b] (b3) at (2,0) {2};
\draw[black]  (b2) to[out=55, in=125, looseness=1.] (b3);
\draw[black]  (b2) to[out=-55, in=-125, looseness=1.] (b3);
\end{tikzpicture}}}:=\sum_{k\geq1}
\, \frac{1}{(2\, k+1)!}\, 
\raisebox{-5ex}{\hbox{\begin{tikzpicture}[scale=0.5, >=stealth']
\tikzstyle{b}=[circle, draw, fill, minimum size=2, inner sep=0.02]
\node [b] (b2) at (0,0) {1};
\node [] (b1) at (1,0.25) {$\vdots$};
\node [] (b4) at (1.2,-1.2) {\footnotesize{$2\, k+1$ edges}};
\node [b] (b3) at (2,0) {2};
\draw[black]  (b2) to[out=55, in=125, looseness=1.] (b3);
\draw[black]  (b2) to[out=-55, in=-125, looseness=1.] (b3);
\end{tikzpicture}}}\, .
\end{eqnarray}
\end{itemize}
\end{Example}

We conclude this section by introducing the \textbf{concatenation} of two graphs into a single (disconnected) graph, and denoted
$\cup:\gra_{M,j}\otimes\gra_{N,k}\to\gra_{M+N,\, j+k}$, as in the following example:
\begin{eqnarray}
\label{figure:concatenation product}
\raisebox{-4ex}{\hbox{
\begin{tikzpicture}[scale=0.5, >=stealth']
\tikzstyle{w}=[circle, draw, minimum size=4, inner sep=1]
\tikzstyle{b}=[circle, draw, fill, minimum size=4, inner sep=1]
\node [ext] (b1) at (0,0) {1};
\node [ext] (b2) at (2,0) {2};
\node [ext] (b3) at (1,-1.73) {3};
\draw[black,->,>=latex]  (b1) to (b2);
\draw (1,0.5) node[anchor=center] {{\small $i$}};
\draw (2,-0.85) node[anchor=center] {{\small $ii$}};
\draw (-0.08,-0.9) node[anchor=center] {{\small $iii$}};
\draw[black,->,>=latex]  (b2) to (b3);
\draw[black,->,>=latex]  (b3) to (b1);
\end{tikzpicture}}}
\quad\cup\quad
\raisebox{-0.8ex}{\hbox{
\begin{tikzpicture}[scale=0.5, >=stealth']
\tikzstyle{w}=[circle, draw, minimum size=4, inner sep=1]
\tikzstyle{b}=[circle, draw, fill, minimum size=4, inner sep=1]
\node [ext] (b4) at (2,0) {1};
\node [ext] (b5) at (4,0) {2};
\draw (3,0.5) node[anchor=center] {{\small $i$}};
\draw[black,->,>=latex]  (b4) to (b5);
\end{tikzpicture}}}
\quad=\quad
\raisebox{-4ex}{\hbox{
\begin{tikzpicture}[scale=0.5, >=stealth']
\tikzstyle{w}=[circle, draw, minimum size=4, inner sep=1]
\tikzstyle{b}=[circle, draw, fill, minimum size=4, inner sep=1]
\node [ext] (b1) at (0,0) {1};
\node [ext] (b2) at (2,0) {2};
\node [ext] (b3) at (1,-1.73) {3};
\draw[black,->,>=latex]  (b1) to (b2);
\draw (1,0.5) node[anchor=center] {{\small $i$}};
\draw (2,-0.85) node[anchor=center] {{\small $ii$}};
\draw (-0.08,-0.9) node[anchor=center] {{\small $iii$}};
\draw[black,->,>=latex]  (b2) to (b3);
\draw[black,->,>=latex]  (b3) to (b1);
\end{tikzpicture}}}
\quad\quad
\raisebox{-0.8ex}{\hbox{
\begin{tikzpicture}[scale=0.5, >=stealth']
\tikzstyle{w}=[circle, draw, minimum size=4, inner sep=1]
\tikzstyle{b}=[circle, draw, fill, minimum size=4, inner sep=1]
\node [ext] (b4) at (2,0) {4};
\node [ext] (b5) at (4,0) {5};
\draw (3,0.5) node[anchor=center] {{\small $iv$}};
\draw[black,->,>=latex]  (b4) to (b5);
\end{tikzpicture}}}\, .
\end{eqnarray}
 The concatenation product can be extended by linearity to yield a product $\cup:\Gra_d(M)\otimes\Gra_d(N)\to\Gra_d(M+N)$, which can furthermore be checked to be associative and to satisfy the ``commutation'' relation:
 \begin{eqnarray}
\nn
 \Gamgraph\cup\Gamgraph'=(-1)^{kk'(1-d)}\Sigma_{M+N}(\Gamgraph'\cup\Gamgraph|\sigma)
 \end{eqnarray}
 with $\Gamgraph\in\gra_{M,k}$, $\Gamgraph'\in\gra_{N,k'}$ and where the permutation $\sigma\in\mathbb S_{M+N}$ is defined as 
\begin{eqnarray}
\nn
 \sigma:=\begin{pmatrix}
1 &&\cdots&&&M+N\\
N+1&\cdots&M+N&1&\cdots&N
\end{pmatrix}
 \end{eqnarray}
so that:
\begin{eqnarray}
\nn
 \sigma\un=\begin{pmatrix}
1 &&\cdots&&&M+N\\
M+1&\cdots&M+N&1&\cdots&M
\end{pmatrix}\, .
 \end{eqnarray}

\noindent We will denote with the same symbol $\cup$ the corresponding concatenation operation of two graphs in $\fGC_d$. The latter can be shown to be:
 \begin{enumerate}
 \item of degree $d$
 \item graded commutative \ie $\Gamgraph\cup\Gamgraph'=(-1)^{(\Gamgraph+d)(\Gamgraph'+d)}\, \Gamgraph'\cup\Gamgraph$
 \item associative \ie $(\Gamgraph\cup\Gamgraph')\cup\Gamgraph''=\Gamgraph\cup(\Gamgraph'\cup\Gamgraph'')$.
 \end{enumerate}
  Furthermore, the differential $\deltabr$ satisfies the Leibniz rule:
  \begin{eqnarray}
  \label{Leibniz}
  \deltabr(\Gamgraph\cup\Gamgraph')=\deltabr\Gamgraph\cup\Gamgraph'+(-1)^{|\Gamgraph|+d}\, \Gamgraph\cup\deltabr\Gamgraph'\, .
  \end{eqnarray}
  
  \begin{Proposition}
  \label{prop:concatenation}
The concatenation product endows $H^{-d}(\fGC_d)$ with a structure of commutative algebra.
\end{Proposition}
\begin{proof}
Let $\gamma$, $\gamma'$ be two non-trivial cocycles of degree $-d$ in $\fGC_d$. The concatenation product being of degree $d$, the concatenation $\gamma\cup\gamma'$ is of degree $-d$. Furthermore, the Leibniz rule \eqref{Leibniz} ensures that $\gamma\cup\gamma'$ is a (necessarily non-trivial) cocycle hence $\gamma\cup\gamma'\in H^{-d}(\fGC_d)$. Lastly, note that the concatenation product is associative and that its graded commutative property ensures that $\Gamgraph\cup\Gamgraph'=\Gamgraph'\cup\Gamgraph$ when restricted to graphs of degree $-d$.
\end{proof}
\subsection{Cohomology of the full graph complex}
\label{section:Cohomology of the full graph complex}
\begin{figure}[h]
  \begin{center}
$\raisebox{-2ex}{\hbox{\begin{tikzpicture}[scale=0.5, >=stealth']
\tikzstyle{w}=[circle, draw, minimum size=4, inner sep=1]
\tikzstyle{b}=[circle, draw, fill, minimum size=2, inner sep=0.02]
\node [b] (b1) at (0,0) {1};
    \draw [>=latex]  (0,0) arc(0:360:0.8);
\end{tikzpicture}}}$
\quad\quad
$\raisebox{-2ex}{\hbox{\begin{tikzpicture}[scale=0.5, >=stealth']
\tikzstyle{b}=[circle, draw, fill, minimum size=2, inner sep=0.02]
\node [b] (b2) at (0,0) {2};
\node [b] (b3) at (2,0) {3};
\draw[black]  (b2) to[out=55, in=125, looseness=1.] (b3);
\draw[black]  (b2) to[out=-55, in=-125, looseness=1.] (b3);
\end{tikzpicture}}}$
\quad\quad
$\raisebox{-3ex}{\hbox{\begin{tikzpicture}[scale=0.5, >=stealth']
\tikzstyle{w}=[circle, draw, minimum size=4, inner sep=1]
\tikzstyle{b}=[circle, draw, fill, minimum size=2, inner sep=0.02]
\node [b] (b1) at (0,0) {1};
\node [b] (b2) at (2,0) {2};
\node [b] (b3) at (1,-1.73) {3};
\draw[black]  (b1) to (b2);
\draw[black]  (b2) to (b3);
\draw[black]  (b3) to (b1);
\end{tikzpicture}}}$
\quad\quad
$\raisebox{-3ex}{\hbox{\begin{tikzpicture}[scale=0.5, >=stealth']
\tikzstyle{w}=[circle, draw, minimum size=4, inner sep=1]
\tikzstyle{b}=[circle, draw, fill, minimum size=2, inner sep=0.02]
\node [b] (b1) at (0,0) {1};
\node [b] (b2) at (2,0) {2};
\node [b] (b4) at (0,2) {3};
\node [b] (b3) at (2,2) {4};
\draw[black]  (b1) to (b2);
\draw[black]  (b2) to (b3);
\draw[black]  (b3) to (b4);
\draw[black]  (b4) to (b1);
\end{tikzpicture}}}$
\quad\quad
$\raisebox{-3ex}{\hbox{\begin{tikzpicture}[scale=0.5, >=stealth']
\tikzstyle{w}=[circle, draw, minimum size=4, inner sep=1]
\tikzstyle{b}=[circle, draw, fill, minimum size=2, inner sep=0.02]
\draw (18:1.3) node [b] (b3) {3};
\draw (90:1.3) node [b] (b2) {2};
\draw (162:1.3) node [b] (b1) {1};
\draw (234:1.3) node [b] (b5) {5};
\draw (306:1.3) node [b] (b4) {4};
\draw[black]  (b1) to (b2);
\draw[black]  (b1) to (b5);
\draw[black]  (b2) to (b3);
\draw[black]  (b3) to (b4);
\draw[black]  (b4) to (b5);
\end{tikzpicture}}}$
  \caption{Loop graphs $L_k$ for $k\in\pset{1,\ldots,5}$}
  \label{figloop}
\end{center}
\end{figure}
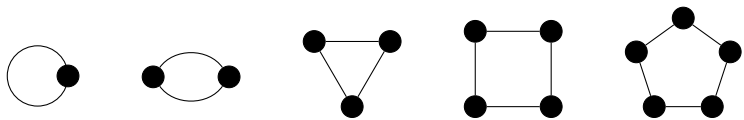

We now collect some known results regarding the cohomology of the full graph complex $\fGC_d$. In the following, we will let $\fGCconn_{d}$ denote the sub-dg Lie algebra of $\fGC_{d}$ spanned by {\it connected} graphs. Furthermore, we define $\GC_d$ as the subcomplex of $\fGCconn_{d}$ spanned by graphs without tadpoles for which all vertices have valence at least 3. The latter subcomplex was introduced [in the case $d=2$] in \cite{Kontsevich1997} and is sometimes referred to as the \textbf{Kontsevich graph complex}.
As noted in \cite{Willwacher2015}, the full graph complex can be described in terms of its connected component\footnote{For any graded vector space $V$, we will let $\widehat{S}(V)$ denote the (completed) symmetric product space of the graded vector space $V$ defined as $\displaystyle\widehat{S}(V):=\prod_{j\geq1}(V^{\otimes j})^{\mathbb S_j}$. } as $\fGC_d=\widehat{S}(\fGCconn_d[-d])[d]$.\footnote{The degree shift by $d$ reflects the degree of the concatenation product \eqref{figure:concatenation product}. } In other words, computing the cohomology of $\fGC_d$ reduces to computing the cohomology of its connected component $\fGCconn_{d}$. The latter admits the following decomposition:
\begin{Theorem}[Kontsevich \cite{Kontsevich1993,Kontsevich1994}, Willwacher \cite{Willwacher2015}]
\label{conncohomo}
The connected part of the full graph complex satisfies:
\begin{eqnarray}
\nn
H^\bullet(\fGCconn_{d})=H^\bullet(\GC_d)\quad\oplus\underset{k\geq1}{\bigoplus_{k=2d+1 \text{ \rm mod }4}}\corps[d-k]
\end{eqnarray}
where the class corresponding to $\corps[d-k]$ is represented by a loop $L_k$ with $k$ edges, \cf {\rm Figure \ref{figloop}}.
\end{Theorem}
For symmetry reasons, the only non-zero loop classes are represented by:
\begin{itemize}
\item \underline{$d$ even}\hspace{1mm} Loops $L_k$ with $k=4j+1$ edges, $j\geq0$
\item \underline{$d$ odd} \hspace{1.2mm} Loops $L_k$ with $k=4j+3$ edges, $j\geq0$.
\end{itemize}

It follows from \Prop{conncohomo} that the cohomology of $\fGCconn_d$ is located in $\GC_d$, up to some known (loop) classes. We now focus on the cohomology of $\GC_d$, for $d=2,3$ (see \eg \cite{Khoroshkin2014,Fresse2017a} for a summary and \cite{Bar-Natan,Khoroshkin2014} for computer generated tables).
\paragraph{Cohomology of $\GC_2$.}

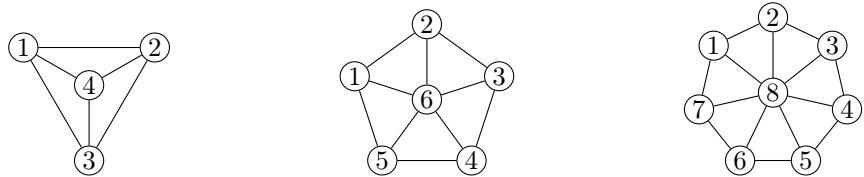
\begin{figure}[h]
  \begin{center}
\raisebox{0ex}{\hbox{\begin{tikzpicture}[scale=0.5, >=stealth']
\tikzstyle{w}=[circle, draw, minimum size=4, inner sep=1]
\tikzstyle{b}=[circle, draw, fill, minimum size=2, inner sep=0.02]
\draw (30:2) node [ext] (b2) {2};
\draw (150:2) node [ext] (b1) {1};
\draw (-90:2) node [ext] (b3) {3};
\draw (0:0) node [ext] (b4) {4};
\draw[black]  (b1) to (b2);
\draw[black]  (b1) to (b3);
\draw[black]  (b1) to (b4);
\draw[black]  (b2) to (b3);
\draw[black]  (b2) to (b4);
\draw[black]  (b3) to (b4);
\end{tikzpicture}}}
\hspace{2cm}
\raisebox{0ex}{\hbox{\begin{tikzpicture}[scale=0.5, >=stealth']
\tikzstyle{w}=[circle, draw, minimum size=4, inner sep=1]
\tikzstyle{b}=[circle, draw, fill, minimum size=2, inner sep=0.02]
\draw (18:2) node [ext] (b3) {3};
\draw (90:2) node [ext] (b2) {2};
\draw (162:2) node [ext] (b1) {1};
\draw (234:2) node [ext] (b5) {5};
\draw (306:2) node [ext] (b4) {4};
\draw (0:0) node [ext] (b6) {6};
\draw[black]  (b1) to (b6);
\draw[black]  (b2) to (b6);
\draw[black]  (b3) to (b6);
\draw[black]  (b4) to (b6);
\draw[black]  (b5) to (b6);
\draw[black]  (b1) to (b2);
\draw[black]  (b1) to (b5);
\draw[black]  (b2) to (b3);
\draw[black]  (b3) to (b4);
\draw[black]  (b4) to (b5);
\end{tikzpicture}}}
\hspace{2cm}
\raisebox{0ex}{\hbox{\begin{tikzpicture}[scale=0.5, >=stealth']
\tikzstyle{w}=[circle, draw, minimum size=4, inner sep=1]
\tikzstyle{b}=[circle, draw, fill, minimum size=2, inner sep=0.02]
\draw (38.57:2) node [ext] (b3) {3};
\draw (90:2) node [ext] (b2) {2};
\draw (141.42:2) node [ext] (b1) {1};
\draw (192.85:2) node [ext] (b7) {7};
\draw (295.71:2) node [ext] (b5) {5};
\draw (244.28:2) node [ext] (b6) {6};
\draw (347.14:2) node [ext] (b4) {4};
\draw (0:0) node [ext] (b8) {8};
\draw[black]  (b1) to (b8);
\draw[black]  (b2) to (b8);
\draw[black]  (b3) to (b8);
\draw[black]  (b4) to (b8);
\draw[black]  (b5) to (b8);
\draw[black]  (b6) to (b8);
\draw[black]  (b7) to (b8);
\draw[black]  (b8) to (b8);
\draw[black]  (b1) to (b2);
\draw[black]  (b2) to (b3);
\draw[black]  (b3) to (b4);
\draw[black]  (b4) to (b5);
\draw[black]  (b5) to (b6);
\draw[black]  (b6) to (b7);
\draw[black]  (b7) to (b1);
\end{tikzpicture}}}
  \caption{Wheel graphs for $j\in\pset{1,2,3}$}
  \label{figwheel}
\end{center}
\end{figure}

One of the major results of  \cite{Willwacher2015} is the following theorem:
\begin{Theorem}[Willwacher \cite{Willwacher2015}]
\label{thmWill}
The cohomology of the Kontsevich graph complex $\GC_2$ satisfies:
\begin{enumerate}
\item {\rm\textbf{Lower bound}}{\rm:} $H^{\leq-1}(\GC_2)=\mathbf{0}$
\item {\rm\textbf{Dominant degree}}{\rm:} $H^{0}(\GC_2)\simeq \mathfrak{grt}_1$ as Lie algebras where $\mathfrak{grt}_1$ stands for the Grothendieck--Teichm\"uller Lie algebra. 
\item {\rm\textbf{Upper bound}}{\rm:} For graphs of Betti number\footnote{A statement of the third item can be found in \cite{Khoroshkin2014} . Note that the first Betti number endows the dg Lie algebra $\fGC_\dime$ with an additional grading. It is defined explicitly as $b=k-N+c$ where $k$ denotes the number of edges, $N$ the number of vertices and $c$ the number of connected components. Relatively to the bigrading given by both $|\Gamgraph|$ and $b$, the graded Lie bracket is of bidegree $0|0$ while the differential is of bidegree $1|1$. }  $b$, the cohomology $H^{\bullet}(\GC_2)$ vanishes in degrees $\geq b-2$. 
\end{enumerate}
\end{Theorem}
Combining Theorems \ref{conncohomo} and \ref{thmWill} leads to a complete characterisation of the connected part of the full graph complex for $d=2$ in low degrees:
\begin{enumerate}
\item $H^{<-1}(\fGCconn_2)=\mathbf{0}$
\item $H^{-1}(\fGCconn_2)=\corps\dl L_1\dr$
\item $H^{0}(\fGCconn_2)=H^{0}(\GC_2)\simeq \mathfrak{grt}_1$ as Lie algebras.
\end{enumerate}
Explicit representatives of classes in the dominant degree $H^{0}(\GC_2)$ can be constructed:
\begin{Theorem}[Willwacher \cite{Willwacher2015}]
\label{propWillwheel}
For every integer $j\geq1$, there exists a non-trivial cocycle $\gamma_{2j+1}\in H^0(\GC_2)$ admitting a non-zero coefficient in front of the wheel with $2j+1$ spokes, \cf {\rm Figure \ref{figwheel}}. 
\end{Theorem}

The Grothendieck--Teichm\"uller Lie algebra $\mathfrak{grt}_1$ is known to contain a series of non-trivial elements $\sigma_3,\sigma_5,\ldots$ indexed by an odd integer\footnote{The odd elements $\pset{\sigma_{2j+1}}_{j\geq1}$ are the homogeneous components of odd degrees of the element $\psi\in\grt_1$ defined such that $g=\exp(\psi)$ is the unique element of $\GRT_1$ sending the Knizhnik--Zamolodchikov associator $\Phi_{{\mathrm{KZ}}}$ to the anti-Knizhnik--Zamolodchikov associator $\Phi_{\overline{\mathrm{KZ}}}$, \cf \eg  \cite{Rossi2014}. }. In fact, the Drinfel'd--Deligne--Ihara conjecture states that there is an isomorphism of Lie algebras between $\mathfrak{grt}_1$ and the (degree completion of) the free Lie algebra generated by the odd elements $\pset{\sigma_{2j+1}}_{j\geq1}$. Part of the conjecture has been proved by F. Brown in \cite{Brown2011} who showed that these elements generate a free Lie subalgebra of $\mathfrak{grt}_1$. In order to fully prove the conjecture, it remains to be shown that this free Lie subalgebra identifies with $\mathfrak{grt}_1$. 

 In \cite{Willwacher2015}, T. Willwacher provides an explicit isomorphism of Lie algebras $H^{0}(\GC_2)\simeq \mathfrak{grt}_1$ under which the series of odd elements $\sigma_{2j+1}\in\mathfrak{grt}_1$ gets mapped to the series of graphs $\gamma_{2j+1}$ in $H^{0}(\GC_2)$. An explicit transcendental formula for the cocycles $\gamma_{2j+1}$ is given in \cite{Rossi2014} as a sum over $\gra_{2j+2,4j+2}$ where the coefficients are provided by explicit converging integrals over the configuration space of $2n$ points in $\mathbb C\setminus \pset{0,1}$. However, a purely combinatorial construction of the $\gamma_{2j+1}$'s is still missing.

\noindent Regarding higher degrees, computer experiments have exhibited sporadic classes in $H^{\geq 3}(\GC_2)$ while it remains a difficult open conjecture (Drinfel'd, Kontsevich) that $H^1(\GC_2)=\mathbf{0}$. 

\paragraph{Cohomology of $\GC_3$.}

The cohomology of the odd graph complex can be characterised in low degrees in a way similar to the even case (see \eg \cite{Bar-Natan1995-2,Khoroshkin2014}):
\begin{enumerate}
\item \textbf{Upper bound}: $H^{\geq -2}(\GC_3)=\mathbf{0}$
\item \textbf{Dominant degree}: The dominant level of the odd graph complex $\GC_3$ is located in degree $-3$. The corresponding cohomology space $H^{-3}(\GC_3)$ can be shown to be spanned by trivalent graphs (\cf Figure \ref{figtrivalent} for examples\footnote{\label{footIHX} Trivalent graphs in the odd graph complex are usually depicted as {\it chord diagrams} where each intersection of three lines stands for a vertex. Note that modding by the IHX relation ensures that the trivalent graphs in Figure \ref{figtrivalent} satisfy the equivalence relations $A\sim2\, B$ and $C\sim4\, D\sim E\sim2\, F$. The tetrahedron graph $B$ is sometimes denoted $t$ in the literature. }) modulo the so-called {\it IHX relation} reading (see \eg  \cite{Bar-Natan1995-2}):
\begin{eqnarray}
\label{IHX}
\raisebox{-3ex}{\hbox{\begin{tikzpicture}[
my angle/.style={draw, <->, angle eccentricity=1.3, angle radius=9mm}
                        ]
\coordinate                     (O)  at (0,0);
\coordinate[]  (M) at (0.5,0);
\coordinate[]  (N) at (0,0.5);
\coordinate[]  (P) at (-0.5,0);
\coordinate[]  (Q) at (0,-0.5);
\coordinate[]  (R) at (0.5,-0.5);
\coordinate[]  (S) at (-0.5,-0.5);
\coordinate[]  (T) at (0.5,0.5);
\coordinate[]  (U) at (-0.5,0.5);
\draw[thick]    (T) -- (U) (R) -- (S) (N) -- (Q);
    \end{tikzpicture}   }}
  \quad  = \quad
    \raisebox{-3ex}{\hbox{\begin{tikzpicture}[
my angle/.style={draw, <->, angle eccentricity=1.3, angle radius=9mm}
                        ]
\coordinate                     (O)  at (0,0);
\coordinate[]  (M) at (0.5,0);
\coordinate[]  (N) at (0,0.5);
\coordinate[]  (P) at (-0.5,0);
\coordinate[]  (Q) at (0,-0.5);
\coordinate[]  (R) at (0.5,-0.5);
\coordinate[]  (S) at (-0.5,-0.5);
\coordinate[]  (T) at (0.5,0.5);
\coordinate[]  (U) at (-0.5,0.5);
\draw[thick]    (S) -- (U) (P) -- (M) (T) -- (R);
    \end{tikzpicture}   }}
 \quad   + \quad
        \raisebox{-3ex}{\hbox{\begin{tikzpicture}[
my angle/.style={draw, <->, angle eccentricity=1.3, angle radius=9mm}
                        ]
\coordinate                     (O)  at (0,0);
\coordinate[]  (M) at (0.5,0);
\coordinate[]  (N) at (0,0.5);
\coordinate[]  (P) at (-0.5,0);
\coordinate[]  (Q) at (0,-0.5);
\coordinate[]  (R) at (0.5,-0.5);
\coordinate[]  (S) at (-0.5,-0.5);
\coordinate[]  (T) at (0.5,0.5);
\coordinate[]  (U) at (-0.5,0.5);
\coordinate[]  (V) at (135:0.25);
\coordinate[]  (W) at (-135:0.25);
\draw[thick]    (U) -- (R) (T) -- (S) (V) -- (W);
    \end{tikzpicture}   }}\, .
    \end{eqnarray}
   The cohomology space $H^{-3}(\GC_3)$ is furthermore endowed with a structure of unital $\mathsf{Com}$-algebra\footnote{The commutative product is defined as follows. Let $\gamma,\gamma'$ be two trivalent graphs. Remove one (arbitrarily chosen) vertex of $\gamma$ so that the resulting graph has now three dangling edges. Repeat the previous operation for $\gamma'$, then pick one dangling edge of the graph obtained from $\gamma$ and connect it to one (arbitrarily chosen) dangling edge of the graph obtained from $\gamma'$. Repeat the operation by connecting the two remaining dangling edges of the graph obtained from $\gamma$ to the ones of $\gamma'$. The outcome is a single trivalent graph. Modding out by the IHX relation ensures that the procedure is independent of both the choice of removed vertices and pairing of dangling edges and that the resulting product is commutative. For example, one can check that $A\cdot B=F$, \cf Figure \ref{figtrivalent}.} where the r\^ole of the unit is played by the $\Theta$-graph $\raisebox{-1.1ex}{\hbox{\begin{tikzpicture}[
my angle/.style={draw, <->, angle eccentricity=1.3, angle radius=9mm},scale=0.45
                        ]
\coordinate                     (O)  at (0,0);
\coordinate[]  (A) at (180:0.6);
\coordinate[]  (B) at (0:0.6);
\draw[thick]    (A) -- (B);
\draw[thick]  (0,0) circle (0.6cm);
    \end{tikzpicture}}}$ . In fact, there is a morphism of commutative algebras:
\begin{eqnarray}
\label{mapH3}
\corps\dl t,\om_0,\om_1,\ldots\dr/ \big(\om_p\, \om_q-\om_0\, \om_{p+q},P\big)\to H^{-3}(\GC_3)
\end{eqnarray}
for a certain (explicitly known) polynomial $P$, \cf \cite{Kneissler2003a,Kneissler2003,Kneissler2003b,Vogel2011}. The map \eqref{mapH3} is conjectured to be an isomorphism, up to a 1-dimensional class represented by the $\Theta$-graph. 
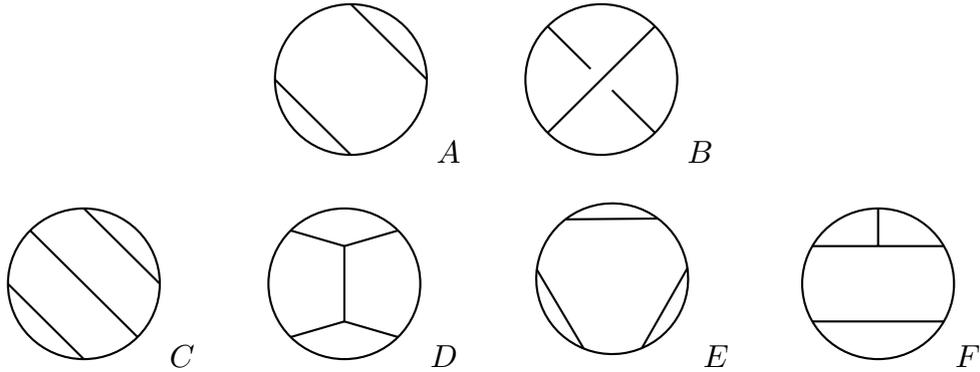
\begin{figure}[h]
  \begin{center}
    $\raisebox{-0.5ex}{\hbox{\begin{tikzpicture}[
my angle/.style={draw, <->, angle eccentricity=1.3, angle radius=9mm}
                        ]
\coordinate                     (O)  at (0,0);
\coordinate[]  (M) at (1,0);
\coordinate[]  (N) at (0,1);
\coordinate[]  (P) at (-1,0);
\coordinate[]  (Q) at (0,-1);
\coordinate[]  (R) at (135:1);
\coordinate[]  (S) at (-45:1);
\draw[thick]    (M) -- (N) (P) -- (Q) ;
\draw[thick]  (0,0) circle (1cm);
    \end{tikzpicture}   }}_{\scalebox{1.2}{$A$}}$
            \quad    \quad
                      $\raisebox{-0.5ex}{\hbox{\begin{tikzpicture}[
my angle/.style={draw, <->, angle eccentricity=1.3, angle radius=9mm}
                        ]
\coordinate                     (O)  at (0,0);
\coordinate[]  (A) at (90:1);
\coordinate[]  (B) at (45:1);
\coordinate[]  (C) at (0:1);
\coordinate[]  (D) at (-45:1);
\coordinate[]  (E) at (-90:1);
\coordinate[]  (F) at (-135:1);
\coordinate[]  (G) at (-180:1);
\coordinate[]  (H) at (135:1);
\coordinate[]  (H1) at (135:0.2);
\coordinate[]  (H2) at (135:-0.2);
\coordinate[]  (I) at (30:1);
\coordinate[]  (J) at (-60:1);
\coordinate[]  (K) at (-150:1);
\coordinate[]  (L) at (120:1);
\coordinate[]  (M) at (60:1);
\coordinate[]  (N) at (-30:1);
\coordinate[]  (P) at (150:1);
\coordinate[]  (Q) at (-120:1);
\coordinate[]  (R) at (165:1);
\coordinate[]  (S) at (-15:1);
\coordinate[]  (T) at (-75:1);
\draw[thick]    (H) -- (H1) (H2) -- (D) (F) -- (B) ;
\draw[thick]  (0,0) circle (1cm);
    \end{tikzpicture}   }}_{\scalebox{1.2}{$B$}}$
    
    \bigskip
        $\raisebox{-0.5ex}{\hbox{\begin{tikzpicture}[
my angle/.style={draw, <->, angle eccentricity=1.3, angle radius=9mm}
                        ]
\coordinate                     (O)  at (0,0);
\coordinate[]  (M) at (1,0);
\coordinate[]  (N) at (0,1);
\coordinate[]  (P) at (-1,0);
\coordinate[]  (Q) at (0,-1);
\coordinate[]  (R) at (135:1);
\coordinate[]  (S) at (-45:1);
\draw[thick]    (M) -- (N) (P) -- (Q) (R) -- (S);
\draw[thick]  (0,0) circle (1cm);
    \end{tikzpicture}   }}_{\scalebox{1.2}{$C$}}$
    \quad    \quad
                 $\raisebox{-0.5ex}{\hbox{   \begin{tikzpicture}[
my angle/.style={draw, <->, angle eccentricity=1.3, angle radius=9mm}
                        ]
\coordinate                     (O)  at (0,0);
\coordinate[]  (M) at (0,0.5);
\coordinate[]  (N) at (0,-0.5);
\coordinate[]  (R) at (135:1);
\coordinate[]  (S) at (45:1);
\coordinate[]  (P) at (-135:1);
\coordinate[]  (Q) at (-45:1);
\draw[thick]    (M) -- (N) (M) -- (R) (M) -- (S) (N) -- (P) (N) -- (Q);
\draw[thick]  (0,0) circle (1cm);
    \end{tikzpicture}    }}_{\scalebox{1.2}{$D$}}$
        \quad  \quad
                $\raisebox{-0.5ex}{\hbox{        \begin{tikzpicture}[
my angle/.style={draw, <->, angle eccentricity=1.3, angle radius=9mm}, rotate=8
                        ]
\coordinate                     (O)  at (0,0);
\coordinate[]  (A) at (90:1);
\coordinate[]  (B) at (45:1);
\coordinate[]  (C) at (0:1);
\coordinate[]  (D) at (-45:1);
\coordinate[]  (E) at (-90:1);
\coordinate[]  (F) at (-135:1);
\coordinate[]  (G) at (-180:1);
\coordinate[]  (H) at (135:1);
\coordinate[]  (I) at (30:1);
\coordinate[]  (J) at (-60:1);
\coordinate[]  (K) at (-150:1);
\coordinate[]  (L) at (120:1);
\coordinate[]  (M) at (60:1);
\coordinate[]  (N) at (-30:1);
\coordinate[]  (P) at (150:1);
\coordinate[]  (Q) at (-120:1);
\coordinate[]  (R) at (165:1);
\coordinate[]  (S) at (-15:1);
\coordinate[]  (T) at (-75:1);
\draw[thick]    (B) -- (L) (R) -- (Q) (T) -- (C);
\draw[thick]  (0,0) circle (1cm);
    \end{tikzpicture}    }}_{\scalebox{1.2}{$E$}}$
        \quad  \quad
                     $\raisebox{-0.5ex}{\hbox{       \begin{tikzpicture}[
my angle/.style={draw, <->, angle eccentricity=1.3, angle radius=9mm}
                        ]
\coordinate                     (O)  at (0,0);
\coordinate[]  (A) at (90:1);
\coordinate[]  (B) at (45:1);
\coordinate[]  (C) at (0:1);
\coordinate[]  (D) at (-45:1);
\coordinate[]  (E) at (-90:1);
\coordinate[]  (F) at (-135:1);
\coordinate[]  (G) at (-180:1);
\coordinate[]  (H) at (135:1);
\coordinate[]  (I) at (30:1);
\coordinate[]  (J) at (-60:1);
\coordinate[]  (K) at (-150:1);
\coordinate[]  (L) at (120:1);
\coordinate[]  (M) at (60:1);
\coordinate[]  (N) at (-30:1);
\coordinate[]  (P) at (150:1);
\coordinate[]  (Q) at (-120:1);
\coordinate[]  (R) at (165:1);
\coordinate[]  (S) at (-15:1);
\coordinate[]  (T) at (-75:1);
\coordinate[]  (U) at (0,0.5);
\draw[thick]    (K) -- (N) (P) -- (I) (U) -- (A);
\draw[thick]  (0,0) circle (1cm);
    \end{tikzpicture}}}_{\scalebox{1.2}{$F$}}$
  \caption{Non-trivial connected trivalent graphs in $\GC_3$ for $N=4$ ($A,B$) and $N=6$ ($C,D,E,F$). }
  \label{figtrivalent}
\end{center}
\end{figure}
\item \textbf{Lower bound}: For graphs of Betti number $b$, the cohomology $H^{\bullet}(\GC_3)$ vanishes in degrees $\leq -b-2$. 
\end{enumerate}
Regarding lower degrees, computer experiments have shown that there exist sporadic classes in $H^{-6}(\GC_3)$.

\subsection{The directed graph complex}
\label{The directed graph complex}
We conclude this review of graph complexes by presenting an important variant of the full graph complex known as the \textbf{full directed graph complex} $\dfGC_d$. Following similar steps as for $\fGC_d$, we start by defining, for all $N\geq1$, the graded vector space $\dGra_\dime(N)$ as:
\begin{itemize}
\item $\dime$ even: \hspace{0.15cm}$\displaystyle\dGra_\dime(N):=\prod_{k\geq0}\big(\corps\dl\gra_{N,k}\dr\otimes_{\mS_k}\sgn_k\big)[k(\dime-1)]$
\item $\dime$ odd: \hspace{0.27cm}$\displaystyle\dGra_\dime(N):=\prod_{k\geq0}\big(\corps\dl\gra_{N,k}\dr_{\mS_k}\big)[k(\dime-1)]$
\end{itemize}
where the subscript stands for taking coinvariants with respect to the diagonal right action of ${\mS_k}$ and the term between brackets denotes degree suspension.

In other words, the definition of $\dGra_\dime(N)$ differs from the one of $\Gra_\dime(N)$ by relaxing the modding out by $\mathbb S_2^{\otimes k}$. 
As a result, we deal with {\it directed graphs} \ie whose edge orientation is fixed. 
Similarly to the undirected case, the set of graded vector spaces $\pset{\dGra_\dime(N)}_{N\geq1}$ assemble to an operad $\dGra_\dime$. There is an injective morphism of operads 
\begin{eqnarray}
\label{morphoperad}
\Or:\Gra_d\hookrightarrow\dGra_d
\end{eqnarray}
called the \textbf{orientation morphism} and defined by sending each {\it undirected} graph into a sum of {\it directed} graphs obtained by interpreting each undirected edge as a sum of directed edges in both directions, \cf Figure \ref{figormorph}. 

Similarly to its undirected counterpart, the deformation complex $\dfGC_d$ of the trivial operad morphism $0:\mLie\{1-d\}\to\dGra_{d}$ is endowed with a pre-Lie structure [\cf Proposition \ref{propprelie}]. The latter can in turn be extended to a dg Lie algebra structure
where the differential is induced by the Maurer--Cartan element:
\begin{eqnarray}
\raisebox{-0.5ex}{\hbox{\begin{tikzpicture}[scale=0.5, >=stealth']
\tikzstyle{b}=[circle, draw, fill, minimum size=2, inner sep=0.02]
\node [b] (b2) at (0,0) {2};
\node [b] (b3) at (2,0) {3};
\draw[black,->,>=latex]   (b2) to (b3);
\end{tikzpicture}}}:=\raisebox{-0.7ex}{\hbox{\begin{tikzpicture}[scale=0.5, >=stealth']
\tikzstyle{b}=[circle, draw, fill, minimum size=2, inner sep=0.02]
\node [ext] (b2) at (0,0) {1};
\node [ext] (b3) at (2,0) {2};
\draw[black,->,>=latex]   (b2) to (b3);
\end{tikzpicture}}}+(-1)^d\, 
\raisebox{-0.7ex}{\hbox{\begin{tikzpicture}[scale=0.5, >=stealth']
\tikzstyle{b}=[circle, draw, fill, minimum size=2, inner sep=0.02]
\node [ext] (b2) at (0,0) {2};
\node [ext] (b3) at (2,0) {1};
\draw[black,->,>=latex]   (b2) to (b3);
\end{tikzpicture}}}\ .\nn
\end{eqnarray} 

We will pursue with the previously introduced notation and denote $\dfGCconn_d$ the sub-dg Lie algebra spanned by connected graphs. The morphism of operads \eqref{morphoperad} induces a morphism of dg Lie algebras:
\begin{eqnarray}
\label{morphdg Lie algebra}
\tOr:\fGCconn_d\hookrightarrow\dfGCconn_d\, .
\end{eqnarray}
\begin{figure}[h]
  \begin{center}
$\tOr\big(\raisebox{-5ex}{\hbox{\begin{tikzpicture}[scale=0.5, >=stealth']
\tikzstyle{w}=[circle, draw, minimum size=4, inner sep=1]
\tikzstyle{b}=[circle, draw, fill, minimum size=2, inner sep=0.02]
\node [b] (b2) at (30:2) {2};
\node [b] (b1) at (150:2) {1};
\node [b] (b3) at (-90:2) {3};
\node [b] (b4) at (0,0) {4};
\draw[black]  (b1) to (b2);
\draw[black]   (b3) to (b1);
\draw[black]   (b4) to (b1);
\draw[black]   (b3) to (b2);
\draw[black]   (b4) to (b2);
\draw[black]   (b3) to (b4);
\end{tikzpicture}}}\big)=24\raisebox{-5ex}{\hbox{\begin{tikzpicture}[scale=0.5, >=stealth']
\tikzstyle{w}=[circle, draw, minimum size=4, inner sep=1]
\tikzstyle{b}=[circle, draw, fill, minimum size=2, inner sep=0.02]
\node [b] (b2) at (30:2) {2};
\node [b] (b1) at (150:2) {1};
\node [b] (b3) at (-90:2) {3};
\node [b] (b4) at (0,0) {4};
\draw[black,->,>=latex]   (b1) to (b2);
\draw[black,->,>=latex]   (b1) to (b4);
\draw[black,->,>=latex]   (b1) to (b3);
\draw[black,->,>=latex]   (b2) to (b4);
\draw[black,->,>=latex]   (b2) to (b3);
\draw[black,->,>=latex]   (b4) to (b3);
\end{tikzpicture}}}
+8\raisebox{-5ex}{\hbox{\begin{tikzpicture}[scale=0.5, >=stealth']
\tikzstyle{w}=[circle, draw, minimum size=4, inner sep=1]
\tikzstyle{b}=[circle, draw, fill, minimum size=2, inner sep=0.02]
\node [b] (b2) at (30:2) {2};
\node [b] (b1) at (150:2) {1};
\node [b] (b3) at (-90:2) {3};
\node [b] (b4) at (0,0) {4};
\draw[black,->,>=latex]   (b1) to (b2);
\draw[black,->,>=latex]   (b1) to (b4);
\draw[black,->,>=latex]   (b1) to (b3);
\draw[black,->,>=latex]   (b4) to (b2);
\draw[black,->,>=latex]   (b3) to (b4);
\draw[black,->,>=latex]   (b2) to (b3);
\end{tikzpicture}}}
+24\raisebox{-5ex}{\hbox{\begin{tikzpicture}[scale=0.5, >=stealth']
\tikzstyle{w}=[circle, draw, minimum size=4, inner sep=1]
\tikzstyle{b}=[circle, draw, fill, minimum size=2, inner sep=0.02]
\node [b] (b2) at (30:2) {2};
\node [b] (b1) at (150:2) {1};
\node [b] (b3) at (-90:2) {3};
\node [b] (b4) at (0,0) {4};
\draw[black,->,>=latex]   (b1) to (b2);
\draw[black,->,>=latex]   (b1) to (b4);
\draw[black,->,>=latex]   (b3) to (b1);
\draw[black,->,>=latex]   (b4) to (b2);
\draw[black,->,>=latex]   (b3) to (b4);
\draw[black,->,>=latex]   (b2) to (b3);
\end{tikzpicture}}}
+8\raisebox{-5ex}{\hbox{\begin{tikzpicture}[scale=0.5, >=stealth']
\tikzstyle{w}=[circle, draw, minimum size=4, inner sep=1]
\tikzstyle{b}=[circle, draw, fill, minimum size=2, inner sep=0.02]
\node [b] (b2) at (30:2) {2};
\node [b] (b1) at (150:2) {1};
\node [b] (b3) at (-90:2) {3};
\node [b] (b4) at (0,0) {4};
\draw[black,->,>=latex]   (b1) to (b2);
\draw[black,->,>=latex]   (b4) to (b1);
\draw[black,->,>=latex]   (b1) to (b3);
\draw[black,->,>=latex]   (b2) to (b4);
\draw[black,->,>=latex]   (b4) to (b3);
\draw[black,->,>=latex]   (b2) to (b3);
\end{tikzpicture}}}$
  \caption{Orientation morphism}
  \label{figormorph}
\end{center}
\end{figure}
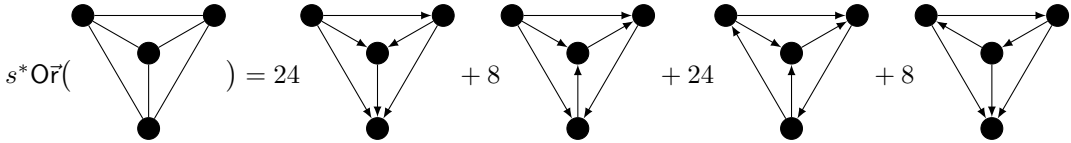

The following result was shown by T. Willwacher in \cite{Willwacher2015}, \cf also \cite{Dolgushev2017}.
\begin{Theorem}
\label{dirQI}
The morphism $\tOr:\fGCconn_d\hookrightarrow\dfGCconn_d$ is a quasi-isomorphism of dg Lie algebras.
\end{Theorem}
Theorem \ref{dirQI} implies that the study of the cohomology of the directed graph complex boils down to the one of the full graph complex, so that essentially nothing new appears when going from undirected to directed graphs. However, the directed graph complex constitutes a useful intermediary when considering representations of the Kontsevich's graph complex, \cf \cite{Buring2018}. Furthermore, the directed graph complex possesses two interesting subcomplexes spanned by {\it oriented} and {\it sourced} graphs, respectively, which have recently been shown to provide some incarnations of the Grothendieck--Teichm\"{u}ller Lie algebra $\grt_1$ in higher dimensions\footnote{That is, for values of $d>2$.}, see \cite{Willwacher2015c,Zivkovic2017a,Zivkovic2017,Merkulov2017,Merkulov2019} for details and \cite{Morand2021} for an application to representations of $\grt_1$ on quasi-Lie bialgebroids.
\section{Stable structures on graded manifolds}
\label{section:Stable structures on graded manifolds}
In the formulation of his ``{\it Formality conjecture}''  \cite{Kontsevich1997}, M. Kontsevich introduced a stable version of the deformation complex of the Schouten algebra of polyvector fields, in the guise of an injective morphism $\fGC_2\hookrightarrow\CE(\Tpoly)$. As shown in \cite{Willwacher2015,Jost2012}, this morphism of dg Lie algebras can be best understood as originating from a morphism of operads $\Gra_2\hookrightarrow\End_{\fonc{\V}}$ where $\V:=T^*[1]\M$ is a $\NP$-manifold of degree 1 whose associated graded Poisson algebra of functions is isomorphic to $\Tpoly$ endowed with the Schouten bracket (up to suspension). The aim of the present section is to generalise Kontsevich's construction from $d=2$ to arbitrary $d>0$. 

In other words, we will introduce a tower of representations $\Gra_d\hookrightarrow \End_{\fonc{\V}}$ with $\V$ an arbitrary $\NP$-manifold of degree $n$, such that $d=n+1$. This tower of morphism of operads will in turn induce a tower of injective morphisms of dg Lie algebras $\fGC_d\hookrightarrow\CE(\Tpolyn)$. Cochains in the image of this map will be called {\it stable} and results from the cohomology of $\fGC_d$ (as recalled in Section \ref{section:Cohomology of the full graph complex}) will allow to classify stable structures on $\NP$-manifolds. 

\subsection{Representations of (un)directed graphs on $\NP$-manifolds}
As for notation, we will let $(\V,\om)$ be a $\NP$-manifold of arbitrary degree $n\in\mathbb N$ with $d=n+1$ and denote $\pbdot_\om$ the associated Poisson bracket of degree $-n$. We will make use of the local presentation of $\NP$-manifolds provided in Section \ref{section:NPmanifolds}. 
By analogy with the $n=1$ case, we will denote $\Tpolyn:=\fonc{\V}[n]$ the $n$-suspension of the graded algebra of functions on $\V$ along the suspension map $s:\Tpolyn\to\fonc{\V}$ of degree $|s|=n$. We will also denote $\CE(\Tpolyn)$ the cohomological Chevalley--Eilenberg bigraded vector space\footnote{\label{footnoteNR}Letting $\alg$ be a graded vector space, the associated cohomological Chevalley--Eilenberg (bi)graded vector space (in the adjoint representation) is defined as 
\begin{eqnarray}
\nn
\CE(\alg):=\bigoplus_{n\in\mathbb Z}\CE^n(\alg)\text{ where }\CE^n(\alg):=\hspace{-0.15cm}\bigoplus_{i+j=n}\Hom^i(\alg^{\w j+1},\alg).
\end{eqnarray}
The latter is endowed with a bigrading: the first degree, denoted $i$, stems from the intrinsic degree of $\alg$ while the second (shifted) degree, denoted $j$, stems form the number of inputs. The total degree is denoted $n=i+j$. The bigraded vector space $\CE(\alg)$ carries a pre-Lie algebra structure through the \textbf{Nijenhuis--Richardson product} defined as
\begin{eqnarray}
\nn
f\circ_{\NR}g=\sum_{\sigma\in \text{\rm Sh}\un(q,p-1)}(-1)^{|\sigma|}(-1)^{(p-1)l}\, \Sigma_{p+q-1}\big(f\circ_1 g\big|\sigma)
\end{eqnarray}
for all bihomogeneous functions $f\in\Hom^{k}(\alg^{\w p},\alg)$ and $g\in\Hom^{l}(\alg^{\w q},\alg)$. 

\noindent The commutator 
\begin{eqnarray}
\label{eqNRbracket}
\br{f}{g}_\NR:=f\circ_{\NR}g-(-1)^{kl}g\circ_{\NR}f
\end{eqnarray} is a graded Lie bracket referred to as the  \textbf{Nijenhuis--Richardson bracket}. Maurer--Cartan elements thereof (\ie elements $\mathfrak m\in\CE^1(\alg)$ satisfying $\br{\mathfrak m}{\mathfrak m}_\NR=0$) are in one-to-one correspondence with \Lieinf-algebra structures on $\alg$. Given a \Lieinf-algebra $\mathfrak m$, we can define a differential $\delta_{\mathfrak m}$ on $\CE(\alg)$ as $\delta_{\mathfrak m}:=\br{\mathfrak m}{\cdot}_\NR$ (which is a derivation of $\brdot_\NR$ as a consequence of the graded Jacobi identity and squares to zero as a consequence of the graded Jacobi identity and $\br{\mathfrak m}{\mathfrak m}_\NR=0$) so that $(\CE(\alg),\delta_{\mathfrak m},\br{\cdot}{\cdot}_\NR)$ is a dg Lie algebra.
} associated with the graded vector space $\Tpolyn$ and $\brdot_{\mathsf S}$ the pullback of the Poisson bracket $\pbdot_\om$ by the suspension map $s:\Tpolyn\to\fonc{\V}$ \ie $\brdot_{\mathsf S}=s\un\circ \pbdot_\om\circ (s\otimes s)$. It can be checked that $\brdot_{\mathsf S}$ is a graded Lie bracket of degree 0, thus endowing $\Tpolyn$ with a (stable) structure of graded Lie algebra. Pursuing with the previous analogy, we will refer to $\big(\Tpolyn,\brdot_{\mathsf S}\big)$ as the $n$-\textbf{Schouten} algebra. This graded Lie structure on $\Tpolyn$ allows to endow $\CE(\Tpolyn)$ with a structure of complex via the Chevalley--Eilenberg differential $\deltaS:=\big[{\brdot_{\mathsf S}},{\cdot}\big]_\NR$ -- where $\brdot_\NR$ is the Nijenhuis--Richardson bracket  \eqref{eqNRbracket} -- associated with the Schouten bracket. 
\begin{Proposition}
\label{propmordir}
The graded algebra of functions on $\V$ is endowed with a structure of a $\dGra_{\npu}$-algebra. 
\end{Proposition}
The corresponding morphism of operads of graded vector spaces will be denoted $\dRep^{(d)}:\dGra_{\npu}\hookrightarrow \End_{\fonc{\V}}$ and defined explicitly as the sequence $\pset{\dRep^{(d)}_N}_{N\geq1}$ \hspace{-2mm}of maps $\dRep^{(d)}_N:\dGra_{\npu}(N)\otimes\fonc{\V}^{\otimes N}\to\fonc{\V}$ reading, for all $\Gamgraph\in\dGra_{\npu}(N)$:
\begin{eqnarray}
\label{Orientationmorphism}
\dRep^{(d)}_N(\Gamgraph)(f_1\otimes \cdots \otimes f_N)=\mu_N\big(\hspace{-0.3cm} \underset{(i,j)\in E_\Gamgraph}{\prod}\bar\Delta_{ij}(f_1\otimes \cdots \otimes f_N)\big)
\end{eqnarray}

 \begingroup 
\setlength{\abovedisplayskip}{10pt}
where:
\setlength{\belowdisplayskip}{10pt}
\endgroup
\begin{itemize}
\item The $f_i$'s are functions on $\V$.
\item The symbol $\mu_N$ denotes the multiplication map on $N$ elements:
\begin{eqnarray}
\nn
\mu_N&:&\fonc{\V}^{\otimes N}\to \fonc{\V}\nn\\
&:&f_1\otimes f_2\otimes \cdots \otimes f_N\mapsto f_1\cdot f_2 \cdots f_N
\end{eqnarray}
\item The product is performed over the set of edges $E_\Gamgraph$. For each edge $(i,j)\in E_\Gamgraph$ connecting vertices labeled by integers $i$ and $j$, the derivative operator $\bar\Delta_{ij}$ is defined as:
\begin{itemize}
\item $\dime$ even: \hspace{0.15cm}$\displaystyle\bar\Delta_{ij}=\frac{\p }{\p x^\mu_{(i)}}\, \frac{\p }{\p p_\mu^{(j)}}+\frac{\p }{\p \psi^{\alpha_k}_{(i)}}\frac{\p }{\p \chi_{\alpha_k}^{(j)}}$\itemnum\label{defDeltadeven}
\item $\dime$ odd: \hspace{0.27cm}$\displaystyle\bar\Delta_{ij}=\frac{\p }{\p x^\mu_{(i)}}\, \frac{\p }{\p p_\mu^{(j)}}+\frac{\p }{\p \psi^{\alpha_k}_{(i)}}\frac{\p }{\p \chi_{\alpha_k}^{(j)}}+\half\frac{\p }{\p \xi^a_{(i)}}\, \kappa^{ab}\, \frac{\p }{\p \xi^b_{(j)}}$\itemnum\label{defDeltadodd}
\end{itemize}
where the sub(super)scripts $(i)$ or $(j)$ indicate that the derivative acts on the $i$-th or $j$-th factor in the tensor product.  
\end{itemize}
\begin{proof}
The maps $\dRep^{(d)}_N$ can be checked to satisfy the three following properties:
\begin{enumerate}
\item $\dRep^{(d)}_1\big(\vcenter{\hbox{\begin{tikzpicture}[scale=0.5, >=stealth']
\node [ext] (b1) at (0,0) {\rm 1};
\end{tikzpicture}}}\big)=\Id_{\fonc{\V}}$
\item $\dRep^{(d)}_{M+N-1}(\Gamgraph\circ^{\dGra}_i\Gamgraph')=\dRep^{(d)}_M(\Gamgraph)\circ_i^{\End}\dRep^{(d)}_N(\Gamgraph')$ for all $\Gamgraph\in\dGra_{\npu}(M)$ and $\Gamgraph'\in\dGra_{\npu}(N)$ where the partial composition maps of the endomorphism operad $\End_{\fonc{\V}}$ take the form:
\begin{eqnarray}
\nn
\theta\circ^\End_i\theta'=\theta\circ\big(1^{\otimes^{i-1}}\otimes\theta'\otimes 1^{\otimes ^{M-i}}\big)
\end{eqnarray}
for all $\theta\in\Hom\big(\fonc{\V}^{\otimes M},\fonc{\V}\big)$ and  $\theta'\in\Hom\big(\fonc{\V}^{\otimes N},\fonc{\V}\big)$. 
\item $\dRep^{(d)}_N\big(\Sigma^\dGra_N(\Gamgraph|\sigma)\big)=\Sigma_N^{\End}\big(\dRep^{(d)}_N(\Gamgraph)|\sigma\big)$ where the endomorphism operad right action reads:
\begin{eqnarray}
\label{actsymgr}
\Sigma^{\End}_N(\theta|\sigma)(f_1,\ldots,f_N):=\theta(f_{\sigma\un_{(1)}},\ldots,f_{\sigma\un_{(N)}})
\end{eqnarray}
for all $f_i\in \fonc{\V}$, $\theta\in\Hom\big(\fonc{\V}^{\otimes N},\fonc{\V}\big)$ and $\sigma\in\mathbb S_N$. 
\end{enumerate} 
The three above properties ensure that the maps $\pset{\dRep^{(d)}_N}_{N\geq1}$ assemble to form a morphism of operads.
\end{proof}
\begin{Example}
\label{exaproof}
Let us exemplify the second item of the previous proof on the following partial composition of graphs for odd $d$:
\begin{eqnarray}
\label{equation:examplegraph}
\raisebox{-0.8ex}{\hbox{
\begin{tikzpicture}[scale=0.5, >=stealth']
\tikzstyle{w}=[circle, draw, minimum size=4, inner sep=1]
\tikzstyle{b}=[circle, draw, fill, minimum size=4, inner sep=1]
\node [ext] (b1) at (0,0) {1};
\node [ext] (b4) at (2,0) {\rm 2};
\draw (1,0.5) node[anchor=center] {{\small $i$}};
\draw[black,->,>=latex]  (b1) to (b4);
\end{tikzpicture}}}
\quad
\circ_2
\quad
(\hspace{-1mm}
\raisebox{-0.8ex}{\hbox{
\begin{tikzpicture}[scale=0.5, >=stealth']
\tikzstyle{w}=[circle, draw, minimum size=4, inner sep=1]
\tikzstyle{b}=[circle, draw, fill, minimum size=4, inner sep=1]
\node [ext] (b1) at (0,0) {1};
\node [ext] (b4) at (1.5,0) {\rm 2};
\end{tikzpicture}}}
)
\quad
=
\quad
\raisebox{-0.8ex}{\hbox{
\begin{tikzpicture}[scale=0.5, >=stealth']
\tikzstyle{w}=[circle, draw, minimum size=4, inner sep=1]
\tikzstyle{b}=[circle, draw, fill, minimum size=4, inner sep=1]
\node [ext] (b1) at (0,0) {1};
\node [ext] (b4) at (2,0) {\rm 2};
\node [ext] (b5) at (4,0) {\rm 3};
\draw (1,0.5) node[anchor=center] {{\small $i$}};
\draw[black,->,>=latex]  (b1) to (b4);
\end{tikzpicture}}}
\quad
+
\quad
\raisebox{-0.8ex}{\hbox{
\begin{tikzpicture}[scale=0.5, >=stealth']
\tikzstyle{w}=[circle, draw, minimum size=4, inner sep=1]
\tikzstyle{b}=[circle, draw, fill, minimum size=4, inner sep=1]
\node [ext] (b1) at (0,0) {1};
\node [ext] (b4) at (2,0) {\rm 3};
\node [ext] (b5) at (4,0) {\rm 2};
\draw (1,0.5) node[anchor=center] {{\small $i$}};
\draw[black,->,>=latex]  (b1) to (b4);
\end{tikzpicture}}}\, .
\end{eqnarray}
Applying the map \eqref{Orientationmorphism} on both graphs appearing on the right-hand side of \eqref{equation:examplegraph} yields:
\begin{multline*}
\dRep^{(d)}_2(\hspace{-1mm}\raisebox{-0.8ex}{\hbox{
\begin{tikzpicture}[scale=0.5, >=stealth']
\tikzstyle{w}=[circle, draw, minimum size=4, inner sep=1]
\tikzstyle{b}=[circle, draw, fill, minimum size=4, inner sep=1]
\node [ext] (b1) at (0,0) {1};
\node [ext] (b4) at (2,0) {\rm 2};
\node [ext] (b5) at (4,0) {\rm 3};
\draw (1,0.5) node[anchor=center] {{\small $i$}};
\draw[black,->,>=latex]  (b1) to (b4);
\end{tikzpicture}}})(f_1\otimes f_2\otimes f_3) =
\frac{\p f_1}{\p x^\mu}\, \frac{\p f_2}{\p p_\mu}f_3 \\
+ \sum_{i=1}^{\half n-1}(-1)^{ik}\frac{\p f_1}{\p \psi^{\alpha_i}}\frac{\p f_2}{\p \chi_{\alpha_i}}f_3+(-1)^{k{n}/{2}}\frac{\p f_1}{\p \xi^a}\, \kappa^{ab}\, \frac{\p f_2}{\p \xi^b}f_3
\end{multline*}
\vspace{-0.8em}
\begin{multline*}
\dRep^{(d)}_2(\hspace{-1mm}\raisebox{-0.8ex}{\hbox{
\begin{tikzpicture}[scale=0.5, >=stealth']
\tikzstyle{w}=[circle, draw, minimum size=4, inner sep=1]
\tikzstyle{b}=[circle, draw, fill, minimum size=4, inner sep=1]
\node [ext] (b1) at (0,0) {1};
\node [ext] (b4) at (2,0) {\rm 3};
\node [ext] (b5) at (4,0) {\rm 2};
\draw (1,0.5) node[anchor=center] {{\small $i$}};
\draw[black,->,>=latex]  (b1) to (b4);
\end{tikzpicture}}})(f_1\otimes f_2\otimes f_3) =
\frac{\p f_1}{\p x^\mu}\, f_2\, \frac{\p f_3}{\p p_\mu} \\
+ \sum_{i=1}^{\half n-1}(-1)^{i(k+l)}\frac{\p f_1}{\p \psi^{\alpha_i}}f_2\frac{\p f_3}{\p \chi_{\alpha_i}}+(-1)^{{n}/{2}(k+l)}\frac{\p f_1}{\p \xi^a}\, f_2\, \frac{\p f_3}{\p \xi_a}\nn
\end{multline*}
for all homogeneous functions $f_1\in \foncg{k}$, $f_2\in \foncg{l}$ and $f_3\in \foncg{m}$ of degree $k$, $l$ and $m$ respectively. 
By comparison, applying the map \eqref{Orientationmorphism} on both graphs appearing on the left-hand side of \eqref{equation:examplegraph} yields:
\begin{eqnarray}
\dRep^{(d)}_2(\hspace{-1mm}\raisebox{-0.8ex}{\hbox{
\begin{tikzpicture}[scale=0.5, >=stealth']
\tikzstyle{w}=[circle, draw, minimum size=4, inner sep=1]
\tikzstyle{b}=[circle, draw, fill, minimum size=4, inner sep=1]
\node [ext] (b1) at (0,0) {1};
\node [ext] (b4) at (2,0) {\rm 2};
\draw (1,0.5) node[anchor=center] {{\small $i$}};
\draw[black,->,>=latex]  (b1) to (b4);
\end{tikzpicture}}})(f_1\otimes f_2)&=&\frac{\p f_1}{\p x^\mu}\, \frac{\p f_2}{\p p_\mu}+\sum_{i=1}^{\half n-1}(-1)^{ik}\frac{\p f_1}{\p \psi^{\alpha_i}}\frac{\p f_2}{\p \chi_{\alpha_i}}+(-1)^{k{n}/{2}}\frac{\p f_1}{\p \xi^a}\, \kappa^{ab}\, \frac{\p f_2}{\p \xi^b}\nn\\
\dRep^{(d)}_2(\hspace{-1mm}\raisebox{-0.8ex}{\hbox{
\begin{tikzpicture}[scale=0.5, >=stealth']
\tikzstyle{w}=[circle, draw, minimum size=4, inner sep=1]
\tikzstyle{b}=[circle, draw, fill, minimum size=4, inner sep=1]
\node [ext] (b1) at (0,0) {1};
\node [ext] (b4) at (1.5,0) {\rm 2};
\end{tikzpicture}}})(f_1\otimes f_2)&=&f_1\cdot f_2\, .\nn
\end{eqnarray}
Partial composition of the corresponding operators in the endomorphism operad $\End_{\fonc{\V}}$ yields:
\begin{align*}
\dRep^{(d)}_2&(\hspace{-1mm}\raisebox{-0.8ex}{\hbox{
\begin{tikzpicture}[scale=0.5, >=stealth']
\tikzstyle{w}=[circle, draw, minimum size=4, inner sep=1]
\tikzstyle{b}=[circle, draw, fill, minimum size=4, inner sep=1]
\node [ext] (b1) at (0,0) {1};
\node [ext] (b4) at (2,0) {\rm 2};
\draw (1,0.5) node[anchor=center] {{\small $i$}};
\draw[black,->,>=latex]  (b1) to (b4);
\end{tikzpicture}}})
\circ_i^{\End}
\dRep^{(d)}_2(\hspace{-1mm}\raisebox{-0.8ex}{\hbox{
\begin{tikzpicture}[scale=0.5, >=stealth']
\tikzstyle{w}=[circle, draw, minimum size=4, inner sep=1]
\tikzstyle{b}=[circle, draw, fill, minimum size=4, inner sep=1]
\node [ext] (b1) at (0,0) {1};
\node [ext] (b4) at (1.5,0) {\rm 2};
\end{tikzpicture}}})(f_1\otimes f_2\otimes f_3) \\
&= \,\, \dRep^{(d)}_2(\hspace{-1mm}\raisebox{-0.8ex}{\hbox{
\begin{tikzpicture}[scale=0.5, >=stealth']
\tikzstyle{w}=[circle, draw, minimum size=4, inner sep=1]
\tikzstyle{b}=[circle, draw, fill, minimum size=4, inner sep=1]
\node [ext] (b1) at (0,0) {1};
\node [ext] (b4) at (2,0) {\rm 2};
\draw (1,0.5) node[anchor=center] {{\small $i$}};
\draw[black,->,>=latex]  (b1) to (b4);
\end{tikzpicture}}})(f_1\otimes \dRep^{(d)}_2(\hspace{-1mm}\raisebox{-0.8ex}{\hbox{
\begin{tikzpicture}[scale=0.5, >=stealth']
\tikzstyle{w}=[circle, draw, minimum size=4, inner sep=1]
\tikzstyle{b}=[circle, draw, fill, minimum size=4, inner sep=1]
\node [ext] (b1) at (0,0) {1};
\node [ext] (b4) at (1.5,0) {\rm 2};
\end{tikzpicture}}})(f_2\otimes f_3)) \\
&= \,\, \frac{\p f_1}{\p x^\mu}\, \frac{\p (f_2\cdot f_3)}{\p p_\mu}\quad+\quad\sum_{i=1}^{\half n-1}(-1)^{ik}\frac{\p f_1}{\p \psi^{\alpha_i}}\frac{\p (f_2\cdot f_3)}{\p \chi_{\alpha_i}}\quad+\quad(-1)^{k{n}/{2}}\frac{\p f_1}{\p \xi^a}\, \kappa^{ab}\, \frac{\p (f_2\cdot f_3)}{\p \xi^b}
\end{align*}
Applying the graded Leibniz rule\footnote{Letting $\psi$ be a coordinate of degree $i$ and $f\in \foncg{k}$, $g\in \foncg{l}$ be homogeneous functions of degree $k$ and $l$ respectively, the graded Leibniz rule reads:
\begin{eqnarray}
\frac{\p(f\cdot g)}{\p\psi}=\frac{\p f}{\p\psi}\cdot g+(-1)^{ik}f\cdot \frac{\p g}{\p\psi}.\nn
\end{eqnarray}} yields the equality. 
\end{Example}
Composing the representation morphism $\dRep^{(d)}:\dGra_{\npu}\hookrightarrow \End_{\fonc{\V}}$ with the orientation morphism \eqref{morphoperad} endows the algebra of functions ${\fonc{\V}}$ with a structure of $\Gra_{\npu}$-algebra
through the morphism $\Rep^{(d)}:\Gra_{\npu}\overset{\Or}{\longhookrightarrow}\dGra_{\npu}\overset{\dRep^{(d)}}{\longrightarrow} \End_{\fonc{\V}}$.

The maps $\pset{\Rep^{(d)}_N}_{N\geq1}$ can be defined explicitly in a form similar to \eqref{Orientationmorphism} as:
\begin{eqnarray}
\label{Orientationmorphism2}
\Rep^{(d)}_N(\Gamgraph)(f_1\otimes \cdots \otimes f_N)=\mu_N\big(\hspace{-0.3cm} \underset{(i,j)\in E_\Gamgraph}{\prod}\Delta_{ij}(f_1\otimes \cdots \otimes f_N)\big)
\end{eqnarray}
where one traded the $\bar\Delta$ operators with:
\begin{eqnarray}
\Delta_{ij} &=& \frac{\p }{\p x^\mu_{(i)}}\, \frac{\p }{\p p_\mu^{(j)}}+\frac{\p }{\p p_\mu^{(i)}}\, \frac{\p }{\p x^\mu_{(j)}}+\frac{\p }{\p \psi^{\alpha_k}_{(i)}}\frac{\p }{\p \chi_{\alpha_k}^{(j)}}+\frac{\p }{\p \chi_{\alpha_k}^{(i)}}\frac{\p }{\p \psi^{\alpha_k}_{(j)}}\label{defDeltadeven2} \\
\Delta_{ij} &=& \frac{\p }{\p x^\mu_{(i)}}\, \frac{\p }{\p p_\mu^{(j)}}-\frac{\p }{\p p_\mu^{(i)}}\, \frac{\p }{\p x^\mu_{(j)}}+\frac{\p }{\p \psi^{\alpha_k}_{(i)}}\frac{\p }{\p \chi_{\alpha_k}^{(j)}}-(-1)^k\frac{\p }{\p \chi_{\alpha_k}^{(i)}}\frac{\p }{\p \psi^{\alpha_k}_{(j)}}+\frac{\p }{\p \xi^a_{(i)}}\, \kappa^{ab}\, \frac{\p }{\p \xi^b_{(j)}}\label{defDeltadodd2}
\end{eqnarray}
where $\dime$ is even or odd respectively.
The differential operator $\Delta_{ij}$ enjoys the following properties\footnote{In contrast, the differential operator $\bar\Delta_{ij}$ only satisfies properties 1-2 consistently with the definition of $\dGra_d$ which omits to mod out by $\mathbb S_2^{\otimes k}$, \cf Section \ref{The directed graph complex}.}:
\begin{enumerate}
\item $|\Delta_{ij}|=1-d$ consistently with the grading of an edge in $\Gra_{d}$.
\item $\Delta_{ij}\, \Delta_{kl}=-(-1)^d\Delta_{kl}\, \Delta_{ij}$ consistently with the fact that permuting two edges in graphs in $\Gra_{d}$ brings a sign only for even $d$.
\item $\Delta_{ij}=(-1)^d\Delta_{ji}$ consistently with the fact that flipping the orientation of an edge in graphs in $\Gra_{d}$ brings a sign only for odd $d$.
\end{enumerate}

The tower of morphisms $\Rep^{(d)}:\Gra_{\npu}\hookrightarrow \End_{\fonc{\V}}$ generalises to all $d$ the Kontsevich morphism for $d=2$, \cf \cite{Kontsevich1997, Willwacher2015,Jost2012} and more recently \cite{Buring2017b,Buring2018,Rutten2018,Kiselev2019}. The morphism $\Rep^{(d)}$ preserves the concatenation product as:
\begin{multline}\label{property}
\Rep^{(d)}_{N+N'}(\Gamgraph\cup\Gamgraph')(f_1\otimes\cdots \otimes f_{N+N'}) = \\ \Rep^{(d)}_N(\Gamgraph)(f_1\otimes\cdots \otimes f_{N})\cdot\Rep^{(d)}_{N'}(\Gamgraph')(f_{N+1}\otimes\cdots \otimes f_{N+N'})
\end{multline}
where the product on the right-hand side is the pointwise product of functions on $\V$.
\subsection{Stable structures on $\fonc{\V}$}
\label{subsection:Stable structures on}
The aim of the present section is to make use of the morphism $\Rep^{(d)}$ in order to define stable structures \big(in the sense of \Defi{stability}\big) on the algebra of functions $\fonc{\V}$. In particular, it was noted earlier (\cf Section \ref{section:NPmanifolds}) that the algebra of functions on $\V$ was naturally endowed with a structure of $\mGer_{d}$-algebra. This statement can be refined as follows:
\begin{Proposition}\label{propGer}The $\mGer_{d}$-algebra structure on $\fonc{\V}$ is stable \ie the action of the operad $\mGer_{d}$ factors through:
\begin{eqnarray}
\nn
\mGer_{d}\overset{i_{d}}{\longhookrightarrow} \Gra_{d}\overset{\Rep^{(d)}}{\longrightarrow} \End_{\fonc{\V}}
\end{eqnarray}
where $i_{d}:\mGer_{d}{\longhookrightarrow} \Gra_{d}$ is the natural embedding of operads defined in {\rm \Prop{propembedGer}}.
\end{Proposition}
\begin{proof}
The statement follows straightforwardly from:
\begin{eqnarray}
\nn
 \Rep^{(d)}_2(\Gammult)(f\otimes g)=f\cdot g\quad,\quad \Rep^{(d)}_2(\Gambr)(f\otimes g)=\pb{f}{g}_\om\, .
\end{eqnarray}
\end{proof}
In the case $d=2$, \Prop{propGer} can be completed\footnote{For $d=2$, T.~Willwacher showed that $\Tpoly$ (together with a choice of Poisson bivector) is in fact endowed with a much larger stable structure, namely a $\mathsf{Br}_\infty$-algebra structure of homotopy braces \cite{Willwacher2011}.  } by stating that the graded algebra of functions $\fonc{\V}$ (isomorphic to the graded algebra of polyvector fields $\Tpoly$ on $\M$) is locally\footnote{This local definition can be made global whenever the underlying manifold $\M$ possesses a volume form. Letting $\M$ be a manifold of dimension $n$, a volume form on $\M$ provides an isomorphism $i:\Tpolyd{\bullet}\to\Om^{n-\bullet}(\M)$ between polyvector fields and differential forms on $\M$. The divergence operator $\Delta=i^{-1}\circ d_\mathsf{dR}\circ i$ defined as the pullback of the de Rham differential along $i$ allows to upgrade the $\mGer_2$-algebra structure on polyvector fields of $\M$ to that of a $\mathsf{BV}$-algebra. In a local chart, the expression of $\Delta$ coincides with the one arising from the tadpole graph due to the fact that a volume form has no local structure. An extension of Kontsevich's formality morphism relating the $\mathsf{BV}$-algebra of polyvector fields and the $\mathsf{BV}_\infty$-algebra of multidifferential operators on manifolds endowed with a volume form has been worked out in \cite{Campos2015}.} endowed with a stable structure of $\mathsf{BV}$-algebra\footnote{A \textbf{Batalin--Vilkovisky algebra} (or $\mathsf{BV}$-algebra for short) is a graded commutative algebra $\big(\alg,\w)$ endowed with a unary operator $\Delta:\alg\to\alg$ satisfying:
\begin{enumerate}
\item $\Delta$ is of degree $-1$
\item $\Delta^2=0$
\item $\Delta(a\w b\w c)-\Delta(a\w b)\w c+\Delta a \w b\w c-(-1)^{|a|}\, a\w \Delta(b\w c)-(-1)^{|b|(|a|-1)}b\w \Delta(a\w c)+(-1)^{|a|}a\w\Delta b\w c+(-1)^{|a|+|b|}\, a\w b\w \Delta c=0$.
\end{enumerate}
A $\mathsf{BV}$-algebra is in particular a $\mGer_2$-algebra where the graded Lie bracket is defined as the obstruction for $\Delta$ to be a derivation of $\w$ \ie $\pb{a}{b}=\Delta(a\w b)-\Delta a\w b-(-1)^{|a|}a\w \Delta b$. The induced bracket can be shown to satisfy the axioms of a $\mGer_2$-algebra as a consequence of the axioms of $\Delta$. Alternatively, a $\mathsf{BV}$-algebra can be defined as a $\mGer_2$-algebra endowed with a unary operator $\Delta$ of degree $-1$ and satisfying $\Delta(a\w b)-\Delta a\w b-(-1)^{|a|}a\w \Delta b=\pb{a}{b}$.
}, \cf \cite{Merkulov2014}. The $\mathsf{BV}$ Laplacian $\Delta$ is then defined as the image of the tadpole:
\begin{eqnarray}
L_1=\vcenter{\hbox{\begin{tikzpicture}[scale=0.5, >=stealth']
\tikzstyle{w}=[circle, draw, minimum size=4, inner sep=1]
\tikzstyle{b}=[circle, draw, fill, minimum size=4, inner sep=1]
\node [ext] (b1) at (0,0) {1};
    \draw [>=latex]  (0,-0.4)arc(157:-157:-1);
\end{tikzpicture}}}\quad\text{ \ie }\quad\Delta:=\Rep^{(2)}_1\big(L_1\big)\, .
\end{eqnarray} 

In the case $d=1$, the classification recalled in Section \ref{section:NPQmanifolds} ensures that a $\NP$-manifold of degree 0 is in fact a (non-graded) symplectic manifold $(\M,\kappa)$. In that case, the chain of morphisms of operads $\mAss\overset{i_1}{\longhookrightarrow} \Gra_{1}\overset{\Rep^{(1)}}{\longrightarrow} \End_{\fonc{\M}}$ endows the algebra of functions on the symplectic manifold $\M$ with a stable associative structure:
 \begin{eqnarray}
 \nn
 f*_\mathsf{GM}g=\Rep^{(1)}_2(\raisebox{-2.2ex}{\hbox{\begin{tikzpicture}[scale=0.5, >=stealth']
\tikzstyle{b}=[circle, draw, fill, minimum size=2, inner sep=0.02]
\node [ext] (b2) at (0,0) {1};
\node [] (b1) at (1,0.25) {$\vdots$};
\node [ext] (b3) at (2,0) {2};
\draw[black]  (b2) to[out=55, in=125, looseness=1.] (b3);
\draw[black]  (b2) to[out=-55, in=-125, looseness=1.] (b3);
\end{tikzpicture}}})(f\otimes g)\, ,
\end{eqnarray}
where the graph $\raisebox{-2.2ex}{\hbox{\begin{tikzpicture}[scale=0.5, >=stealth']
\tikzstyle{b}=[circle, draw, fill, minimum size=2, inner sep=0.02]
\node [ext] (b2) at (0,0) {1};
\node [] (b1) at (1,0.25) {$\vdots$};
\node [ext] (b3) at (2,0) {2};
\draw[black]  (b2) to[out=55, in=125, looseness=1.] (b3);
\draw[black]  (b2) to[out=-55, in=-125, looseness=1.] (b3);
\end{tikzpicture}}}\in\Gra_1(2)$ is defined in \eqref{Moyalgraph}, \cf \cite{Khoroshkin2014}. The induced associative product is the \textbf{Groenewold--Moyal} product  \cite{Groenewold1946,Moyal1949} reading explicitly as:
\begin{eqnarray}
\label{GMproduct}
(f*_\mathsf{GM}g)(\xi):=\exp\Big(\epsilon\, \kappa^{ab}\frac{\p}{\p{\zeta^a}}\frac{\p}{\p{\eta^b}}\Big) f({\zeta})\, g(\eta)\Big|_{\zeta=\eta=\xi} 
\end{eqnarray}
where $\epsilon$ is a formal parameter.
\subsubsection{Stable cochains of the Chevalley--Eilenberg algebra}
The following proposition can be seen as a corollary of our main \Prop{propmordir} and defines a stable graph model for $\CE(\Tpolyn)$ which generalises Kontsevich's model \eqref{morphismfGCCEn} in $d=2$ to arbitrary $d\geq1$:

\begin{Proposition}
\label{propmorphism}
The morphism of operads $\Rep^{(d)}:\Gra_{d}\hookrightarrow \End_{\fonc{\V}}$ induces a morphism of dg Lie algebras:
\begin{eqnarray}
\label{univmordg Lie algebra}
\sRep^{(d)}:\big(\fGC_{d},\deltabr,\brdot\big)\hookrightarrow \big(\CE(\Tpolyn),\deltaS,\brdot_{\NR}\big)\, .
\end{eqnarray}
\end{Proposition}
\begin{proof}
The proof follows straightforwardly from the equivariance \eqref{actsymgr} of the morphism $\Rep^{(d)}$ and from the equality $\sRep^{(d)}(\Gambr)=\brdot_{\mathsf S}$. 
\end{proof}
Pursuing with the terminology introduced in \Defi{stability}, Chevalley--Eilenberg cochains in the image of \eqref{univmordg Lie algebra} will be referred to as {\it stable}. In other words, the dg Lie algebra of graphs $\big(\fGC_{d},\deltabr,\brdot\big)$ provides a stable version of the Chevalley--Eilenberg dg Lie algebra $\big(\CE(\Tpolyn), \deltaS, $ $\brdot_{\NR}\big)$. The former thus controls the deformation theory -- in the stable setting\footnote{Note that {\it not} all conceivable deformations of $\Tpolyn$ are universal, \ie defined in terms of graphs. For example, letting $H\in\Om^3(\M)$ be a closed 3-form on a manifold $\M$, one can define a {\it non-stable} deformation of $\Tpolyd{(1)}=\Tpoly$ by defining a higher bracket of arity $3$ as $l_3\in\Hom^{-1}\big(\Tpoly^{\w3},\Tpoly\big)\in\CE^1(\Tpoly)$ as $l_3(X_1,X_2,X_3)=H({X_1},{X_2},{X_3})$, where the $X_i$'s are polyvector fields. 
       Denoting $l_2$ the usual Schouten bracket, the triplet $(\mathcal T_\poly(\M),l_2,l_3)$ forms a \Lieinf-algebra. Associated Maurer--Cartan elements are so-called \textbf{twisted Poisson structures} \cite{Klimcik2002,Severa2001} \ie bivectors $\pi\in\field{\w^2T\M}$ satisfying $\displaystyle\br{\pi}{\pi}_{\rm S}=\frac{1}{3}H(\pi,\pi,\pi)$.
       The latter can be interpreted as Dirac structures for the standard Courant algebroid on $\M$ twisted by $H$, \cf \cite{Severa2002}.
} -- of the  $n$-Schouten Lie algebra as a \Lieinf-algebra. The two following corollaries make this fact explicit:

\begin{Corollary}
\label{corMC}
Maurer--Cartan elements for the dg Lie algebra $\big(\fGC_{d},\deltabr,\brdot\big)$ are mapped via $\sRep^{(d)}$ to stable deformations of the graded Lie algebra $(\Tpolyn,\brdot_{\mathsf S})$ as a \Lieinf-algebra. 
\end{Corollary}
\begin{Example}[Groenewold--Moyal commutator]
\label{GMc}
Let $(\M,\kappa)$ be a symplectic manifold. 

\noindent The Maurer--Cartan element \eqref{MoyalgraphMC} prolongating the $\Theta$-graph is mapped via $\sRep^{(1)}$ to the (essentially unique) stable deformation of $(\fonc{\M},\pbdot_\kappa)$ as a Lie algebra\footnote{Note that $n=0$ so that the ordinary (\ie non-graded) Lie algebra $(\fonc{\M},\pbdot_\kappa)$ identifies with the $0$-Schouten algebra $(\Tpolyd{(0)},\brdot_{\mathsf S})$. In this case, the deformation complex $\big(\CE(\Tpolyn),\deltaS,\brdot_{\NR}\big)$ controls the deformation theory of $(\fonc{\M},\pbdot_\kappa)$ as an ordinary Lie algebra.} where the Poisson bracket is deformed into the \textbf{Groenewold--Moyal commutator} $\br{f}{g}_\mathsf{GM}:=f*_\mathsf{GM}g-g*_\mathsf{GM}f$ on $\fonc{\M}$ constructed from \eqref{GMproduct}, \cf \cite{Khoroshkin2014}.
\end{Example}
We refer to Example \ref{exalinf} for an example of deformation of the $3$-Schouten algebra as a genuine \Lieinf-algebra. The Lie algebra $H^0(\fGC_\dime)$ being pro-nilpotent\footnote{See \cite{Dolgushev2011} for a proof of the  $d=2$ case. The proof for all $d$ is identical.}, one defines the pro-unipotent group $\exp\big(H^0(\fGC_\dime)\big)$ as  \cite{WillwachernotesGRT}:
 \begin{itemize}
 \item Group elements are elements of $H^0(\fGC_\dime)$, viewed as a set. 
 \item The unit is $\mathbf 0\in H^0(\fGC_\dime)$.
 \item The inverse map sends $\Gamgraph$ to $-\Gamgraph$.
 \item The group operation is defined as $ \Gamgraph_1 \cdot \Gamgraph_2=\text{BCH}(\Gamgraph_1,\Gamgraph_2)$ where $\text{BCH}$ stands for the Baker--Campbell--Hausdorff formula. 
 \end{itemize}
\begin{Corollary}
\label{thmcoho}
The pro-unipotent group $\exp\big(H^0(\fGC_{d})\big)$ acts via \Lieinf-automorphisms defined up to equivalence on the $n$-Schouten algebra. Such \Lieinf-automorphisms will be referred to as stable. 
\end{Corollary}

\begin{proof}The proof is identical to the one of the case $d=2$ (\cf \eg Theorem 1. in \cite{Jost2012}, based on \cite{Willwacher2015}) that we review for completeness. 
Let $\Gamgraph\in H^0(\fGC_{d})$ be a non-trivial cocycle in $\fGC_{d}$. The morphism $\sRep^{(d)}$ of dg Lie algebras introduced in \Prop{propmorphism} maps $\Gamgraph$ to a zero degree Chevalley--Eilenberg cocycle for the $(d-1)$-Schouten algebra denoted  $\sRep^{(d)}(\Gamgraph)$. In other words, $\sRep^{(d)}(\Gamgraph)\in H^0\big(\CE(\Tpolyn)\big)$ is a \Lieinf-derivation of $\big(\Tpolyn,\brdot_{\mathsf S}\big)$. 
This ensures that $\exp\big(\sRep^{(d)}(\Gamgraph)\big)$ is a \Lieinf-automorphism of $\big(\Tpolyn,\brdot_{\mathsf S}\big)$. Furthermore, since $\exp\big(H^0(\fGC_{d})\big)$ is pro-unipotent, to any element $\Gamgraphgroup\in\exp\big(H^0(\fGC_{d})\big)$ one can associate a unique element $\Gamgraph\in H^0(\fGC_{d})$ so that $\Gamgraphgroup=\exp(\Gamgraph)$. We can thus define a \Lieinf-action via its Taylor coefficients:
\begin{eqnarray}
\mathcal U_N&:&\exp\big(H^0(\fGC_{d})\big)\times\Tpolyn{}^{\w N}\to\Tpolyn\nn\\
&:&\big(\Gamgraphgroup,X_1,\ldots,X_N\big)\mapsto\exp\big(\sRep^{(d)}(\Gamgraph)\big)(X_1,\ldots,X_N)\nn
\end{eqnarray}
for all $N\geq1$. Finally, we note that two equivalent cocycles $\gamma$ and $\gamma'$ induce homotopic \Lieinf-automorphisms.
\end{proof}

\subsubsection{Classification of stable structures}
We now make use of the results regarding cohomology of the full graph complex as reviewed in Section \ref{section:Cohomology of the full graph complex} in order to provide a classification of stable structures on the $(d-1)$-Schouten algebra for all $d\geq1$ -- where the term stable structures will refer to\footnote{We will focus on stable structures obtained from {\it connected} graphs.}:
\begin{enumerate}
\item Stable \Lieinf-automorphisms of the $(d-1)$-Schouten algebra
\item Stable deformations of the $(d-1)$-Schouten algebra as a \Lieinf-algebra.
\end{enumerate}

Since the case $d=2$  has already been addressed in the literature, we treat it separately:
\begin{enumerate}
\item Recall from Theorem \ref{thmWill} that there exists an isomorphism of Lie algebras $H^{0}(\fGCconn_2)\simeq \mathfrak{grt}_1$. By Corollary \ref{thmcoho}, it follows that the Grothendieck--Teichm\"uller group $\GRT_1:=\exp(\grt_1)$ acts  via \Lieinf-automorphisms on the Schouten algebra $\Tpoly$, see  \cite{Willwacher2015,Jost2012}. 

\item As noted earlier, it is a difficult open conjecture (Drinfel'd, Kontsevich) that $H^{1}(\fGCconn_2)=\mathbf{0}$.\footnote{As noted in \cite{Merkulov2019}, although the cohomology of $\GC_2$ in degree 1 is conjectured to be trivial, a choice of Drinfel'd associator is necessary in order to convert cocycles of degree 1 in $\GC_2$ into coboundaries of degree 0 so that an iterative procedure can exist but cannot be trivial.}
If the conjecture holds, then there are no stable deformations of the Schouten algebra as a \Lieinf-algebra\footnote{If true, the statement only holds in the stable setting, \cf footnote \ref{orientedregimefootnote} for a statement in the oriented setting. } \ie $(\Tpoly,\brdot_{\mathsf S})$ is rigid as a stable  \Lieinf-algebra. 
\end{enumerate}

We now turn to the case $d\neq2$. The following classification of the cohomology of the connected part of the full graph complex in low degrees is obtained from the various bounds collected in Section \ref{section:Cohomology of the full graph complex}:
\begin{Lemma}
\label{lemmacoho}
Let us denote $L_k$ the loop graph with $k$ edges. 
The cohomology of the (connected part of) the full graph complex $\fGCconn_d$ in low degrees for all $d\neq2$ is given by:

\begin{itemize}
\item  \underline{Degree $0$}: $H ^0(\fGCconn_{4j+3})=\corps\dl L_{4j+3}\dr$ for all $j\geq0$ and trivial otherwise.

\item \underline{Degree $1$}: $H^1(\fGCconn_1)=\corps\dl\Theta\dr$, $H ^1(\fGCconn_{4j+4})=\corps\dl L_{4j+5}\dr$ for all $j\geq-1$ and trivial otherwise.

\item \underline{Degree $2$}: $H^2(\fGCconn_{4j+1})=\corps\dl L_{4j+3}\dr$ for all $j\geq0$ and trivial otherwise.

\end{itemize}
\end{Lemma}
Note that the only non-loop cocycle in this classification is given by the $\Theta$-graph. The stable \Lieinf-structure induced by the Maurer--Cartan element \eqref{MoyalgraphMC} prolongating the latter is given by the Groenewold--Moyal bracket on symplectic manifolds --\cf Example \ref{GMc} -- which constitutes the unique\footnote{\label{orientedregimefootnote}Departing from the stable regime to the oriented regime, we note that the incarnation in $d=2$ of the $\Theta$-graph induces the Kontsevich--Shoikhet cocycle in $H^1(\GCor_2)$ whose prolongation to a Maurer-Cartan element is mapped to the Kontsevich--Shoikhet \Lieinf-algebra structure deforming the Schouten algebra of infinite dimensional polyvector fields, \cf \cite{Shoikhet2008a,Willwacher2015c}. Further, the incarnation of the $\Theta$-graph in $d=3$ yields a potential obstruction to the quantization of Lie bialgebroids \cite{Morand2021}. } stable structure in dimension $d=1$. It follows that the only stable structures in dimension $d>2$ are induced by loop classes\footnote{Stable structures induced from the Grothendieck--Teichm\"uller algebra $\grt_1$ only occur in dimension $d=2$. However, departing from the stable to the (multi)-oriented setting allows to generate stable structures from $\grt_1$ in dimension $d>2$, see \eg \cite{Morand2021} for an example on (quasi)-Lie bialgebroids in $d=3$.}:
\begin{Proposition}[Loop induced stable structures]
\label{propapp}
\hfill
\begin{enumerate}
\item Let $(\V,\om)$ be a $\NP$-manifold of (even) degree $n=4j+2$, $j\geq0$. The loop cocycle $L_{4j+3}$ induces a 1-dimensional family of stable \Lieinf-automorphisms of the $n$-Schouten algebra $\big(\Tpolyn,\brdot_{\mathsf S}\big)$.
\item Let $(\V,\om)$ be a $\NP$-manifold of (odd) degree $n=4j+3$, $j\geq0$. The loop cocycle $L_{4j+5}$ induces a 1-dimensional family of stable deformations of the $n$-Schouten algebra $\big(\Tpolyn,\brdot_{\mathsf S}\big)$ as a \Lieinf-algebra. 
\end{enumerate}
\end{Proposition}
\begin{proof}
The first statement is merely a rephrasing of the first item of Lemma \ref{lemmacoho}, given the definition of stable automorphism provided in Corollary \ref{thmcoho}. 
As for the second statement, since $H ^2(\fGCconn_{4j+4})$ vanishes for all $j\geq0$, there is no obstruction to the prolongation of the loop cocycle $L_{4j+5}$ into a formal Maurer--Cartan element\footnote{For any (graded) vector space $V$, we will denote $V\eps$ the (graded) vector space of formal series in the formal parameter $\epsilon$ with coefficients in $V$. Consistently, whenever $(V,\brdot)$ is a graded Lie algebra (or more generally a \Lieinf-algebra), we will denote $\MC(V\eps)$ the set of {\it formal} Maurer--Cartan elements \ie elements $\m\in\epsilon\, V^{1}\eps$ satisfying $\br{\m}{\m}=0$.} $\mathfrak m_{4j+5}\in\MC(\fGC_{4j+4}\eps)$. Corollary \ref{corMC} then ensures that the formal Maurer--Cartan element $\mathfrak m_{4j+5}$ is mapped via $\sRep^{(d)}$ to a stable deformation of the graded Lie algebra $(\Tpolyn,\brdot_{\mathsf S})$ as a \Lieinf-algebra\footnote{The only non-trivial higher brackets $l_m$ of the \Lieinf-algebra induced by the loop cocycle $L_{4j+5}$ have arities $m={p(4j+3)+2}$ 
for all $p\geq0$ with $l_2=\brdot_{\mathsf S}$. }. 
\end{proof}
\Prop{propapp} thus provides two mechanisms for generating stable structures on graded manifolds of specific degrees.
The following example illustrates the procedure in the odd degree case (for $j=0$): 

\begin{Example}
\label{exalinf}
Let $(\V,\om)$ be a $\NP$-manifold of degree 3 coordinatised by, \cf \eg \cite{Ikeda2010,Liu2016a,Grutzmann2010}:
\[
\pset{\underset{0}{x^\mu},\underset{1}{\psi^{\alpha}},\underset{2}{\chi_{\alpha}},\underset{3}{p_\mu}}\, .
\]
The pentagon graph:
\[
L_5:=\raisebox{-3ex}{\hbox{\begin{tikzpicture}[scale=0.5, >=stealth']
\tikzstyle{w}=[circle, draw, minimum size=4, inner sep=1]
\tikzstyle{b}=[circle, draw, fill, minimum size=2, inner sep=0.02]
\draw (18:1.3) node [b] (b3) {3};
\draw (90:1.3) node [b] (b2) {2};
\draw (162:1.3) node [b] (b1) {1};
\draw (234:1.3) node [b] (b5) {5};
\draw (306:1.3) node [b] (b4) {4};
\draw[black]  (b1) to (b2);
\draw[black]  (b1) to (b5);
\draw[black]  (b2) to (b3);
\draw[black]  (b3) to (b4);
\draw[black]  (b4) to (b5);
\end{tikzpicture}}}
\]
can be promoted to a formal Maurer--Cartan element $\mathfrak m_{5}\in\MC(\fGC_4\eps)$ reading:
\[
\mathfrak m_{5}:=\epsilon^5\, L_5+\epsilon^8\, \mathfrak m^{(8)}_{5}+\cdots+\epsilon^{3p+2}\, \mathfrak m^{(3p+2)}_{5}+\cdots
\]
 which induces a \Lieinf-algebra structure on the shifted graded algebra of functions $\Tpolyd{(3)}:=\fonc{\V}[3]$ with non-vanishing brackets $l_2, l_5, l_8, \ldots,l_{3p+2}$, $p\geq0$ such that, for all $X_i\in\Tpolyd{(3)}$:
\[ \begin{array}{rcl}
{l_2(X_1,X_2)} &=& {\sRep^{(4)}\big(\Gambr\big)(X_1,X_2)=\br{X_1}{X_2}_{\mathsf S}} \\
{l_5(X_1,\ldots,X_5)} &=& {\sRep^{(4)}\big(L_5
\big)(X_1,\ldots,X_5)} \\
&=& {s\un \mu_5\Big(\Delta_{12}\, \Delta_{23}\, \Delta_{34}\, \Delta_{45}\, \Delta_{51}\big(s(X_1),\ldots,s(X_5)\big)\Big)} \\
{} &{\vdots}& {} \\
{l_{3p+2}(X_1,\ldots,X_{3p+2})} &=& {\sRep^{(4)}\big( \mathfrak m^{(3p+2)}_{5}\big)(X_1,\ldots,X_{3p+2})}
\end{array} \]
where: 
\[
\Delta_{ij}=\frac{\p }{\p x^\mu_{(i)}}\, \frac{\p }{\p p_\mu^{(j)}}+\frac{\p }{\p p_\mu^{(i)}}\, \frac{\p }{\p x^\mu_{(j)}}+\frac{\p }{\p \psi^{\alpha}_{(i)}}\frac{\p }{\p \chi_{\alpha}^{(j)}}+\frac{\p }{\p \chi_{\alpha}^{(i)}}\frac{\p }{\p \psi^{\alpha}_{(j)}}\,.
\]
\end{Example}
\subsection{Stable Hamiltonian deformations}
\label{section:Stable deformations}
Up to now, our attention has been focused on $\NP$-manifolds. We will now consider additional structures on graded manifolds by focusing on $\NPQ$-manifolds \ie symplectic graded manifolds endowed with a Hamiltonian structure -- \cf Section \ref{section:NPQmanifolds} -- and discuss how the previously developed machinery can be used in order to generate Hamiltonian deformations from graph cocycles. More precisely, we discuss two separate mechanisms that allow to map graph cocycles to Hamiltonian deformations. The first one is a generalisation to any $d$ of the mechanism first identified by Kontsevich for $d=2$ \cite{Kontsevich1997} mapping elements from $H^0(\fGC_{d})$ to {\it stable} deformations (hence quantizations) of Hamiltonian functions [see \Prop{propunivhamdef} below]. The second mechanism is novel and allows to map graphs in $H^{-d}(\fGC_{d})$ to {\it Weyl} Hamiltonian deformations [see Corollary \ref{cor:WeylHamiltonian} below]. This new mechanism will prove to be particularly relevant for $d=3$,  mapping trivalent graphs to Weyl deformations of Courant algebroids. 

We will consider a  $\NP$-manifold $(\V,\om)$ of arbitrary degree $n\in\mathbb N$ with $d=n+1$. 
We will use the notation $[\mathcal{F}]$ to designate equivalence classes of functions $\mathcal F'\sim\mathcal F$ of $\V$ modulo gauge transformations [see \Defi{defi:gauge transfo}].
The set of Hamiltonian functions\footnote{Recall from Section \ref{section:NPQmanifolds} that a Hamiltonian function is a function $\cH\in\mathscr C^{\infty|d}{(\V)}$ such that $\pb{\cH}{\cH}_\om=0$. } on $(\V,\om)$ will be denoted $\Ham$. 
We start by defining the notion of \textbf{Hamiltonian deformation}:
\begin{Definition}[Hamiltonian deformation]
Let $\cH\in\Ham$ be a Hamiltonian function on $\V$. 

A Hamiltonian deformation of $\cH$ is a formal power series $\cH_*\in\epsilon\, \mathscr C^{\infty|d}{(\V)}\eps$ such that 
\begin{enumerate}
\item $\cH_*$ is nilpotent with respect to the graded Poisson bracket \ie $\pb{\cH_*}{\cH_*}_\om=0$.
\item The first order of the expansion of $\cH_*$ in terms of the formal parameter coincides with $\cH$ or equivalently $\cH=\frac{1}{\epsilon}\cH_*\big|_{\epsilon=0}$.
\end{enumerate}

The set of Hamiltonian deformations of a given Hamiltonian function $\cH$ will be denoted $\FHam_\cH$. A map $[\Ham]\to[\FHam]$ which assigns to each equivalence class [see Section \ref{section:NPQmanifolds}] $[\cH]\in\Ham$ an equivalence class $[\FHam_\cH]$ will be referred to as a \textbf{Hamiltonian deformation map}.
\end{Definition}
\begin{Example}[Formal Poisson structures]
Recall from Example \ref{exasymp} that Hamiltonian functions on $\NP$-manifolds of degree 1 are in one-to-one correspondence with Poisson structures on the corresponding base manifold. In this context, a Hamiltonian deformation is a map $[\pi]\mapsto[\pi_*]$ sending each Poisson bivector $\pi$ in the equivalence class $[\pi]$ (under diffeomorphism of the base manifold) to a formal Poisson structure $\pi_*=\epsilon\pi+\epsilon^2\pi_{(2)}+\cdots+\epsilon^k\pi_{(k)}+\cdots$ where the $\pi_{(i)}$'s are bivectors such that $\br{\pi_*}{\pi_*}_{\mathsf S}=0$ and the equivalence relation on $[\pi_*]$ is given by formal diffeomorphisms.
\end{Example}
The first defining condition of a Hamiltonian deformation ensures that the pullback  $\m:=s\un(\cH_*)$ is a formal Maurer--Cartan element of the graded Lie algebra $\big(\Tpolyn\eps,\brdot_{\mathsf S}\big)$ \ie $\m\in\MC(\Tpolyn\eps)$. This fact, combined with Corollary \ref{thmcoho}, yields the following Proposition:

\begin{Proposition}
\label{propunivhamdef}
There is a canonical map $H^0(\fGC_{d})\to([\Ham]\to[\FHam])$ mapping cocycles in the zeroth graph cohomology to Hamiltonian deformation maps. Such Hamiltonian deformation maps will be referred to as stable. 
\end{Proposition}
\begin{proof}
Recall from Corollary \ref{thmcoho} that to each cocycle $\Gamgraph\in H^0(\fGC_{d})$, one can associate a homotopy class of \Lieinf-automorphisms $\mathcal U^\Gamma=\exp\big(\sRep^{(d)}(\Gamgraph)\big)$ of $\big(\Tpolyn,\brdot_{\mathsf S}\big)$, where we denoted $\Gamma:=\exp(\gamma)$. The latter induces a bijective map between equivalence classes of formal Maurer--Cartan elements of $\big(\Tpolyn\eps,\brdot_{\mathsf S}\big)$ as $\hat{\mathcal U}^\Gamma:\MC(\Tpolyn\eps)\, /\sim\ \iso\ \MC(\Tpolyn\eps)\, /\sim\, :\, [\mathfrak m]\mapsto [\hat{\mathcal U}^\Gamma(\mathfrak m)]$ where $\hat{\mathcal U}^\Gamma(\mathfrak m):=\sum_{k=1}^\infty \frac{1}{k!}\, {\mathcal U}^\Gamma_k(\mathfrak m^{\otimes k})$.
Let $\cH\in\Ham$ be a Hamiltonian function. We will denote $\m_\cH\in\MC(\Tpolyn\eps)$ the canonical formal Maurer--Cartan element defined as $\m_\cH:=\epsilon\, s\un(\cH)$. The latter is mapped via $\hat{\mathcal U}^\Gamma$ to $\hat{\mathcal U}^\Gamma(\m_\cH)=\sum_{k=1}^\infty \frac{\epsilon^k}{k!}\, {\mathcal U}^\Gamma_k(s\un(\cH)^{\otimes k})$. It follows from the above reasoning that $\hat{\mathcal U}^\Gamma(\m_\cH)$ is a formal Maurer--Cartan element. Finally, we define $\cH_*:=s\big(\hat{\mathcal U}^\Gamma(\m_\cH)\big)$. Since $\hat{\mathcal U}^\Gamma_1$ is the identity of $\Tpolyn$, then $\cH=\frac{1}{\epsilon}\cH_*\big|_{\epsilon=0}$ and hence $\cH_*$ is a Hamiltonian deformation of $\cH$. We conclude that the map $[\mathsf{Ham}]\to[\mathsf{FHam}]:[\cH]\mapsto [\cH_*]$ is a Hamiltonian deformation map. Two equivalent cocycles $\gamma$ and $\gamma'$ (\ie such that there exists a degree $-1$ cochain $\chi$ such that $\gamma'-\gamma=\delta\chi$) induce homotopic \Lieinf-automorphisms thus yielding the same Hamiltonian deformation map. 
\end{proof}
The previous Proposition takes full advantage of the previously recalled Kontsevich machinery mapping graph cocycles to \Lieinf-automorphisms. As we will see, there exists an alternative way to map graph cocycles to Hamiltonian deformations, in a way that does not lift to a \Lieinf-automorphism. Before addressing them in the full case, we make a detour by considering infinitesimal deformations (or flows).

\paragraph{Infinitesimal deformations.}

Let $(\V,\om)$ be a $\NP$-manifold of arbitrary degree $n\in\mathbb N$ and denote $d=n+1$.

\begin{Proposition}
\label{propmapcoH}
There is a canonical map $H^\bullet(\fGC_{d})\to\big([\Ham]\to [H^{\bullet+d}(\fonc{\V}|\Q)]$\big). 
\end{Proposition}
Explicitly, cocycles of degree $p\in \mathbb Z$ in $\fGC_d$ with $N$ vertices are sent to applications mapping any Hamiltonian function $\cH\in\Ham$ to a degree $p+d$ cocycle:
\begin{eqnarray}
\label{eqRepC}
\mathcal C_{\gamma}(\cH):=\Rep^{(d)}_N(\Gamgraph)(\cH^{\otimes N})
\end{eqnarray}
in the cohomological complex induced on $\foncg{\bullet}$ by the associated\footnote{Note that the class $[H^{\bullet+d}(\fonc{\V}|\Q)]$ is independent of the choice of representative in the gauge class $[\cH]$ used to define $\Q$.
} homological vector field $\Q:=\pb{\cH}{\cdot}_\om$. 
\begin{proof}
Let $\gamma$ be a non-trivial cocycle of degree $p$ in $\fGC_d$ with $N$ vertices and $k$ edges. The cocycle condition $\deltabr\Gamgraph=0$ ensures that $\sRep^{(d)}(\Gamgraph)$ is a degree $p$ Chevalley--Eilenberg cocycle for the graded Lie algebra $\big(\Tpolyn,\brdot_{\mathsf S}\big)$. Explicitly, denoting $\mathcal U^\Gamma_N$ the component of $\sRep^{(d)}(\Gamgraph)$ acting on $N$ inputs\footnote{Or equivalently the first non-trivial Taylor coefficient, beside the identity, of the \Lieinf-automorphism $\mathcal U^\Gamma$, with $\Gamma:=\exp(\gamma)$, \cf the proof of \Prop{propunivhamdef}.}, the Chevalley--Eilenberg cocycle condition can be expressed as:
\begin{eqnarray}
\label{eqproofHam}
\big[{\brdot_{\mathsf S}},{\mathcal U^\Gamma_N}\big]_\NR&=&
\sum_{\sigma\in \text{\rm Sh}\un(N,1)}\hspace{-0.3cm}(-1)^{|\sigma|}(-1)^{p}\Sigma_{N+1}\big(\brdot_{\mathsf S}\circ_1 \mathcal U^\Gamma_N\big|\sigma)\\
\nn
&-&\sum_{\sigma\in \text{\rm Sh}\un(2,N-1)}
(-1)^{|\sigma|}\Sigma_{N+1}\big(\mathcal U^\Gamma_N\circ_1 \brdot_{\mathsf S}\big|\sigma)=0\, .
\end{eqnarray}
Acting on $s\un(\cH)^{\otimes{N+1}}$, the second term vanishes due to $\pb{\cH}{{\cH}}_\om=0$, thus yielding:
\begin{eqnarray}
\nn
\Big[{s\un(\cH)},{\mathcal U^\Gamma_N(s\un(\cH)^{\otimes{N}})}\Big]_{\mathsf S}=0\ .
\end{eqnarray}
 Denoting $\mathcal C_{\gamma}(\cH):=s\big( \mathcal U^\Gamma_N(s\un(\cH)^{\otimes{N}})\big)=\Rep^{(d)}_N(\Gamgraph)(\cH^{\otimes N})$ leads to the cocycle equation $\pb{\cH}{\mathcal C_{\gamma}(\cH)}_\om=0$. The latter function is of degree $|\mathcal C_{\gamma}(\cH)|=|\Gamgraph|+dN$ in $\fonc{\V}$ where $|\Gamgraph|$ stands for the degree of $\Gamgraph$ in $\Gra_d$ \ie $|\Gamgraph|=k(1-d)$. Since $\Gamgraph$ is assumed to be of degree $p$ in $\fGC_d$, then $d(N-1)+k(1-d)=p$ and thus $|\Gamgraph|=k(1-d)=p+d(1-N)$ so that $|\mathcal C_{\gamma}(\cH)|=p+d$. A reasoning similar to the above shows that two equivalent cocycles $\gamma'-\gamma=\delta\chi$ yield cohomologically equivalent functions $\mathcal C_{\gamma'}(\cH)-\mathcal C_{\gamma}(\cH)= \pb{\cH}{\mathcal C_{\chi}(\cH)}$ where we denoted $\mathcal C_{\chi}(\cH):=\Rep^{(d)}_{N-1}(\chi)(\cH^{\otimes N-1})$. Finally, two gauge equivalent Hamiltonians, \ie differing infinitesimally by a coboundary $\pb{\cH}{f}_\om$, yield two equivalent functions, differing by a coboundary $\pb{\cH}{\Rep^{(d)}_N(\Gamgraph)(f\otimes\cH^{\otimes N-1})}_\om$ and a gauge term $\pb{f}{\mathcal C_{\gamma}(\cH)}_\om$, as can be shown by letting \eqref{eqproofHam} act on $f\otimes\cH^{\otimes {N}}$.
\end{proof}
Of particular interest for us will be the two subcases $p=0$ and $p=-d$ allowing to map graph cocycles to flows on the space of Hamiltonian functions and Weyl rescalings of Hamiltonian functions, respectively.

\paragraph{Stable Hamiltonian flows $(p=0)$.}

\begin{Definition}[Hamiltonian flow]
A map $[\Ham]\to [H^d(\fonc{\V}|\Q)]$ mapping (gauge classes of) Hamiltonian functions $\cH\in\Ham$ to (gauge classes of) degree $d$ elements in the cohomology of the complex induced on $\foncg{\bullet}$ by the homological vector field $\Q:=\pb{\cH}{\cdot}_\om$ will be called a Hamiltonian flow.
\end{Definition}

The set of Hamiltonian flows will be denoted $\mathsf{Hflow}$. The following corollary of \Prop{propmapcoH} can be seen as a linearisation of \Prop{propunivhamdef}:
\begin{Corollary}
\label{propmapcoHd=0}
There is a canonical map $H^0(\fGC_{d})\to\mathsf{Hflow}$. Hamiltonian flows in the image of this map are called stable. 
\end{Corollary}
In plain words, the above corollary allows to map graph cocycles of degree $0$ to infinitesimal deformations of Hamiltonian functions -- or equivalently to flows on the space of Hamiltonian functions -- thus generalising to all $d$ the Kontsevich's construction of stable flows on the space of Poisson manifolds from cocycles in $H^0(\GC_2)$ [\cf Section 5.3 in \cite{Kontsevich1997} and Section \ref{section:Poisson manifolds} below].

As noted in Section \ref{section:Cohomology of the full graph complex}, the zeroth cohomology is the dominant degree of the cohomology of $\GC_2$, being isomorphic to the infinite dimensional Grothendieck--Teichm\"uller Lie algebra $\mathfrak{grt}_1$. The previous construction for $d=2$ thus allows to generate infinitely many Hamiltonian flows on the space of Poisson bivectors. Explicit examples of stable Hamiltonian flows for $d=2,3$ will be exhibited in Section \ref{section:Applications}.

\paragraph{Stable Weyl flows $(p=-d)$.}
\begin{Definition}[Weyl flow]
A map $[\Ham]\to [H^0(\fonc{\V}|\Q)]$ mapping (gauge equivalence classes of) Hamiltonian functions $\cH\in\Ham$ to (gauge equivalence classes of) degree $0$ elements in the cohomology of the complex induced on $\foncg{\bullet}$ by the homological vector field $\Q:=\pb{\cH}{\cdot}_\om$ will be called a Weyl factor.
\end{Definition}
Explicitly, Weyl flows map Hamiltonian functions $\cH\in\Ham$ to functions\footnote{Note that there are no functions of degree $-1$ in $\fonc{\V}$ hence the definition of Weyl factors does not involve equivalence classes.} $\Om\in\fonc{\M}$ -- called \textbf{Weyl factors} -- satisfying $\Q[{\Om}]=0$ [or equivalently $\pb{\cH}{\Om}_\om=0$ with $\cH\in \Ham$ the Hamiltonian function associated to $\Q$]. The set $H^0(\fonc{\V}|\Q)$ of Weyl factors is a commutative algebra with product the pointwise product of functions on the base manifold. 

We now show that cocycle graphs of degree $-d$ induce Weyl factors as follows:

\begin{Proposition}
\label{propconfflow}
Let $(\mathcal V,\pbdot_\om,\Q)$ be a $\NPQ$-manifold. There is a canonical morphism of commutative algebras $H^{-d}(\fGC_{d})\to H^0(\fonc{\V}|\Q)$ mapping cocycles in the degree $-d$ graph cohomology to Weyl factors. Such Weyl factors will be referred to as stable. 
\end{Proposition}
Indeed, it follows from \Prop{propmapcoH} for $p=-d$ that there exists a canonical map $H^{-d}(\fGC_{d})\to\big([\Ham]\to [H^{0}(\fonc{\V}|\Q)]\big)$, reading, for any particular Hamiltonian function $\cH\in\Ham$ as:
\begin{eqnarray}
H^{-d}(\fGC_{d})&\to& H^0(\fonc{\V}|\Q)\nn\\
\gamma&\mapsto& \Om_{\gamma}(\cH):=\Rep^{(d)}_N(\Gamgraph)(\cH^{\otimes N})\label{eq:defWeyl}
\end{eqnarray} 
where $\Omega_{\gamma}(\cH):=\mathcal C_{\gamma}(\cH)$ is defined as in eq.\eqref{eqRepC}, with $N$ the number of vertices of the cocycle $\gamma$.
\begin{proof}
It remains to be shown that the above map $\gamma\mapsto \Om_\gamma(\cH)$ is a morphism of commutative algebras which is ensured by the fact that $\Rep^{(d)}$ preserves the concatenation product \eqref{property}. 
\end{proof}

\subsection{Weyl Hamiltonian deformations}

  \begin{Corollary}
  \label{cor:WeylHamiltonian}
  There are infinitely many maps $H^{-d}(\fGC_{d})\to([\Ham]\to[\FHam])$ mapping cocycles in the graph cohomology of degree $-d$ to Hamiltonian deformation maps. Such Hamiltonian deformation maps will be referred to as Weyl. 
\end{Corollary}
Let $\gamma$ be a non-trivial cocycle of degree $-d$ in $\fGC_d$ with $N>0$ vertices. Given a Hamiltonian function $\cH\in\Ham$, we can define the Weyl factor $\Om_\gamma(\cH)$ as in \eqref{eq:defWeyl}. 
The latter can be used to define a Hamiltonian deformation map $[\Ham]\to[\FHam]:[\cH]\mapsto[\cH_*]$ where $\cH_*$ reads as:
\begin{eqnarray}
\nn
\cH_*=\epsilon\cH+\sum_{k=1}^{\infty}a_k\, \epsilon^{kN+1}\, \Om_\gamma(\cH)^k\cdot \cH
\end{eqnarray}
where the coefficients $a_k\in\corps$ are arbitrary.

 \begin{proof}
 We need to check that $\cH_*$ is a Hamiltonian deformation of $\cH$. The condition $\pb{\cH_*}{\cH_*}_\om=0$ follows straightforwardly from $\pb{\cH}{\cH}_\om=0$ ($\cH$ being Hamiltonian), $\pb{\cH}{\Om_\gamma(\cH)}_\om=0$ (\Prop{propmapcoH}) and $\pb{\Om_\gamma(\cH)}{\Om_\gamma(\cH)}_\om=0$ (holds identically from degree consideration). The second condition $\cH=\frac{1}{\epsilon}\cH_*\big|_{\epsilon=0}$ follows directly from $N>0$.
  \end{proof}
Contradistinctly to stable Hamiltonian deformation maps [\Prop{propunivhamdef}], Weyl Hamiltonian deformation maps do not descend from \Lieinf-automorphisms of $\Tpolyn$ and thus constitute a novel way to define consistent deformations from graph cocycles. 

Similarly as before, one can make use of the results reviewed in Section \ref{section:Cohomology of the full graph complex} in order to classify stable Weyl factors. In fact one can check that stable Weyl factors originating from connected graphs only occur in $d=3$:
 
\begin{Proposition}
 $H^{-d}(\fGCconn_d)=\mathbf{0}$ for all $d\neq3$.
\end{Proposition}
\begin{proof}
Let $\gamma$ be a graph of $\fGCconn_d$ with $N$ vertices and $k$ edges. The corresponding degree is given by $|\Gamgraph|_d=\dime(N-1)+k(1-\dime)$. Imposing $|\gamma|_d=-d$ thus yields the constraint: $dN=k(d-1)$. Setting $d=0$ in the latter constraint imposes $k=0$, hence the corresponding graph is either disconnected or the trivial graph or the graph with one unique vertex, which is never a cocycle. We will then assume $d>0$. The above constraint thus yields $N=k-\frac{k}{d}$ which is only integer for either $k=0$ (ruled out by the above reasoning) or $k\geq d$. Hence we will be focusing on the cases $0<d\leq k$.
We first observe that graphs satisfying the above constraint cannot be (non-trivial) loops as setting $k=N$ would yield $N=0$ for all $d$. Hence Theorem \ref{conncohomo} ensures that $H^{-d}(\fGCconn_d)=H^{-d}(\GC_d)$. We now distinguish between even and odd $d$. 
\begin{itemize}
\item \emph{Even $d$}: As noted earlier, the cohomology $H^{\bullet}(\GC_d)$ for even $d$ is isomorphic to $H^{\bullet}(\GC_2)$, although the isomorphism is not degree preserving. According to the bounds on cohomology degree displayed in Theorem \ref{thmWill}, for $\gamma$ to be a non-trivial cocycle, one needs $1<|\gamma|_2<b-2$, where $b$ stands for the first Betti number associated to $\gamma$. In the case $d=2$, we have $|\gamma|_2=-2$ hence the first inequality ensures that $H^{-2}(\GC_2)$ is empty. For $d\geq4$, the second part of the inequality reads $|\gamma|_2<b-2\Leftrightarrow 3N-2k-1<0$. Multiplying by $d>0$ on both sides and using the constraint yields $k<\frac{d}{d-3}$. Combining with $0<d\leq k$ yields $d<\frac{d}{d-3}$ which admits no solution for $d>0$. Hence $H^{-d}(\GC_d)$ is empty for even $d$.

\item \emph{Odd $d$}: In the odd case, the cohomology $H^{\bullet}(\GC_d)$ is isomorphic to $H^{\bullet}(\GC_3)$. According to the bounds on cohomology degree displayed in Section \ref{section:Cohomology of the full graph complex}, the inequality $-b-2<|\gamma|_3<-2$ needs to hold for $\gamma$ to be a non-trivial cocycle. The second inequality reads more explicitly $3N-2k-1<0$. Multiplying by $d>0$ on both sides and using the constraint yields $k(d-3)<d$. Assuming $d\neq3$, and using the previous reasoning ensures that $H^{-d}(\GC_d)$ is empty for odd $d\neq3$. \qedhere
\end{itemize}
\end{proof}
It follows that stable Weyl factors induced by connected graphs only occur
 in dimension $d=3$. The corresponding cohomology space is then $H^{-3}(\GC_3)$ which happens to coincide with the dominant degree of $\GC_3$ [\cf Section \ref{section:Cohomology of the full graph complex}], thus allowing to generate infinitely many stable Weyl factors for Courant algebroids from trivalent graphs modulo IHX relations \big(see Figure \ref{figtrivalent} and eq.\eqref{IHX}\big). Explicit examples thereof will be displayed in Section \ref{section:Courant algebroids}.

\section{Stable deformations of symplectic Lie $n$-algebroids}
\label{section:Applications}
The aim of the present section is to illustrate some of the machinery developed in Section \ref{section:Stable structures on graded manifolds} to the case of $\NPQ$-manifolds of degrees 1 and 2. As recalled in Example \ref{exasymp}, the associated geometric notions (\ie symplectic Lie $1,2$-algebroids) identify with the one of Poisson manifolds and Courant algebroids, respectively. 

\subsection{Poisson manifolds ($n=1$)}
\label{section:Poisson manifolds}

As shown in \cite{Roytenberged.Contemp.Math.Vol.315Amer.Math.Soc.ProvidenceRI2002}, $\NP$-manifolds $\V$ of degree $1$ are in bijective correspondence with ordinary smooth manifolds $\M$ via the identification of $\V$ with the shifted cotangent bundle $T^*[1]\M$. The tower of fibrations \eqref{fibration} thus reduces to the vector bundle structure $T^*[1]\M\to \M$. The graded Poisson algebra of functions on $T^*[1]\M$ is isomorphic to the Gerstenhaber algebra of polyvector fields $\Tpoly$ and Hamiltonian functions are in bijection with Poisson bivectors on $\M$.

The representation morphism $\Rep^{(2)}:\Gra_2\hookrightarrow\End_{\fonc{T^*[1]\M}}$ of the 2-dimensional graph operad $\Gra_2$ on the space of functions of the shifted cotangent bundle was first introduced by M. Kontsevich in \cite[Section 5.2]{Kontsevich1997}  and reads as \eqref{Orientationmorphism2} with $\Delta$ given by [\cf eq.\eqref{defDeltadeven2}]:
\begin{eqnarray}
\nn
\Delta_{ij}=\frac{\p }{\p x^\mu_{(i)}}\, \frac{\p }{\p p_\mu^{(j)}}+\frac{\p }{\p p_\mu^{(i)}}\, \frac{\p }{\p x^\mu_{(j)}}\, .
\end{eqnarray}
Following the leitmotiv of Section \ref{section:Stable structures on graded manifolds}, the representation morphism $\Rep^{(2)}$ can be used in order to induce stable structures on $\M$. In particular, using the isomorphism $H^{0}(\GC_2)\simeq \mathfrak{grt}_1$, it follows from Corollary \ref{thmcoho} that the Grothendieck--Teichm\"uller group $\GRT_1:=\exp(\grt_1)$ acts  via \Lieinf-automorphisms on the Schouten algebra $\Tpoly$, see  \cite{Willwacher2015,Jost2012}. 

At the linear level, Corollary \ref{propmapcoHd=0} ensures that cocycles in  $H^{0}(\GC_2)$ yield stable flows on the space of Poisson bivectors. In other words, given a manifold $\M$ and a cocycle $\Gamgraph\in H^{0}(\GC_2)$ with $N$ vertices, one can define a map $\pi\mapsto\dot\pi$ mapping Poisson bivectors $\pi\in\field{\w^2T\M}$ (thus satisfying $\br{\pi}{\pi}_{\mathsf S}=0$) on $\M$ to stable {\it Lichnerowicz cocycles} \ie bivectors $\dot\pi\in\field{\w^2T\M}$ satisfying $\delta_\pi\dot\pi:=\br{\pi}{\dot\pi}_{\mathsf S}=0$. Concretely, this is done by first defining the function $\cH=\frac{1}{2}\pi^{\mu\nu}(x)p_\mu\,  p_\nu$ -- which can be checked to be Hamiltonian (\ie $\pb{\cH}{\cH}_\om=0$) as a consequence of the fact that $\pi$ is Poisson -- and then define the function $\dot{\cH}:=\Rep^{(2)}_N(\Gamgraph)(\cH^{\otimes N})$ -- satisfying $\pb{\cH}{\dot\cH}_\om=0$ as a consequence of $\delta\Gamgraph=0$ and $\pb{\cH}{\cH}_\om=0$. Finally, one defines $\dot\pi$ as the principal symbol of the function $\dot\cH$ \ie $\dot\cH=\frac{1}{2}\dot\pi^{\mu\nu}p_\mu\,  p_\nu$.

The simplest example of the previous construction is given by the {\it tetrahedral flow} introduced in \cite[Section 5.3]{Kontsevich1997} and further studied in \cite{Bouisaghouane2016,Bouisaghouane2016a}. The latter is induced by the tetrahedron graph $\gamma_3\in H^{0}(\GC_2)$ [\cf Proposition \ref{propWillwheel} and Figure \ref{figwheel}]. An explicit expression for the map $\pi\mapsto\dot\pi$ can be obtained by first using the orientation morphism \eqref{morphdg Lie algebra} on $\gamma_3$ as in Figure \ref{figormorph}, yielding a linear combination of four directed graphs. Decorating vertices with copies of the Hamiltonian function $\cH$ and interpreting edges as differential operators $\bar\Delta$ [see eq.\eqref{defDeltadeven}], the first two graphs vanish since they include vertices with more than two outgoing edges. The two remaining graphs yield the following local expression\footnote{Although the two terms of eq.\eqref{tetracocycle} already appeared in \cite{Kontsevich1997}, the relative factor $1:6$ was only recently obtained in \cite{Bouisaghouane2016,Bouisaghouane2016a} where it was also shown to constitute the unique choice allowing for the cocycle property to hold. } for the Lichnerowicz cocycle $\dot \pi$ associated with the Poisson bivector:
\begin{eqnarray}
\dot\pi^{\mu \nu}=\partial_{\epsilon}{\pi^{\alpha \beta}} \partial_{\alpha}{\pi^{\gamma \delta}} \partial_{\gamma}{\pi^{\epsilon \lambda}} \partial_{\beta \delta \lambda}{\pi^{\mu \nu}}+6\, \partial_{\epsilon}{\pi^{\alpha \beta}} \partial_{\alpha}{\pi^{\gamma \delta}} \partial_{\gamma \lambda}{\pi^{\epsilon [\mu}} \partial_{\beta \delta}{\pi^{\nu] \lambda}}.\label{tetracocycle}
\end{eqnarray}
Furthermore, it follows from \Prop{propunivhamdef} that the cocycle $\dot\pi$ can be promoted to a full formal Maurer-Cartan element in $(\Tpoly\eps,\delta_\pi,\brdot_{\mathsf S})$ thus yielding a stable formal Poisson structure\footnote{We refer to \cite{Alm2012} for results regarding the Lichnerowicz cohomology associated with stable deformations $\pi_*$ of Poisson manifolds.  
} $\pi_*=\pi+\epsilon^4\dot\pi+\cdots$  such that $\br{\pi_*}{\pi_*}_{\mathsf S}=0$. 

Note that the tetrahedral flow is only the first and simplest example of an infinite set of stable flows on the space of Poisson bivectors provided by elements in the Grothendieck--Teichm\"uller algebra $\grt_1$. We refer in particular to \cite{Buring2017c} and  \cite{Buring2017a} for results regarding the flows associated with the pentagon $\gamma_5$ and heptagon graphs $\gamma_7$, respectively. 

\paragraph{Relation to quantization.}
Before concluding with the $n=1$ case, we recall known results regarding the deformation quantization problem for Poisson manifolds. Our emphasis will be on the classification problem for formality morphisms and how the above results regarding stable deformations of Poisson structures can be used to shed light on the matter. Such considerations will hopefully provide guiding lines in order to address cases for which the deformation quantization problem is less well understood [\cf Section \ref{section:Courant algebroids} for a related discussion on Courant algebroids]. First, recall from the Introduction that Kontsevich's solution to the deformation quantization problem for Poisson manifolds involves a formality morphism \eqref{QI} \ie a \Lieinf quasi-isomorphism between the Schouten algebra on $\Tpoly$ and the Hochschild dg Lie algebra of multidifferential operators $\Dpoly$. As emphasised earlier, Kontsevich's formality morphism is {\it stable} in a precise sense introduced in \cite{Dolgushev2011}. The set of (homotopy classes\footnote{We refer to \cite{Dolgushev2007b} for a definition of  the notion of homotopy equivalence between  \Lieinf-morphisms. } of) stable formality morphisms of the form \eqref{QI} will be denoted $\SQI$.   
A first incarnation of the Grothendieck--Teichm\"uller group as playing a classification r\^ole for $\SQI$ stems from a construction due to D. Tamarkin in his formulation of an alternative proof to Kontsevich's formality theorem \cite{Tamarkin2007,Hinich2000}. The latter provides a bijective map $\mathcal U:\DAss\iso\SQI$ where $\DAss$ stands for the set of Drinfel'd associators. As mentioned earlier, the set $\DAss$ is a $\GRT_1$-torsor thus providing an (implicit) action of the Grothendieck--Teichm\"uller group on $\SQI$. However, Tamarkin's map is far from being explicit making it difficult to precisely characterise the corresponding $\GRT_1$-action on quantization procedures. The situation has been clarified by V. A. Dolgushev who showed  in \cite{Dolgushev2011} that the set $\SQI$ is naturally endowed with a regular action of the pro-unipotent group $\exp\big(H^0(\GC_{2})\big)$. This result, combined with T. Willwacher's isomorphism $H^0(\GC_{2})\simeq\grt_1$ \cite{Willwacher2015} defines a regular $\GRT_1$-action on $\SQI$, so that both sides of Tamarkin's map $\mathcal U:\DAss\iso\SQI$ are $\GRT_1$-torsors. It has furthermore been shown in \cite{Dolgushev2014a} that Tamarkin's map is equivariant with respect to the action of $\GRT_1$ \ie $\mathcal U$ is a bijection of $\GRT_1$-torsors. Under this bijection, the (homotopy class of) Kontsevich's morphism with standard (or harmonic) propagator \cite{Kontsevich:1997vb} is mapped to the Alekseev--Torossian associator \cite{Alekseev2010,Severa2009} (see \cite{Dolgushev2014a}) while the (homotopy class of) Kontsevich's morphism with logarithmic propagator \cite{KontsevichLett.Math.Phys.48:35-721999,Alekseev2014} is mapped \cite{Rossi2014} to the Knizhnik--Zamolodchikov associator \cite{Drinfeld1991}. 

In this picture, one can argue that the map assigning to each group element in $\exp\big(H^0(\GC_{2})\big)\simeq \GRT_1$ a homotopy class of \Lieinf-automorphisms of $\Tpoly$ [\cf \cite{Willwacher2015,Jost2012} and Corollary \ref{thmcoho} for its generalisation to all $d$] constitutes a useful intermediate step allowing a precise characterisation of the arbitrariness in quantization procedures. In order to illustrate this, we let $\Phi\in\DAss$ be a Drinfel'd associator and denote $[\mathcal U_\Phi]:\Tpoly\iso\Dpoly$ the homotopy class of formality morphisms associated with $\Phi$ through Tamarkin's procedure. Let furthermore $\Gamgraphgroup\in \exp\big(H^0(\GC_{2})\big)\simeq \GRT_1$. The following diagram commutes (in the category of \Lieinf-algebras with homotopy classes of \Lieinf quasi-isomorphisms as morphisms):
\begin{eqnarray}
\begin{split}
\xymatrix{
\Tpoly\ar[d]_{[{\mathcal U}^{\Gamgraphgroup}]}\ar[rd]^{[{\mathcal U}_{\Phi\cdot\Gamgraphgroup}]}\ar[r]^{[{\mathcal U}_{\Phi}]}&\Dpoly\ar[d]^{[{\mathcal U}_{\mathcal D}(\Phi,\Gamgraphgroup)]}\\
\Tpoly\ar[r]^{[{\mathcal U}_{\Phi}]}&\Dpoly
}\label{diag1}
\end{split}
\end{eqnarray}
where:
\begin{itemize}
\item $[\mathcal U_{\Phi\cdot\Gamgraphgroup}]:\Tpoly\iso\Dpoly$ denotes the homotopy class of formality morphisms associated with the Drinfel'd associator $\Phi\cdot\Gamgraphgroup$.
\item $[\mathcal U^\Gamgraphgroup]:\Tpoly\iso\Tpoly$ denotes the homotopy class of \Lieinf-automorphisms of $\Tpoly$ associated with the element $\Gamgraphgroup\in\GRT_1$ through Corollary \ref{thmcoho} (for $d=2$).
\item $[{\mathcal U}_{\mathcal D}(\Phi,\Gamgraphgroup)]:\Dpoly\iso\Dpoly$ denotes the homotopy class of \Lieinf quasi-isomorphisms from $\Dpoly$ to itself associated with the pair $(\Phi,\Gamgraphgroup)$ and defined through $[{\mathcal U}_{\mathcal D}(\Phi,\Gamgraphgroup)]=[\mathcal U_{\Phi}\circ{\mathcal U}^{\Gamgraphgroup}\circ\mathcal U_{\Phi}\un]$, where $\mathcal U_{\Phi}\un$ is a homotopy inverse\footnote{Here, we use the fact that a \Lieinf-morphism is a quasi-isomorphism if and only if it is a homotopy equivalence [Section 3.7 of \cite{Keller2007} for a statement as well as references therein]. This implies that there exists a (non-canonical)  \Lieinf quasi-isomorphism $\mathcal U_{\Phi}\un:\Dpoly\to\Tpoly$ such that $\mathcal U_{\Phi}\circ \mathcal U_{\Phi}\un\sim \Id_{\Dpoly}$ and $\mathcal U_{\Phi}\un\circ \mathcal U_{\Phi}\sim\Id_{\Tpoly}$.  Contrarily to its counterpart $\mathcal U^\Gamma$, the family of \Lieinf-automorphisms ${\mathcal U}_{\mathcal D}$ does depend on the existence of a Drinfel'd associator (although the explicit choice does not matter due to the equivariance relation ${\mathcal U}_{\mathcal D}(\Phi\cdot \Gamgraphgroup',\Gamgraphgroup)={\mathcal U}_{\mathcal D}(\Phi,\text{Ad}_{\Gamgraphgroup'}\Gamgraphgroup)$ with $\text{Ad}_{\Gamgraphgroup'}\Gamgraphgroup=\Gamgraphgroup'\cdot\Gamgraphgroup\cdot\Gamgraphgroup'{}\un$).} of the representative $\mathcal U_{\Phi}$.
\end{itemize}
The regularity of the action of $\GRT_1$ on $\SQI$ can be restated as follows: for any pair of homotopy classes of stable formality morphisms $[\mathcal U_{\Phi}]$ and $[\mathcal U_{\Phi'}]$, there exists a unique element $\Gamma\in\GRT_1$ such that $[\mathcal U_{\Phi'}]=[\mathcal U_{\Phi}\circ \mathcal U^\Gamma]$, for any representatives $\mathcal U_{\Phi}$ and $\mathcal U^\Gamma$. The space $\SQI$ of homotopy classes of stable formality morphisms can then be fully explored by composition with \Lieinf-automorphisms of $\Tpoly$ induced from $\GRT_1$. Such a reasoning can also be shown to hold at the level of quantization maps. Indeed, each arrow appearing in Diagram \ref{diag1} is a (homotopy class of) \Lieinf quasi-isomorphisms and thus induces a bijection between (equivalence classes of) Maurer--Cartan sets [\cf footnote \ref{footDGnew}] as: 
\begin{eqnarray}
\nn
\xymatrix{
\quad\mathsf{FPoiss}\quad\ar[d]_{[\hat{\mathcal U}^{\Gamgraphgroup}]}\ar[r]^{[\hat{\mathcal U}_{\Phi}]}\ar[rd]^{[\hat{\mathcal U}_{\Phi\cdot\Gamgraphgroup}]}&\quad\mathsf{Star}\quad\ar[d]^{[\hat{\mathcal U}_{\mathcal D}(\Phi,\Gamgraphgroup)]}\\
\quad\mathsf{FPoiss}\quad\ar[r]^{[\hat{\mathcal U}_{\Phi}]}&\quad\mathsf{Star}\quad
}
\quad\quad
\xymatrix{
\quad[\pi]\quad\ar@{|->}[d]_{[\hat{\mathcal U}^{\Gamgraphgroup}]}\ar@{|->}[r]^{[\hat{\mathcal U}_{\Phi}]}\ar@{|->}[rd]^{[\hat{\mathcal U}_{\Phi\cdot\Gamgraphgroup}]}&\quad[*]\quad\ar@{|->}[d]^{[\hat{\mathcal U}_{\mathcal D}(\Phi,\Gamgraphgroup)]}\\
\quad[\pi']\quad\ar@{|->}[r]^{[\hat{\mathcal U}_{\Phi}]}&\quad[*']\quad
}
\end{eqnarray}
Mimicking the above reasoning allows to span the whole space of stable quantization maps ${\hat{\mathcal U}_{\Phi}}$ by composition with stable deformation maps ${\hat{\mathcal U}^{\Gamgraphgroup}}$ induced from $\GRT_1$ [\cf \Prop{propunivhamdef}]. The resulting characterisation of the action of $\GRT_1$ on stable quantization maps in terms of stable deformations has the merit to make certain features relatively explicit. As an example, it follows from the previous reasoning that formal Poisson structures $[\pi]$ which are insensitive to stable deformations\footnote{That is, such that ${\hat{\mathcal U}^{\Gamgraphgroup}}([\pi])=[\pi]$ for all $\Gamgraphgroup\in\GRT_1$.} admit canonical quantizations \ie their quantum class is unique\footnote{In other words, the associated (class of) star products $[*]=\hat{\mathcal U}_\Phi([\pi])$ does not depend on the choice of Drinfel'd associator $\Phi$. }. Straightforward reasoning on the number of derivatives involved in stable deformations [see \eg eq.\eqref{tetracocycle}] entails that Poisson bivectors whose local description is at most quadratic in coordinates admit a unique local quantization. This is in particular the case for constant Poisson bivectors (and in particular for symplectic manifolds in Darboux coordinates) which are uniquely (locally) quantized by the {Groenewold--Moyal} star product \cite{Groenewold1946,Moyal1949}, \cf eq.\eqref{GMproduct}. Slightly less trivial is the Kostant--Souriau--Kirillov Poisson bracket -- defined on the dual of any Lie algebra -- which is linear in coordinates. The latter admits two known quantizations, namely the Gutt \cite{Gutt1983,Esposito2015a} and Kontsevich \cite{Kontsevich:1997vb} star products. According to the previous reasoning, these two star products must belong to the same equivalence class. However, they do not coincide, as shown for example in \cite{Kontsevich:1997vb,Kathotia2009,Shoikhet2007,Dito2007,Arnal1999}. Rather, they are related via an isomorphism given by the Duflo map (\cf \eg Theorem 14 in \cite{Felder2008}).
\subsection{Courant algebroids ($n=2$)}
\label{section:Courant algebroids}

The present section applies the results of Section \ref{section:Stable structures on graded manifolds} to symplectic Lie $2$-algebroids. The latter notion identifies with the one of Courant algebroids that we now review, following the presentation \ala Dorfman (\cf \eg \cite{Jurco2016} for details and \cite{Kosmann-Schwarzbach2013} for a historical account).
\begin{Definition}[Courant algebroid]
\label{defCourantD}
A Courant algebroid is a quadruplet $(E,\dldot,\mD,\brdot_E)$ where:
\begin{itemize}
\item The pair $(E,\dldot)$ is a \textbf{pseudo-Euclidean vector bundle} \ie
\begin{itemize}
\item $E\to\M$ is a vector bundle over the smooth manifold $\M$. We will denote $(\foncm,\cdot)$ the commutative associative algebra of functions on $\M$ and $*:\foncm\otimes\fE\to\fE$ the module structure on fibers of $E$. The latter satisfies the associativity relation $f*(g*X)=(f\cdot g)*X$ for all $f,g\in\foncm$ and $X\in\fE$.
\item The map $\dldot:\fE\otimes\fE\to\foncm$ satisfies the following conditions:
\begin{enumerate}
\item $\fonc{\M}$-bilinear \ie $\dlr{f* X}{Y}_E=\dlr{X}{f* Y}_E=f\cdot \dlr{X}{Y}_E$
\item symmetric \ie $\dlr{X}{Y}_E=\dlr{Y}{X}_E$
\item non-degenerate \ie $\dlr{X}{Y}_E=0$ for all $Y\in\fE$ $\Leftrightarrow$ $X=0$.
\end{enumerate}
A bilinear form satisfying these conditions will be referred to as a \textbf{fiber-wise metric}. 
\end{itemize}
\item The pair $(\mD,\brdot_E)$ is a \textbf{Courant--Dorfman structure} on $(E,\dldot)$ \ie
\begin{itemize}
\item $\brdot_E:\fE\otimes\fE\to\fE$ is a $\corps$-bilinear form on the fibers of $E$ called the \textbf{Dorfman bracket}. 
\item $\mD:\foncm\to\fE$ is a $\corps$-linear derivation \ie $\mD(f\cdot g)=f*\mD f+g*\mD f$ for all $f,g\in\foncm$. The derivation $\mD$ defines a $\foncm$-linear map $\rho:\fE\to\vf$ called the \textbf{anchor} as $\rho_X[f]=\dlr{X}{\mD f}_E$ for all $f\in\foncm$, $X\in\fE$.
\end{itemize}
such that the following conditions are satisfied:
\begin{enumerate}
\item The Dorfman bracket satisfies the Jacobi identity in its Leibniz form:
\begin{eqnarray}
\nn
\br{X}{\br{Y}{Z}_E}_E=\br{\br{X}{Y}_E}{Z}_E+\br{Y}{\br{X}{Z}}_E\text{ for all }X,Y,Z\in\field{E}
\end{eqnarray}
so that the pair $(\fE,\brdot_E)$ is a $\corps$-Leibniz algebra.

\item The symmetric part of the Dorfman bracket is controlled by the derivation $\mD$ as:
\begin{eqnarray}
\nn
\br{X}{Y}_E+\br{Y}{X}_E=\mathcal D\dlr{X}{Y}_E\text{ for all }X,Y\in\fE\, .
\end{eqnarray}

\item The fiber-wise metric $\dldot$ is compatible with the Courant--Dorfman structure $(\mD,$ $\brdot_E)$, \ie
\begin{eqnarray}
\nn
\dlr{X}{\mD\dlr{Y}{Z}_E}_E=\dlr{\br{X}{Y}_E}{Z}_E+\dlr{Y}{\br{X}{Z}_E}_E=0\text{ for all }X,Y,Z\in\fE\, .
\end{eqnarray}
\end{enumerate}
\end{itemize}
\end{Definition}
Introducing a basis $\pset{e_a}_{a=1, \ldots, \dim E}$ of the space of sections $\fE$ allows to provide a component expression of the Courant algebroid maps as follows:
\begin{itemize}
\item In components, the fiber wise metric reads $\dlr{X}{Y}_E=\kappa_{ab}\, X^a\, Y^b$ where the constant matrix $\kappa$ satisfies:
\begin{enumerate}
\item $\kappa$ is symmetric \ie $\kappa_{ab}=\kappa_{ba}$.
\item $\kappa$ admits an inverse $\kappa\un$ such that $\kappa^{ac}\kappa_{cb}=\delta^a_b$ with $\delta$ the Kronecker delta.
\end{enumerate}
\item The component expression for the Courant--Dorfman structure $(\mD,\brdot_E)$ is captured by a pair $(\rho_a{}^\mu,T_{abc})$, where $T_{abc}$ is totally skewsymmetric. Explicitly, we have:
\begin{itemize}
\item $\mathcal D$-map: $\mathcal Df=\kappa^{ab}\rho_b{}^\mu\, \p_\mu f\, e_a$
\item Anchor: $\rho_X[f]=X^a\rho_a{}^\mu\p_\mu f$
\item Dorfman bracket: $\br{X}{Y}_E=\big(\rho_X[Y^a]-\rho_Y[X^a]-T_{bc}{}^aX^bY^c+\kappa^{ab}\rho_b{}^\mu\p_\mu X^c\kappa_{cd}Y^d\big) e_a$
\end{itemize}
where indices are raised and lowered with $\kappa$.

It can be checked that the defining conditions of a Courant algebroid are satisfied if and only if the pair $(\rho_a{}^\mu,T_{abc})$ satisfies the set of conditions:
\begin{enumerate}
\item ${\mathcal C_1}^{\mu \nu}:=\rho_{a}{}^{\mu}\kappa^{ab}\rho_{b}{}^{\nu}=0$
\item ${\mathcal C_2}^{\mu}_{a b}:=\rho_{ c}{}^{\mu}\kappa^{cd}T_{d a b}+2\, \rho_{[a}{}^{\lambda}\, \p_{\lambda}\rho_{b]}{}^{\mu}=0$
\item ${\mathcal C_3}_{a b c d}:=\frac{1}{4}T_{e [a b}\kappa^{ef}T_{c d] f}+\frac{1}{3}\rho_{[a}{}^{\mu}\, \p_\mu T_{bcd]}=0$.
\end{enumerate}
\end{itemize}
Comparing this set of constraints with \eqref{ConstraintCA1}-\eqref{ConstraintCA3} allows to relate Courant algebroids with symplectic Lie $2$-algebroids (or $\NPQ$-manifolds of degree 2). The precise nature of this relation is articulated in the following theorem:

\begin{Theorem}[D. Roytenberg \cite{Roytenberged.Contemp.Math.Vol.315Amer.Math.Soc.ProvidenceRI2002}]
\label{thmRoytenberg}
\hfill
\begin{itemize}
\item $\NP$-manifolds of degree 2 are in bijective correspondence with pseudo-Euclidean vector bundles. 
\item $\NPQ$-manifolds of degree 2 are in bijective correspondence with Courant algebroids. 
\end{itemize}
\end{Theorem}

The Poisson algebra of functions associated with a given $\NP$-manifold $\V$ of degree 2 \big(or equivalently the $2$-Schouten algebra $\Tpolyd{(2)}=\fonc{\V}[2]$\big) was referred to as the {\it Rothstein algebra}  in \cite{Keller2008}. The latter can be interpreted as the deformation complex of Hamiltonian functions on $\V$.
Via the second point of Theorem \ref{thmRoytenberg}, this can be rephrased as saying that the Rothstein algebra controls the deformation theory of Courant--Dorfman structures $(\mD,\brdot_E)$ on the pseudo-Euclidean vector bundle $(E,\dldot)$ -- where $E$ is defined by the fibration \eqref{fibration} as $\M\leftarrow E[1]\leftarrow\V$ -- according to the following sequence of bijective correspondences:
\begin{eqnarray}
(\mD,\brdot_E)\Leftrightarrow(\rho_a{}^\mu,T_{abc})\Leftrightarrow\cH=\rho_a{}^\mu\, \xi^a p_\mu+\frac16\, T_{abc}\, \xi^a\xi^b\xi^c\label{bijcor}
\end{eqnarray}
where the right-hand side makes use of the local set of coordinates $\Big\{\underset{0}{x^\mu},\underset{1}{\xi^a},\underset{2}{p_\mu}\Big\}$ [Example \ref{exasymp}].

The supergeometric interpretation of Courant algebroids provided by Theorem \ref{thmRoytenberg} will allow us to apply the results of Section \ref{section:Stable deformations} in order to generate new stable deformation formulas for Courant--Dorfman structures $(\mD,\brdot_E)$ on a given pseudo-Euclidean vector bundle $(E,\dldot)$.\footnote{In this sense, the procedure does not deform the {\it full} Courant algebroid structure since the bilinear form $\dldot$ remains undeformed. } As noted in Lemma \ref{lemmacoho}, the zeroth cohomology of the connected part of the full Kontsevich graph complex in $d=3$ is one dimensional and spanned by the triangle class  \ie 
$H^{0}(\fGCconn_3)=\corps\dl L_3\dr$, \cf Figure \ref{figloop}. This result ensures \big(\cf Proposition \ref{propapp}\big) that there exists a {\it unique} stable deformation of Courant algebroids that we now explicitly characterise. 
 
 Letting $(E,\dldot)$ be a pseudo-Euclidean vector bundle, we use the bijective correspondence \eqref{bijcor} in order to associate to each Courant--Dorfman structure $(\mD,\brdot_E)$ on $(E,\dldot)$ the corresponding Hamiltonian function $\cH=\rho_a{}^\mu\, \xi^a p_\mu+\frac16\, T_{abc}\, \xi^a\xi^b\xi^c$ with associated homological vector field $\Q:=\pb{\cH}{\cdot}_\om$. Via Corollary \ref{propmapcoHd=0}, the stable Hamiltonian flow associated with $L_3$ is defined as:
  \begin{eqnarray}
  \nn
\dot\cH=\Rep^{(3)}_3\big(L_3\big)(\cH^{\otimes3})=\mu_3\big(\Delta_{12}\, \Delta_{23}\, \Delta_{31}\, (\cH^{\otimes 3})\big)
\end{eqnarray}
where the expression of the operator $\Rep_3^{(3)}$ follows \eqref{Orientationmorphism2} with $\Delta$ given by [\cf eq.\eqref{defDeltadodd2}]:
\begin{eqnarray}
\nn
\Delta_{ij}=\frac{\p }{\p x^\mu_{(i)}}\, \frac{\p }{\p p_\mu^{(j)}}-\frac{\p }{\p p_\mu^{(i)}}\, \frac{\p }{\p x^\mu_{(j)}}+\frac{\p }{\p \xi^a_{(i)}}\, \kappa^{ab}\, \frac{\p }{\p \xi^b_{(j)}}.
\end{eqnarray}
Explicitly, the triangle Hamiltonian flow maps any Hamiltonian function $\cH$ towards the associated Rothstein cocycle $\dot\cH\in H^3(\fonc{\V}|\Q)$ defined as $\dot\cH=\dot\rho_a{}^\mu\, \xi^a p_\mu+\frac16\, \dot T_{abc}\, \xi^a\xi^b\xi^c$ where:
\begin{align*}
\dot\rho_a{}^\mu &=
\Rep^{(3)}_3\big(L_3 \big)(\cH^{\otimes3})_{a}{}^\mu \quad = \quad
\raisebox{-3ex}{\hbox{\begin{tikzpicture}[scale=0.5, >=stealth']
\node  (b1) at (0,0) {$\rho_\bullet\hspace{-1pt}$};
\node  (b2) at (2,0) {$\rho_\bullet{}^\mu$};
\node  (b3) at (1,-1.73) {$\rho_a$};
\draw[black,>=latex]  (b1) to (b2);
\draw[black,->,>=latex]  (b1) to (b3);
\draw[black,->,>=latex]  (b3) to (b2);
\end{tikzpicture}}}
+
\raisebox{-3ex}{\hbox{\begin{tikzpicture}[scale=0.5, >=stealth']
\node  (b1) at (0,0) {$\rho_\bullet\hspace{-3pt}$};
\node  (b2) at (2,0) {$\rho_\bullet{}^\mu\hspace{-5pt}$};
\node  (b3) at (1,-1.73) {$T_{a\bullet\bullet}$};
\draw[black,->,>=latex]  (b1) to (b2);
\draw[black,>=latex]  (b2) to (b3);
\draw[black,>=latex]  (b3) to (b1);
\end{tikzpicture}}} \\
\dot T_{abc} &=
6 \,\, \Rep^{(3)}_3\big(L_3 \big)(\cH^{\otimes3})_{abc} \quad = \quad
\raisebox{-3ex}{\hbox{\begin{tikzpicture}[scale=0.5, >=stealth']
\node  (b1) at (0,0) {$\rho_a\hspace{-1pt}$};
\node  (b2) at (2,0) {$\rho_b$};
\node  (b3) at (1,-1.73) {$\rho_c$};
\draw[black,->,>=latex]  (b1) to (b2);
\draw[black,->,>=latex]  (b2) to (b3);
\draw[black,->,>=latex]  (b3) to (b1);
\end{tikzpicture}}}
-
\raisebox{-3ex}{\hbox{\begin{tikzpicture}[scale=0.5, >=stealth']
\node  (b1) at (0,0) {$\rho_a$};
\node  (b2) at (2,0) {$\rho_b$};
\node  (b3) at (1,-1.73) {$\rho_c$};
\draw[black,->,>=latex]  (b2) to (b1);
\draw[black,->,>=latex]  (b3) to (b2);
\draw[black,->,>=latex]  (b1) to (b3);
\end{tikzpicture}}}
-
\raisebox{-3ex}{\hbox{\begin{tikzpicture}[scale=0.5, >=stealth']
\node  (b1) at (0,0) {$T_{a\bullet\bullet}\hspace{-5pt}$};
\node  (b2) at (2,0) {$T_{b\bullet\bullet}\hspace{-5pt}$};
\node  (b3) at (1,-1.73) {$T_{c\bullet\bullet}$};
\draw[black,>=latex]  (b2) to (b1);
\draw[black,>=latex]  (b3) to (b2);
\draw[black,>=latex]  (b1) to (b3);
\end{tikzpicture}}} \\
& \quad \quad \quad \quad \quad + \quad 3
\raisebox{-3ex}{\hbox{\begin{tikzpicture}[scale=0.5, >=stealth']
\node  (b1) at (0,0) {$\rho_\bullet$};
\node  (b2) at (2,0) {$\rho_a$};
\node  (b3) at (1,-1.73) {$T_{bc\bullet}$};
\draw[black,->,>=latex]  (b1) to (b2);
\draw[black,->,>=latex]  (b2) to (b3);
\draw[black,>=latex]  (b1) to (b3);
\end{tikzpicture}}}
+ \quad 3
\raisebox{-3ex}{\hbox{\begin{tikzpicture}[scale=0.5, >=stealth']
\node  (b1) at (0,0) {$\rho_\bullet$};
\node  (b2) at (2,0) {$T_{a\bullet\bullet}$};
\node  (b3) at (1,-1.73) {$T_{bc\bullet}$};
\draw[black,>=latex]  (b2) to (b1);
\draw[black,>=latex]  (b3) to (b2);
\draw[black,->,>=latex]  (b1) to (b3);
\end{tikzpicture}}}
+ \quad \text{skewsym. } (a-b-c)
\end{align*}
Here, the directed arrows stand for space-time derivatives while undirected arrows represent contractions of fiber indices with the non-degenerate symmetric bilinear form $\kappa$. The local expression of the Hamiltonian flow induced by the triangle cocycle can be equivalently expressed in components as:
\begin{align}
\dot\rho_a{}^\mu &=
\rho_b{}^\lambda\, \p_\lambda \rho_a{}^\nu\, \p_\nu\rho^{b|\mu}+\rho_b{}^\lambda\, \p_\lambda\rho_c{}^\mu\, T_{a}{}^{bc} \label{HamiltonianflowCA1} \\
\dot T_{abc} &=
\partial_{\mu}{\rho_a{}^{\nu}}\, \partial_{\nu}{\rho_b{}^{\lambda}}\,  \partial_{\lambda}{\rho_c{}^{\mu}} -\partial_{\mu}{\rho_a{}^{\lambda}}\,  \partial_{\nu}{\rho_b{}^{\mu}}\,  \partial_{\lambda}{\rho_c{}^{\nu}} -T_{a}{}^{d e} T_{b d f} T_{c e}{}^{ f} \nn \\
& \quad +3\, \rho_d{}^{\mu} \partial_{\mu}{\rho_{[a}{}^{\nu}} \partial_{\nu}{T_{b c]}{}^{ d}} +3\, \rho_{d}{}^\mu T_{[a}{}^{d e} \partial_{\mu}{T_{b c] e}} \nn
\end{align}
where indices are raised and lowered with $\kappa$. Consistently with Proposition \ref{propmapcoH}, it can be checked that $\pb{\cH}{\dot\cH}_\om=0$ modulo the relations \eqref{ConstraintCA1}-\eqref{ConstraintCA3} coming from $\pb{\cH}{\cH}_\om=0$. 

It should be emphasised that the situations corresponding to $d=2$ and $d=3$ are drastically different. In the case $d=2$, the zeroth cohomology is the ``dominant'' degree \ie contains an infinite number of non-trivial classes leading to infinitely many stable deformations of Poisson manifolds, \cf Section \ref{section:Poisson manifolds}. On the contrary, for $d=3$, the zeroth cohomology is one-dimensional and thus yields a unique stable deformation of Courant--Dorfman structures given by \eqref{HamiltonianflowCA1}. The ``dominant'' degree of $\fGC_3$ being $-3$ \big(\cf Section \ref{section:Cohomology of the full graph complex}\big), it would be desirable to find a construction mapping elements of $H^{-3}(\fGC_3)$ to Courant--Dorfman deformations. This is achieved through Weyl factors, as we now show.

\paragraph{Weyl factors for Courant algebroids.}
Recall from Section \ref{section:Cohomology of the full graph complex} that $H^{-3}(\fGC_3)$ is spanned by trivalent graphs modulo IHX relations [see Figure \ref{figtrivalent} and eq.\eqref{IHX}]. \Prop{propconfflow} ensures that each element $\gamma\in H^{-3}(\fGC_3)$ is mapped to a stable Weyl factor for Courant algebroids. Explicitly, given a trivalent graph $\gamma$, Hamiltonian functions $\cH$ can be mapped to Weyl factors $\Om_\gamma(\cH)\in H^0(\fonc{\V}|\Q)\simeq\Ker\mD$. The explicit local expression of the Weyl factor $\Om_\gamma(\cH)$ associated to a given graph $\gamma\in H^{-3}(\fGC_3)$ with $N$ vertices is given by $\Om_\gamma(\cH):=\Rep^{(3)}_N(\Gamgraph)(\cH^{\otimes N})$. We now exemplify this construction by displaying the Weyl factors associated to the simplest trivalent graphs. The simplest example of trivalent graph is given by the ``$\Theta$'' graph $\raisebox{-1.1ex}{\hbox{\begin{tikzpicture}[
my angle/.style={draw, <->, angle eccentricity=1.3, angle radius=9mm},scale=0.45
                        ]
\coordinate                     (O)  at (0,0);
\coordinate[]  (A) at (180:0.6);
\coordinate[]  (B) at (0:0.6);
\draw[thick]    (A) -- (B);
\draw[thick]  (0,0) circle (0.6cm);
    \end{tikzpicture}}}$ being the only connected trivalent graph with $N=2$ vertices. The latter yields the following Weyl factor:
\begin{equation}\label{Thetaflow}
\begin{aligned}
\Om_\Theta(\cH) := \Rep^{(3)}_2(\Theta)(\cH^{\otimes 2}) &=
\raisebox{-3ex}{\hbox{\begin{tikzpicture}[scale=0.5, >=stealth']
\tikzstyle{w}=[circle, draw, minimum size=4, inner sep=1]
\tikzstyle{b}=[circle, draw, fill, minimum size=4, inner sep=1]
\node  (b2) at (0,0) {$T_{\bullet\bullet\bullet}$};
\node (b3) at (3,0) {$T_{\bullet\bullet\bullet}$};
\draw[black][>=latex]  (b2) to[out=45, in=135, looseness=1.] (b3);
\draw[black][>=latex]  (b2) to (b3);
\draw[black][>=latex]  (b2) to[out=-45, in=-135, looseness=1.] (b3);
\end{tikzpicture}}}
+
6\, 
\raisebox{-3.5ex}{\hbox{\begin{tikzpicture}[scale=0.5, >=stealth']
\tikzstyle{w}=[circle, draw, minimum size=4, inner sep=1]
\tikzstyle{b}=[circle, draw, fill, minimum size=4, inner sep=1]
\node  (b2) at (0,0) {$\rho_\bullet$};
\node (b3) at (3,0) {$\rho_\bullet$};
\draw[black][>=latex,->]  (b2) to[out=55, in=130, looseness=1.] (b3);
\draw[black][>=latex]  (b2) to (b3);
\draw[black][>=latex,<-]  (b2) to[out=-50, in=-125, looseness=1.] (b3);
\end{tikzpicture}}} \\
&=
T_{a b c}\,  T^{a b c}+6\,  \partial_{\nu}{\rho_a{}^{\lambda}}\,  \partial_{\lambda}{\rho}^{a|\nu}
\end{aligned}
\end{equation}
and the equality $\mD\, \Om_\Theta=0$ follows from $\pb{\cH}{\cH}_\om=0$. In particular, this ensures that $\dot \cH_\Theta:=\Om_\Theta(\cH)\cdot \cH$ is a Rothstein cocycle.\footnote{It can be checked by brute-force computation that the vector space of universal Hamiltonian flows for $N=3$ is of dimension 2 and spanned by the triangle flow \eqref{HamiltonianflowCA1} and the Weyl $\Theta$-flow defined from \eqref{Thetaflow}.}

The next to simplest case is given by the graphs $A$ and $B$ from Figure \ref{figtrivalent} for $N=4$ yielding:
\begin{eqnarray}
\Om_A(\cH)&:=&\Rep^{(3)}_4\big(A\big)(\cH^{\otimes4})= T_{a b c}\, T^{a b d}\, T^{c e f}\, T_{d e f}+4\,  \partial_{\mu}{\rho_a{}^{\nu}}\,  \partial_{\nu}{\rho_b{}^{\mu}}\,  \partial_{\lambda}{\rho^{a|\rho}}\,  \partial_{\rho}{\rho^{b|\lambda}}\nn\\ 
&&\hspace{3cm}- 8\,  \partial_{\mu}{\rho_a{}^{\nu}}\, \partial_{\nu}{\rho^{a|\lambda}}\, \partial_{\lambda}{\rho_b{}^{\rho}}\, \partial_{\rho}{\rho^{b|\mu}}+4\,  \partial_{\mu}{\rho^{a|\nu}} \, \partial_{\nu}{\rho_d{}^{\mu}} T_{a b c}\,  T^{d b c}\nn\\
\Om_B(\cH)&:=&\Rep^{(3)}_4\big(B\big)(\cH^{\otimes4})\nn\\
&&\hspace{-0.65cm}=T_{a b c}\, T^a{}_{d e}\, T^{b d f}\, T^{c e}{}_{f} -8\,  \partial_{\mu}{\rho_a{}^{\nu}}\, \partial_{\nu}{\rho_b{}^{\lambda}}\, \partial_{\lambda}{\rho_c{}^{\mu}}\, T^{a b c} - 6\,  \partial_{\mu}{\rho_a{}^{\nu}}\, \partial_{\nu}{\rho_b{}^{\lambda}}\, \partial_{\lambda}{\rho^{a|\rho}}\, \partial_{\rho}{\rho^{b|\mu}}.\nn
\end{eqnarray}
Together with $\Om_\Theta(\cH)^2$ (corresponding to the disconnected graph $\gamma=\Theta\cup\Theta$), these are the only Weyl factors available for $N=4$. Note however that the trivalent graphs $A$ and $B$ can be related through the IHX relation \eqref{IHX} as $A\sim2\, B$, \cf footnote \ref{footIHX}. This ensures that their respective Weyl factors are related via $\Om_A(\cH)=2\, \Om_B(\cH)$ where the correspondence can be shown by making use of the constraints \eqref{ConstraintCA1}-\eqref{ConstraintCA3} coming from $\pb{\cH}{\cH}_\om=0$. 
\paragraph{Relation to quantization.}
Remarkably, Kontsevich's original quantization formula (\ie with standard propagator) can be interpreted  \cite{Kontsevich:1997vb,Cattaneo2000,CattaneoFelder2001a,CattaneoFelder2001} as a 3-point function in the path integral quantization of a 2-dimensional topological field theory
-- the \textbf{Poisson $\sigma$-model}, introduced in \cite{Ikeda1993,Ikeda1993a,Schaller1994} -- whose source is of dimension $d=2$ and whose target is the (shifted) cotangent bundle associated to the Poisson manifold. The graphs appearing in Kontsevich's stable formula can therefore be interpreted from the point of view of quantum field theory as Feynman diagrams associated with the quantification of the Poisson $\sigma$-model. 

As mentioned previously, the Poisson $\sigma$-model constitutes the first rung of an infinite ladder of AKSZ $\sigma$-models \cite{AlexandrovKontsevichSchwarzEtAl1997} associating to any symplectic Lie $n$-algebroid a topological field theory of dimension $d=n+1$. An interesting open problem concerns the possibility of generalising such interplay between deformation quantization results (on the algebraic side) and quantization of AKSZ-type of models (on the field theoretic side) to higher values of $n$. 

For $n=2$, the relevant AKSZ $\sigma$-model was constructed by D. Roytenberg in \cite{Roytenberg2007b} (\cf\cite{Ikeda2003} for an earlier derivation from consistent deformations of a Chern--Simons gauge theory coupled with a $0$-dimensional $BF$ theory). Such model associates to any Courant algebroid a canonical 3-dimensional topological field theory -- the \textbf{Courant $\sigma$-model}. From the field theory side, quantization of the Courant $\sigma$-model within the Batalin--Vilkovisky formalism \cite{Batalin1984} has been considered in \cite{HofmanParkCommun.Math.Phys.249:249-2712004,HofmanPark2007} (\cf also \cite{Ikeda2001} for a discussion of observables in general AKSZ $\sigma$-models).

On the algebraic side, a possible candidate for the quantum notion associated with Courant algebroids is given by \textbf{vertex algebroids}, as introduced in \cite{Gorbounov2000} from truncation of vertex algebras \cite{Gorbounov2000a}. Indeed, it was shown in \cite{Bressler2005} that the semi-classicalisation of (commutative) vertex algebroids yields a Courant algebroid. This suggests a formulation of a deformation quantization problem for Courant algebroids, similar to the one formulated in \cite{Berezin1975,BayenFlatoFronsdalEtAl1978} for Poisson manifolds. 

Although it is outside of the scope of the present paper to address the quantization problem for Courant algebroids, we note that some insights can be gained\footnote{We refer to \cite{Morand2021} for more details on the partition of deformation quantization problems according to graph cohomology. } from the classification of graph cocycles in $H^\bullet(\fGC_3)$:
\begin{enumerate}
\item $H^1(\fGC_3)=\mathbf{0}$: {\it The existence of stable formality morphisms for Courant algebroids is unobstructed.}
\item $H^0(\fGC_3)=\corps$: {\it The space of stable formality morphisms for Courant algebroids is of dimension $1$. }
\end{enumerate}

The first statement asserts that the $2$-Schouten algebra $(\Tpolyd{(2)},\brdot_{\mathsf S})$ is rigid -- at least in the ``stable setting'' -- \ie it does not admit non-trivial deformations as a \Lieinf-algebra. In other words, the graded Lie algebra $\Tpolyd{(2)}$ is intrinsically formal\footnote{A graded Lie algebra $(\alh,\brdot)$ is said to be \textbf{intrinsically formal} if any \Lieinf-algebra $(\alg,l)$ restricting to the dg Lie algebra $(\alh,0,\brdot)$ (\ie with vanishing differential) in cohomology \big(\ie $H(\alg,l)=(\alh,0,\brdot)$\big) is formal. } in the stable setting. The homotopy transfer theorem thus ensures that, given a dg Lie algebra (or \Lieinf-algebra) $\Dpolyd{(2)}$ such that $H^\bullet(\Dpolyd{(2)})$ is isomorphic to $\Tpolyd{(2)}$ as a graded Lie algebra and a quasi-isomorphism of complexes $\mathcal U_1:\Tpolyd{(2)}\iso\Dpolyd{(2)}$ given by stable formul\ae, the ``HKR-type''  map $\mathcal U_1$ can always be prolongated to a full stable \Lieinf quasi-isomorphism $\mathcal U:\Tpolyd{(2)}\iso\Dpolyd{(2)}$, which in turn would provide a quantization map for Courant algebroids.

Recall that such reasoning\footnote{We refer to Section 4.1 of \cite{Morand2021} for more details. } constituted the initial rationale behind the introduction of the graph complex $\fGC_2$ in \cite{Kontsevich1997}. However, in this case, it is a hard open conjecture that $H^1(\fGC_2)=\mathbf{0}$ so that M. Kontsevich had to rely on different methods in order to prove his formality theorem for Poisson manifolds. On the contrary, for Courant algebroids, it is straightforward to show the rigidity of the $2$-Schouten algebra (see above) so that one can use the original Kontsevich approach to prove a formality theorem for Courant algebroids. 

From the second statement, we learn that such morphism is not unique, but rather that formality morphisms form a $1$-dimensional space. Consequently, there should exist a one-parameter family of stable quantization maps for Courant algebroids\footnote{As in the Poisson case, some Courant algebroids are insensitive to deformations so that their associated quantum class is unique \ie independent of the parameter. This can be shown to be the case of {\it exact} Courant algebroids [see \eg \cite{Jurco2016} for a definition] which are insensitive to the triangular deformation defined by \eqref{HamiltonianflowCA1}. }. This is again in sharp contrast with the Poisson case for which $H^0(\fGC_2)$ is infinite-dimensional (being isomorphic to the Grothendieck--Teichm\"{u}ller Lie algebra $\grt_1$) and thus the space of formality morphisms (and consequently also the one of stable quantization maps for Poisson manifolds) forms an infinite-dimensional space (in bijective correspondence with the space of Drinfel'd associators). 

However, as shown in Corollary \ref{cor:WeylHamiltonian}, there is a way to consistently deform Courant algebroids via trivalent vertices (modulo IHX relations) using Weyl deformation maps. Composing such deformations with a given quantization map yields an infinite-dimensional space of quantizations for Courant algebroids. In other words, despite the fact that the space of stable formality morphisms is finite-dimensional (of dimension 1), the space of universal quantizations of a given Courant algebroid is infinite-dimensional, the arbitrariness being encoded into trivalent graphs (on top of the triangle graph).

\section*{Acknowledgements}
\noindent We are indebted to Thomas Basile for numerous stimulating discussions regarding various aspects of the present work as well as for valuable feedback on a preliminary version of the manuscript. We are also grateful to Vasily A. Dolgushev, Noriaki Ikeda and Hsuan--Yi Liao for useful exchanges and to Cl\'ement Berthiere for precious help in dealing with \LaTeX\,  issues. Finally, we would like to deeply thank the anonymous referee whose comments and suggestions greatly helped to improve the quality of the present paper. This work was supported by Brain Pool Program through the National Research Foundation of Korea (NRF) funded by the Ministry of Science and ICT (2018H1D3A1A01030137) and by Basic Science Research Program through the National Research Foundation of Korea (NRF) funded by the Ministry of Education (NRF-2020R1A6A1A03047877 and NRF-2022R1I1A1A01071497).

\bibliographystyle{plain}

\end{document}